\documentclass[11pt, twoside]{amsart}

\usepackage[utf8]{inputenc}
\usepackage[T1]{fontenc}
\usepackage{amssymb,amsmath,amstext}
\usepackage{hyperref, graphicx}

\newtheorem{theor}{Theorem}[section]
\newtheorem{lem}[theor]{Lemma}
\newtheorem{defin}[theor]{Definition}

\newtheorem{prop}[theor]{Proposition} 

\newtheorem{notation}[theor]{Notation}
\newtheorem{exam}[theor]{Example}

\newtheorem{cor}[theor]{Corollary}
\newtheorem{rem}[theor]{Remark}

\newtheorem{assump}[theor]{Assumption}

\newtheorem{termin}[theor]{Terminology}

\numberwithin{equation}{section}

\newcommand{\mr}{\mathrm}

\newcommand{\mb}{\mathbf}

\newcommand{\es}{\emptyset}

\newcommand{\uhr}{\upharpoonright}

\newcommand{\nts}{\negthickspace}
\newcommand{\uhrc}{\nts \upharpoonright \nts}

\newcommand{\dist}{\mathrm{dist}}

\newcommand{\mcA}{\mathcal{A}}
\newcommand{\mcB}{\mathcal{B}}
\newcommand{\mcC}{\mathcal{C}}

\newcommand{\mcG}{\mathcal{G}}

\newcommand{\mcN}{\mathcal{N}}

\newcommand{\mcT}{\mathcal{T}}

\newcommand{\mbB}{\mathbf{B}}

\newcommand{\mbE}{\mathbf{E}}

\newcommand{\mbT}{\mathbf{T}}
\newcommand{\mbU}{\mathbf{U}}

\newcommand{\mbW}{\mathbf{W}}
\newcommand{\mbX}{\mathbf{X}}
\newcommand{\mbY}{\mathbf{Y}}
\newcommand{\mbZ}{\mathbf{Z}}

\newcommand{\mbbE}{\mathbb{E}}
\newcommand{\mbbG}{\mathbb{G}}
\newcommand{\mbbP}{\mathbb{P}}
\newcommand{\mbbN}{\mathbb{N}}
\newcommand{\mbbR}{\mathbb{R}}

\newcommand{\msfC}{\mathsf{C}}
\newcommand{\msfN}{\mathsf{N}}

\newcommand{\rng}{\mathrm{rng}}

\title[Random expansions]
{Random expansions of finite structures \\ with bounded degree}

\author{Vera Koponen}

\address{Vera Koponen, Department of Mathematics, Uppsala University, Sweden.}

\email{vera.koponen@math.uu.se}

\date{22 October, 2025}

\begin{document}

\maketitle

\begin{abstract}
We consider finite relational signatures $\tau \subseteq \sigma$, a sequence of finite
base $\tau$-structures $(\mathcal{B}_n : n \in \mbbN)$ the cardinalities of which tend to infinity and such
that, for some number $\Delta$, the degree of (the Gaifman graph of) every $\mathcal{B}_n$ is at most $\Delta$.
We let $\mathbf{W}_n$ be the set of all expansions of $\mathcal{B}_n$ to $\sigma$ and we consider a 
probabilistic graphical model, a concept used in machine learning and artificial intelligence,
to generate a probability distribution $\mbbP_n$ on $\mathbf{W}_n$ for all $n$.
We use a many-valued ``probability logic'' with truth values in the unit interval to express probabilities within probabilistic 
graphical models and to express queries on $\mathbf{W}_n$.
This logic uses aggregation functions (e.g. the average) instead of quantifiers and it can express all queries 
(on finite structures) that
can be expressed with first-order logic since the aggregation functions maximum and minimum can be
used to express existential and universal quantifications, respectively.
The main results concern asymptotic elimination of aggregation functions
(the analogue of almost sure elimination of quantifiers for two-valued logics with quantifiers)
and the asymptotic distribution of truth values of formulas, the analogue of logical convergence results 
for two-valued logics. The structure theory that is developed for sequences $(\mathcal{B}_n : n \in \mbbN)$
as above may be of independent interest.
\end{abstract}

\tableofcontents

\section{Introduction}

\subsection*{Logical convergence laws}

Since the pioneering work of Glebskii, Kogan, Liogonkii and Talanov \cite{Gle} in 1969 and, independently, Fagin \cite{Fag},
logical convergence laws and the related notion of almost sure (or asymptotic) elimination of quantifiers, or aggregation functions,
have been studied in various contexts, that is, for various types of structures, various logics, and various probability distributions;
see for example \cite{ABFN, Bur, HK, Hill, Kai, KL, KPR, KoVa, Lyn, McC, MT, Zhu, SS} which is a far from complete list,
but it gives an idea of the variety of the results. 
Such results have implications of a more practical nature.
If a convergence law holds in a given context, then it is possible to use random sampling of structures with a 
fixed sufficiently large domain 
to estimate the probability that a sentence (of the logic) is true in a random structure with any large enough domain.
If one can prove a result about, say, almost sure elimination of quantifiers, then this typically implies a convergence law,
and the proof of almost sure elimination of quantifiers describes (in all cases that I am aware of) 
how to eliminate quantifiers, step by step, until one gets a quantifier-free formula that is ``almost surely''
equivalent to the initial formula.
The elimination procedure does not depend on the domain size, 
but only on the initial formula and the formalism that defines the probability distribution. 
The probability that the quantifier-free formula that is eventually produced holds for any choice of parameters (elements
from a domain) can now be computed (independently of the domain size) by only using the formula and the formalism that defines
the probability distribution.

\subsection*{Base structures with bounded degree}

Most convergence laws (e.g. those mentioned above) consider a context in which all relations are uncertain, that is,
described by a probabilistic model. 
But an agent may operate in a context in which some relations, or properties, are certain while other are uncertain.
For example, in many contexts there are certain cardinality constraints. 
For example, in a country, no matter how populous, there is only one president, or prime minister.
In an oligopolic market, there may be, say, at most 10 (holding) companies which dominate the market.
The various probabilistic models that have been considered when studying logical convergence laws are not suitable
for modelling such cardinality constraints. 
Rather than trying to describe this situation by a
probabilistic model we can simply consider random expansions
(by ``uncertain relations'') of a ``base structure'' with unary relation symbols $P$ and $Q$ where
the intrepretation of $P$ is a singleton set and the interpretation of $Q$ has cardinality at most 10.

As another example, consider a genealogical tree. If the tree represents many generations then the tree has many vertices,
but in reality there is a fixed finite bound on the number of children that any vertex has.
The kind of probability distributions that have been considered in finite model theory or in
Statistical Relational Artificial Intelligence \cite{BKNP, DKNP} are not suitable for modelling such cardinality constraints.
Instead we can view the genealogical tree as a base structure. 
We can now consider uncertain relations the probabilities of which depend on some probabilistic model and
on the underlying genealogical tree.
So the probability that some person $x$ has a particular medical condition can 
depend on both the relatives of $x$ (which are certain relationships) and medical conditions of the relatives that may with
some probability transfer to $x$ (uncertain relationships).

As a final example, the probability distributions typically considered when studying random graphs
are unsuitable for describing a road network
as some sort of random graph,
partly because in practice there is some fixed finite bound (a cardinality constraint) 
on the number of choices that one can make in a crossing,
and partly because it may be more natural to view a road network as fixed (in the shorter run) rather than as 
changing randomly. 
Instead the road network can be viewed as a base structure and various uncertain properties,
say the level of congestion on a piece of a road can be assigned probabilities, conditioned, say, on how the
road network looks locally around this piece of road, and possibly also on other properties that are modelled probabilistically.

Strictly mathematial examples, inspired by the above intuitive examples, of base structures
for which the main
results of this article apply are given in 
Section~\ref{Examples of sequences of base structures};
these examples include ``paths'', ``grids'', and Galton-Watson trees.
In the sequel $\tau$ will be a finite relational signature and $\mbB = (\mcB_n : n \in \mbbN^+)$
will be a sequence of finite $\tau$-``base structures'' $\mcB_n$ such that the (Gaifman) degree 
(Definition~\ref{definition of degree}) of each $\mcB_n$ is bounded by a
fixed $\Delta \in \mbbN$ (not depending on $n$), and the cardinality of the domain of $\mcB_n$ tends to infinity as $n$
tends to infinity.
Then we consider a finite relational signature $\sigma$ such that $\tau \subseteq \sigma$ and let 
$\mbW_n$ be the set of all $\sigma$-structures that are expansions of $\mcB_n$.

\subsection*{Logic and probability distributions}

Aggregation functions such as the average (of a finite sequence of reals) are useful for analyzing data
but usually do not return the values 0 or 1 as output.
Since this work is partially motivated by the aspiration to combine logical and probabilistic methods in
the fields of data mining, machine learning and artificial intelligence (AI), we will consider a many valued
logic, which we call $PLA^*$
(Definition~\ref{syntax of PLA*}), with (truth) values in the unit inverval $[0, 1]$ 
and which employs aggregation functions instead of quantifiers.
Since the first-order quantifiers can be expressed by using the aggregation functions maximum and minimum it follows that
$PLA^*$ subsumes first-order logic.

Motivated by the field of Statistical Relational Artifical intelligence (SRAI)
\cite{BKNP, DKNP}, which combines the logical and probabilistic approaches
to AI, we will use a kind of (para\-metrized/\-lifted) Bayesian network, or probabilistic graphical model (PGM)
\cite{BKNP, KMG, Koller}, to determine a 
probability distribution $\mbbP_n$ on $\mbW_n$ 
(intuitively speaking, the set of ``possible worlds that are based on $\mcB_n$'') for all $n$.
The kind of PGM that we use will be called a {\em $PLA^*(\sigma)$-network (based on $\tau$)} 
(Definition~\ref{definition of PLA-network})
because it consists of
a directed acyclic graph with vertex set $\sigma \setminus \tau$ and, for each vertex $R \in \sigma \setminus \tau$, 
a $PLA^*$ formula $\theta_R$
which expresses the probability of $R$ conditioned on its parents.
With a $PLA^*(\sigma)$-network one can model both dependencies and independencies between atomic relations 
(from a logician's point of view) or between 0/1-valued random variables (from a probabilist's point of view).
One can model probabilities which do not depend on $n$ and probabilities which depend on $n$.
For example, all probabilities of an edge relation considered by Shelah and Spencer in \cite{SS} can be modelled by
a $PLA^*(\sigma)$-network where $\sigma \setminus \tau$ contains a relation symbol of arity 2.

\subsection*{Challenges and results, informally}

Let $\varphi(x_1, \ldots, x_k)$ be a formula of $PLA^*$
with free variables $x_1, \ldots, x_k$ and which uses only relation symbols from $\sigma$.
For a finite $\sigma$-structure $\mcA$ (in the sense of first-order logic) and elements $a_1, \ldots, a_k$ from the domain of $\mcA$
we let $\mcA\big(\varphi(a_1, \ldots, a_k)\big)$ denote the value of $\varphi(a_1, \ldots, a_k)$ in $\mcA$.
Let $\mbbG$ be a $PLA^*(\sigma)$-network and, for each $n \in \mbbN^+$, let $\mbbP_n$ be the probability distribution 
on $\mbW_n$ which is induced by $\mbbG$.
For an interval $I \subset [0, 1]$ and $a_1, \ldots, a_k$ from the domain of $\mcB_n$ we can now ask what the probability
is that $\mcA\big(\varphi(a_1, \ldots, a_k)\big) \in I$  for a random $\mcA \in \mbW_n$.
In principle we can compute this probability by, for each $\mcA \in \mbW_n$, 
computing $c_\mcA := \mcA\big(\varphi(a_1, \ldots, a_k)\big)$ and $\alpha_\mcA := \mbbP_n(\mcA)$, and then adding all 
$\alpha_\mcA$ for which $c_\mcA \in I$, but the time needed is in general exponential in the cardinality of the domain of $\mcB_n$
which is assumed to tend to infinity as $n\to\infty$.

So we look for other methods to compute, or at least estimate, the above probability.
If the formula $\varphi(x_1, \ldots, x_k)$ above is aggregation-free, i.e. does not use any aggregation function,
then $\mcA\big(\varphi(a_1, \ldots, a_k)\big)$ can be computed without inspecting any other elements in the domain than
$a_1, \ldots, a_k$, so the computation is independent of $n$.
If the formula $\varphi(x_1, \ldots, x_k)$ uses aggregation functions only in such a way that they can be
evaluated by only inspecting a {\em bounded part of the domain},
or a ``{\em neighbourhood  of $a_1, \ldots, a_k$}'', that is determined
by $a_1, \ldots, a_k$, then $\mcA\big(\varphi(a_1, \ldots, a_k)\big)$
can still be computed in time that is independent of $n$.
In both cases, the probability that $\mcA\big(\varphi(a_1, \ldots, a_k)\big) \in I$, for a random $\mcA \in \mbW_n$,
can be computed by only computing the probabilities that various subsequences of a bounded part of the domain
satisfy certain atomic relations that are determined by $\varphi(x_1, \ldots, x_k)$; 
the time needed depends only on the $PLA^*(\sigma)$-network $\mbbG$ and $\varphi(x_1, \ldots, x_k)$, and not on $n$.

This motivates to look for conditions on a sequence $\mbB = (\mcB_n : n \in \mbbN^+)$ of base structures with degree
bounded by some fixed $\Delta$, 
a $PLA^*(\sigma)$-network $\mbbG$, and a formula $\varphi(x_1, \ldots, x_k)$ (of $PLA^*$)
which imply that $\varphi$ can be ``reduced'' to a ``simpler'' formula $\psi(x_1, \ldots, x_k)$, such that
\begin{itemize}
\item[(a)] $\psi(x_1, \ldots, x_k)$ uses  aggregation functions only in such a way that they can be
evaluated by only inspecting a bounded part of the domain that is determined
by $a_1, \ldots, a_k$, and

\item[(b)] for every $\varepsilon > 0$, if $n$ is large enough, then, with probability at least $1 - \varepsilon$,
$\big|\mcA\big(\varphi(a_1, \ldots, a_k)\big) - \mcA\big(\psi(a_1, \ldots, a_k)\big)\big| \leq \varepsilon$.
(We will say that $\varphi$ and $\psi$ are {\em asymptotically equivalent}.)
\end{itemize}
If the reduction of $\varphi$ to $\psi$ with properties (a) and (b) depends only on $\varphi$ and $\mbbG$ then, 
for every $\varepsilon > 0$ and all sufficiently large $n$,
we can estimate the probability that $\mcA\big(\varphi(a_1, \ldots, a_k)\big) \in I$ for a random $\mcA \in \mbW_n$
with error at most $\varepsilon$ and the time needed is independent from $n$.

The first condition that will be imposed on $\mbB$, Assumption~\ref{properties of the base structures},
is, roughly speaking, that for every $\tau$-structure $\mcC$
and all $k, \lambda \in \mbbN^+$,
either there is $m \in \mbbN$ such that, for all sufficiently large $n$, the number of elements from $\mcB_n$ 
which have a ``$\lambda$-neighbourhood'' that
is isomorphic to $\mcC$ is at most $m$,
or the number of such elements grows faster than every logarithm as $n \to \infty$.
Later, in Assumption~\ref{relative frequency of tau-closure types}, 
we also assume, intuitively speaking, that 
if $\mcC  \subset \mcC'$ are $\tau$-structures and $\lambda < \lambda'$, then the number
of elements of $\mcB_n$ that have a $\lambda'$-neighbourhood that is isomorphic to $\mcC'$ divided by the number of 
elements of $\mcB_n$ that have a $\lambda$-neighbourhood that is isomorphic to $\mcC$ converges as $n \to \infty$.
All examples of Section~\ref{Examples of sequences of base structures} satisfy both
 Assumption~\ref{properties of the base structures}
 and 
 Assumption~\ref{relative frequency of tau-closure types}.
 
 We will also impose conditions on the formulas used by a $PLA^*(\sigma)$-network, and identical conditions
 on a formula that expresses a query, denoted $\varphi$ above.
 One of these conditions is that the formulas use only continuous aggregation functions
 such as the average, or, under stronger assumptions
 (which hold if the probability of every atomic relation does not depend on $n$) admissible
 (or ``semicontinuous'') aggregation functions, for example maximum and minimum.
 Assumptions will also be made on so-called ``conditioning formulas'' which specify over which elements, or tuples, in the
 domain an aggregation function ranges when used in a formula (see Definition~\ref{syntax of PLA*}).
 Since probabilities of uncertain relations may depend on the underlying structure of the certain relations, 
 represented by the sequence of structures $\mbB$, 
 it seems like we must make some assumptions on the  conditioning formulas, or else we would have to
 impose stronger conditions on the sequence of base structures $\mbB$ or on the $PLA^*(\sigma)$-network $\mbbG$.
 
 The first main result, Theorem~\ref{main result about strongly unbounded aggregations},
 shows, roughly speaking, that if 
 Assumption~\ref{properties of the base structures} 
 holds and all conditioning formulas used by $\mbbG$ and by $\varphi(x_1, \ldots, x_k)$ 
 are essentially conjunctions of atomic $\tau$-formulas or negations of such
 (with the possible exception that bounded aggregations may be
 used to talk about ``rare elements'' if such exist),
 then $\varphi$ is asymptotically equivalent to a formula $\psi$ with only bounded aggregations,
so (a) and (b) above are satisfied,
 and $\varphi$ satisfies a convergence law.
 Corollary~\ref{corollary to main result about strongly unbounded aggregations}
 tells that if, in addition, $\mbbG$ is such that probabilities of atomic relations do not depend on $n$ 
 and $\varphi(x_1, \ldots, x_k)$ uses only admissible aggregation functions (which are more general than the continuous ones and
 include maximum and minimum), then we get the same conclusions as in
 Theorem~\ref{main result about strongly unbounded aggregations}.
Theorem~\ref{main result about aggregations in general with extra assumption} 
has a similar statement as Theorem~\ref{main result about strongly unbounded aggregations}.
But since Theorem~\ref{main result about aggregations in general with extra assumption} 
assumes that both Assumption~\ref{properties of the base structures} 
and Assumption~\ref{relative frequency of tau-closure types} (about $\mbB$)  hold it
applies to {\em more} $PLA^*$-networks and {\em more} formulas $\varphi(x_1, \ldots, x_k)$, in particular, it allows more
general conditioning formulas,
which talk about {\em larger neighbourhoods},
to be used than in Theorem~\ref{main result about strongly unbounded aggregations}.
Corollary~\ref{corollary to main result about aggregations in general with extra assumption}
shows that if, in addition, probabilities of atomic relations do not depend on $n$ 
then we can allow $\varphi$ to contain admissible aggregation functions and we still get the same
conclusions as in Theorem~\ref{main result about aggregations in general with extra assumption}.
Examples~\ref{example for corollary to theorem with no extra assumption}
and~\ref{example of results with additional assumption}
describe some concrete contexts in which the main results can be applied.

Besides the main results, the results in Sections~\ref{the base sequence}
and~\ref{properties of closure types}
may be of some independent interest as they develop a theory of the asymptotics of sequences of finite structures
with bounded degree.

\subsection*{Related work}

In \cite{KT} Koponen and Tousinejad considers the same general situation as described above,
{\em but} in \cite{KT} it is assumed that, for a fixed $\Delta$, each $\mcB_n$ is a tree
the height of which is bounded by $\Delta$, but there is no bound on the number of children that a vertex may have, 
so the degree of $\mcB_n$ is not bounded. The overall approach to asymptotically eliminating aggregation functions
(as described in Section~\ref{A general approach} below) is the
same in \cite{KT} as here, {\em but} the technical work done in 
Sections \ref{convergence and balance}, \ref{Proving convergence in the inductive step}, 
and \ref{Finding the balance in the inductive step}
is completely different from the corresponding work in \cite{KT}, since we consider base structures of a different kind here.
Also, in this study it is necessary to do much more preparatory work,
in Sections~\ref{the base sequence} and~\ref{properties of closure types},
about the base structures $\mcB_n$.
In \cite{KK} Koponen and Karlsson study expansions of linear {\em pre}orders and convergence laws
for first-order logic and its expansion with proportion quantifiers.

It seems like Lynch  was first to study random expansions of nontrivial structures
and he formulated a condition ($k$-extendibility) on the sequence of base structures $\mbB$ that guarantees that a
convergence law holds for first-order logic and the uniform probability distribution;
he also applied his result to get convergence laws for some 
sequences $\mbB$ of base structures \cite{Lyn}.
Shelah \cite{She02} and Baldwin  \cite{Bal03} considered a context that can be expressed by using
a signature $\tau = \{E\}$, where $E$ is a binary relation symbol, and base structures
$\mbB = (\mcB_n : n \in \mbbN^+)$ where each $\mcB_n$ is a directed graph isomorphic to a directed path of length $n$.
Let $R$ be another binary relation symbol.
Shelah proved a first-order convergence law for random expansions to $\sigma = \{E, R\}$ of the structures in $\mbB$ where 
the probability of an (undirected) $R$-edge between two vertices equals
the distance (in the underlying $E$-path $\mcB_n$) between the vertices raised to $-\alpha$ where $\alpha \in (0, 1)$ is irrational.
Baldwin \cite{Bal03} proved a first-order zero-one law for random expansions to $\sigma = \{E, R\}$ of 
the structures in $\mbB$ where the probability of an (undirected) $R$-edge is $n^{-\alpha}$ for an irrational 
$\alpha \in (0, 1)$ where $n$ is the length of the underlying directed path (represented by $\mcB_n$).
Lynch \cite[Corollary~2.16]{Lyn} and later
Abu Zaid, Dawar, Grädel and Pakusa \cite{ADGP} and Dawar, Grädel and Hoelzel \cite{DGH} 
proved first-order (and $L_{\infty, \omega}^\omega$) convergence laws for expansions of $\mcB_n$
(and the uniform distribution) where 
$\mcB_n$ is a product of finite cyclic groups. In \cite{DGH} it is shown that if $\mcB_n$ is an $n$-fold product of
a linear order with exactly two elements, then a first-order convergence law for expansions of $\mcB_n$ 
(and the uniform distribution) does {\em not} hold.

Studies of convergence laws where the probability distribution or logic is inspired by concepts from machine learning and AI,
but which do not consider underlying base structures (so all relations are modelled probabilistically), include
\cite{Jae98a, Jae98b, CM, Kop20, Wei21, KW1, KW2, Adam-Day1, Adam-Day2, Wei24, KT, Kop25}, 
in chronologial order.

\subsection*{Structure of the article}

Section~\ref{Preliminaries} clarifies basic notation and terminology that will be used and 
recalls a couple of probability theoretic results.
Section~\ref{Probability logic with aggregation functions}
defines the notion of aggregation function, the syntax and semantics of $PLA^*$ and some related notions.

Section~\ref{A general approach} defines
the notion of {\em asymptic equivalence of formulas}
(Definition~\ref{definition of asymptotically equivalent formulas})
and describes a general approach, from \cite{KW3},
of asymptotic elimination of aggregation functions. 
The idea is that if we can find a ``basic sublogic'' $L_0$ of $PLA^*$ that contains only ``well behaved and simple'' formulas 
in the sense of Assumption~\ref{assumptions on the basic logic},
then every {\em continuous} aggregation function 
(see Definition~\ref{definition of strong admissibility}) 
can be ``asymptotically eliminated''.
Theorem~\ref{general asymptotic elimination} (a consequence of \cite{KW3}) formulates this conclusion.
By repeating such ``asymptotic eliminations'' of aggregation functions one finally gets
an``$L_0$-basic formula'' which is essentially a boolean combination of formulas in 
$L_0$ and which is asymptotically equivalent to the original formula.
Thus the approach later in the article will be to find a set $L_0$ of ``simple'', or ``basic'', formulas that satisfy
the conditions of Assumption~\ref{assumptions on the basic logic}.

In Section~\ref{the base sequence} we state the exact assumptions that we make on the sequence of base structures
$\mbB = (\mcB_n : n \in \mbbN^+)$
(Assumption~\ref{properties of the base structures}). 
Besides the assumptions that for some $\Delta \in \mbbN$ all $\mcB_n$ have degree at most $\Delta$ we need to 
assume that the base structures $\mcB_n$ behave in a sufficiently uniform way. 
The assumptions on $\mbB$ involves the notion of {\em (un)bounded neighbourhood type}, and later we also
define a notion of {\em (un)bounded closure type}.
Then some results about neighbourhood and closure types are proved (which will be essential later).
Section~\ref{Examples of sequences of base structures}
gives examples of sequences of base structures that satisfy the conditions of 
Assumption~\ref{properties of the base structures}.
Then Section~\ref{properties of closure types}
proceeds with a more detailed study of bounded, unbounded, uniformly (and strongly) unbounded closure types.
The results will be used later.
Recall that all base structures $\mcB_n$ are $\tau$-structures for some finite relational signature $\tau$.

In Section~\ref{Probability distributions}
we consider a larger finite relational signature $\sigma \supset \tau$ and let $\mbW_n$ be the set of all
expansions to $\sigma$ of $\mcB_n$.
We define the notion of {\em $PLA^*$-network} and explain how it induces a probability distribution $\mbbP_n$ on $\mbW_n$
for every $n$. 
Example~\ref{examples of PLA-networks} gives an idea of the expressivity of $PLA^*$-networks in terms of what kind
of distributions they can induce on $\mbW_n$.
The ``global structure'' of proof of the main results 
(stated in Section~\ref{Asymptotic elimination of aggregation functions}) 
is an
induction on the ``maximal path rank'' of the underlying DAG of the $PLA^*$-network that induces the distributions $\mbbP_n$.
The maximal path rank of a DAG is the length of the longest directed path in the DAG,
and we have the convention that a DAG with empty vertex set has maximal path rank $-1$.
(A DAG with at least one vertex but no edges has maximal path rank 0.)
The bulk of the proof consists of proving, by induction on the maximal path rank, 
that a ``{\em convergence}'' condition and a {``\em balance}'' condition
holds for closure types (defined in Section~\ref{the base sequence}), which are 0/1-valued formulas that express,
for some $\lambda \in \mbbN$, 
what the ``$\lambda$-closure'' around some elements looks like.

In Section~\ref{convergence and balance}
the notions of ``{\em convergent pairs of formulas}'' and ``{\em balanced triples of formulas}'' are defined.
The reason for considering pairs and triples of  formulas is that we need to study, on the one hand,
the probability that a formula is satisfied given that another formula is satisfied, and on the other hand,
the frequency of tuples that satisfy a formula among the set of tuples that satisfy another formula, under the assumption
that some ``passive'' parameters satisfy a third formula.
The relevance of balanced triples is that one part of 
Assumption~\ref{assumptions on the basic logic}, which is used for the asymptotic elimination of aggregation functions,
is  that certain triples are balanced.
The relevance of convergent pairs is that we use it for proving results about balanced triples.

The base case of the induction (when the maximal path rank is $-1$, or equivalently, when $\sigma = \tau$)
is taken care of in Section~\ref{convergence and balance}.
This section also formulates the induction hypothesis 
(Assumption~\ref{induction hypothesis}) that will be used in the inductive step.
In Sections~\ref{Proving convergence in the inductive step}
and~\ref{Finding the balance in the inductive step}
the inductive step of the proof is carried out;
the ``convergence property'' is taken care of in Section~\ref{Proving convergence in the inductive step}
and the ``balance property'' in Section~\ref{Finding the balance in the inductive step}.
Finally, in Section~\ref{Asymptotic elimination of aggregation functions}
we put together the pieces from earlier sections to get our main theorems.

\section{Preliminaries}\label{Preliminaries}

\noindent
For basics about first-order logic and finite model theory see for example \cite{EF}.
Structures in the sense of first-order logic are denoted by calligraphic letters $\mcA, \mcB, \mcC, \ldots$ and their
domains (universes) by the corresponding noncalligraphic letter $A, B, C, \ldots$.
Finite sequences (tuples) of objects are denoted by $\bar{a}, \bar{b}, \ldots, \bar{x}, \bar{y}, \ldots$.
We usually denote logical variables by $x, y, z, u, v, w$. 
Unless stated otherwise, when $\bar{x}$ is a sequence of variables we assume that $\bar{x}$ does not repeat a variable.
But if $\bar{a}$ denotes a sequence of elements from the domain of a structure then repetitions may occur
(unless something else is said).

We let $\mbbN$ and $\mbbN^+$ denote the set of nonnegative integers and the set of positive integers, respectively.
For a set $S$, $|S|$ denotes its cardinality, and for a finite sequence $\bar{s}$, $|\bar{s}|$ denotes its length
and $\rng(\bar{s})$ denotes the set of elements in $\bar{s}$.
For a set $S$, $S^{<\omega}$ denotes the set of finite nonempty sequences (where repetitions are allowed) of elements from $S$,
so $S^{<\omega} = \bigcup_{n\in \mbbN^+}S^n$.
In particular, $[0, 1]^{<\omega}$ denotes the set of all finite nonempty sequences of reals from the
unit interval $[0, 1]$.

A signature (vocabulary) is called {\em finite (and) relational} if it is finite and contains only relation symbols.
Let $\sigma$ be a signature and let $\mcA$ be a $\sigma$-structure.
If $\tau \subset \sigma$ then $\mcA \uhrc \tau$ denotes the {\em reduct} of $\mcA$ to $\tau$.
If $B \subset A$ then $\mcA \uhrc B$ denotes the {\em substructure} of $\mcA$ generated by $B$.
We let $FO(\sigma)$ denote the set of first-order formulas that can be constructed with the signature $\sigma$.

Suppose that $\mcG$ is a {\em directed acyclic graph (DAG)}. 
As a convention we allows the vertex set of a DAG to be empty.
The {\em maximal path rank}, or just mp-rank, of $\mcG$, denoted $\mr{mp}(\mcG)$, 
is the length of the longest directed path in $\mcG$
if its vertex set is nonempty. 
If the vertex set is empty then we stipulate that the maximal path rank is $-1$.

A random variable will be called {\em binary} if it can only take the value $0$ or $1$.
The following is a direct consequence of \cite[Corollary~A.1.14]{AS} which in turn follows from the
Chernoff bound \cite{Che}:

\begin{lem}\label{independent bernoulli trials}
Let $Z$ be the sum of $n$ independent binary random variables, each one with probability $p$ of having the value 1.
For every $\varepsilon > 0$ there is $c_\varepsilon > 0$, depending only on $\varepsilon$, such that the probability that
$|Z - pn| > \varepsilon p n$ is less than $2 e^{-c_\varepsilon p n}$.
(If $p = 0$ the same statement holds if `$2e^{-c_\varepsilon p n}$' is replaced by `$e^{-n}$'.)
\end{lem}

\noindent
The following is a straightforward corollary (the full proof of which is given in \cite{KW2}):

\begin{cor}\label{independent bernoulli trials, second version}
Let $p \in [0, 1]$ and let $\varepsilon > 0$.
Let $Z$ be the sum of $n$ independent binary random variables $Z_1, \ldots, Z_n$, where for each $i = 1, \ldots, n$ the
probability that $Z_i$ equals 1 belongs to the interval $[p-\varepsilon, p+\varepsilon]$.
Then there is $c_\varepsilon > 0$, depending only on $\varepsilon$, such that the probability that
$Z > (1 + \varepsilon)(p + \varepsilon) n$ or 
$Z < (1 - \varepsilon)(p - \varepsilon) n$ 
is less than $2 e^{-c_\varepsilon p n}$.
\end{cor}

\noindent
The following lemma follows easily from the definition of conditional probability.

\begin{lem}\label{basic fact about conditional probabilities}
Suppose that $\mbbP$ is a probability measure on a set $\Omega$.
Let $X \subseteq \Omega$ and $Y \subseteq \Omega$ be measurable.
Also suppose that $Y = Y_1 \cup \ldots \cup Y_k$, $Y_i \cap Y_j = \es$ if $i \neq j$, and that each $Y_i$ is measurable.
If $\alpha \in [0, 1]$, $\varepsilon > 0$, and
$\mbbP(X \ | \ Y_i) \in [\alpha - \varepsilon, \alpha + \varepsilon]$ for all $i = 1, \ldots, k$,
then $\mbbP(X \ | \ Y) \in [\alpha - \varepsilon, \alpha + \varepsilon]$.
\end{lem}

\section{Probability logic with aggregation functions}\label{Probability logic with aggregation functions}

\noindent
In this section we define the syntax and semantics of $PLA^*$ ({\em probability logic with aggregation functions})
which is a logic with (truth) values in $[0, 1]$ that we consider both for expressing queries and for 
expressing probabilities. 
Example~\ref{example of page rank}
shows how it can express queries and
Example~\ref{examples of PLA-networks}
shows how it can express probabilities in the context of a $PLA^*$-network
(defined in Section~\ref{Probability distributions}).
A quite similar logic has been considered by Jaeger in for example \cite{Jae98a}. 
Variants of $PLA^*$, called $PLA$ and $PLA^+$, where considered 
by the author and Weitkämper in \cite{KW1} and \cite{KW2}, respectively.
$PLA^*$, first given in \cite{KW3}, differs from these by allowing a more flexible use of aggregation functions.
Recall that $[0, 1]^{<\omega}$ denotes the set of all finite nonempty sequences of reals in 
the unit interval $[0, 1]$.

\begin{defin}\label{definition of connective and aggregation function} {\rm
Let $k \in \mbbN^+$.\\
(i) A function $\msfC : [0, 1]^k \to [0, 1]$ will also be called a {\bf \em ($k$-ary) connective}.\\
(ii) A function $F : \big([0, 1]^{<\omega}\big)^k \to [0, 1]$ which is symmetric in the following sense
will be called a {\bf \em ($k$-ary) aggregation function}:
if $\bar{p}_1, \ldots, \bar{p}_k \in [0, 1]^{<\omega}$ and, for $i = 1, \ldots, k$,
$\bar{q}_i$ is a reordering of the entries in $\bar{p}_i$,
then $F(\bar{q}_1, \ldots, \bar{q}_k) = F(\bar{p}_1, \ldots, \bar{p}_k)$.
}\end{defin}

\noindent
The functions defined in the next definition are continuous and when restricted to $\{0, 1\}$
(as opposed to the interval $[0, 1]$) they have the usual meanings  of $\neg$, $\wedge$, $\vee$, and
$\rightarrow$, respectively. (The definitions correspond to the semantics of Lukasiewicz logic 
(see for example \cite[Section~11.2]{Ber}, or \cite{LT}).

\begin{defin}\label{special connectives} {\rm
Let
\begin{enumerate}
\item $\neg : [0, 1] \to [0, 1]$ be defined by $\neg(x) = 1 - x$,
\item $\wedge : [0, 1]^2 \to [0, 1]$ be defined by $\wedge(x, y) = \min(x, y)$,
\item $\vee : [0, 1]^2 \to [0, 1]$ be defined by $\vee(x, y) = \max(x, y)$, and
\item $\rightarrow : [0, 1]^2 \to [0, 1]$ be defined by $\rightarrow(x, y) = \min(1, \ 1 - x + y)$.
\end{enumerate}
}\end{defin}

\noindent
We now define some common aggregation functions.

\begin{defin}\label{examples of aggregation functions} {\rm
For $\bar{p} = (p_1, \ldots, p_n) \in [0, 1]^{<\omega}$, define
\begin{enumerate}
\item $\max(\bar{p})$ to be the {\em maximum} of all $p_i$,
\item $\min(\bar{p})$ to be the {\em minimum} of all $p_i$,
\item $\mr{am}(\bar{p}) = (p_1 + \ldots + p_n)/n$, so `am' is the {\em arithmetic mean}, or {\em average}, 
\item $\mr{gm}(\bar{p}) = \big(\prod_{i=1}^n p_i\big)^{(1/n)}$, so `gm' is the {\em geometric mean}, and
\item for every $\beta \in (0, 1]$, $\mr{length}^{-\beta}(\bar{p}) = |\bar{p}|^{-\beta}$, and
\end{enumerate}
}\end{defin}

\noindent
All aggregaton functions above are unary, that is, they take only one sequence as input.
But there are useful aggregation functions of higher arities, i.e. taking two or more sequences as input,
as shown in Examples~5.5 -- 5.7 in \cite{KW1} and in Example~6.4 in \cite{KW1}.

{\bf \em For the rest of this section we fix a finite and relational signature $\sigma$.}

\begin{defin}\label{syntax of PLA*}{\bf (Syntax of $PLA^*$)} {\rm 
We define formulas of $PLA^*(\sigma)$, as well as the set of free variables  of a formula $\varphi$, 
denoted $Fv(\varphi)$, as follows.
\begin{enumerate}
\item  For each $c \in [0, 1]$, $c \in PLA^*(\sigma)$ (i.e. $c$ is a formula) and $Fv(c) = \es$. 
We also let $\bot$ and $\top$
denote $0$ and $1$, respectively.

\item For all variables $x$ and $y$, `$x = y$' belongs to $PLA^*(\sigma)$ and $Fv(x = y) = \{x, y\}$.

\item For every $R \in \sigma$, say of arity $r$, and any choice of variables $x_1, \ldots, x_r$, $R(x_1, \ldots, x_r)$ belongs to 
$PLA^*(\sigma)$ and  $Fv(R(x_1, \ldots, x_r)) = \{x_1, \ldots, x_r\}$.

\item If $k \in \mbbN^+$, $\varphi_1, \ldots, \varphi_k \in PLA^*(\sigma)$ and
$\msfC : [0, 1]^k \to [0, 1]$ is a continuous connective, then 
$\msfC(\varphi_1, \ldots, \varphi_k)$ is a formula of $PLA^*(\sigma)$ and
its set of free variables is $Fv(\varphi_1) \cup \ldots \cup Fv(\varphi_k)$.

\item Suppose that $\varphi_1, \ldots, \varphi_k \in PLA^*(\sigma)$,
$\chi_1, \ldots, \chi_k \in PLA^*(\sigma)$,
$\bar{y}$ is a sequence of distinct variables,
and that $F : \big( [0, 1]^{<\omega} \big)^k \to [0, 1]$ is an aggregation function.
Then 
\[
F(\varphi_1, \ldots, \varphi_k : \bar{y} : 
\chi_1, \ldots, \chi_k)
\]
is a formula of $PLA^*(\sigma)$ and its set of free variables is
\[
\big( \bigcup_{i=1}^k Fv(\varphi_i)\big) \setminus \rng(\bar{y}),
\] 
so this construction binds the variables in $\bar{y}$.
The construction $F(\varphi_1, \ldots, \varphi_k : \bar{y} : \chi_1, \ldots, \chi_k)$ will be called an
{\em aggregation (over $\bar{y}$)} and the formulas $\chi_1, \ldots, \chi_k$ are called the {\em conditioning formulas} of this aggregation.
\end{enumerate}
}\end{defin}

\begin{defin}\label{definition of subformula} {\rm
(i) A formula in $PLA^*(\sigma)$ without free variables is called a {\bf \em sentence}.\\
(ii) In part~(4) of Definition~\ref{syntax of PLA*} the formulas $\varphi_1, \ldots, \varphi_k$
are called {\bf \em subformulas} of $\msfC(\varphi_1, \ldots, \varphi_k)$.\\
(iii) In part~(5) of Definition~\ref{syntax of PLA*} the formulas $\varphi_1, \ldots, \varphi_k$ and
$\chi_1, \ldots, \chi_k$ are called {\em subformulas} of
$F(\varphi_1, \ldots, \varphi_k : \bar{y} : \chi_1, \ldots, \chi_k)$.
We also call $\chi_1, \ldots, \chi_k$ {\bf \em conditioning subformulas} of 
$F(\varphi_1, \ldots, \varphi_k : \bar{y} : \chi_1, \ldots, \chi_k)$.\\
(iv) We stipulate the following transitivity properties:
If $\psi_1$ is a subformula of $\psi_2$ and $\psi_2$ is a subformula of $\psi_3$,
then $\psi_1$ is a subformula of $\psi_3$. 
If $\psi_1$ is a conditioning subformula of $\psi_2$ and $\psi_2$ is a subformula of $\psi_3$, then $\psi_1$
is a conditioning subformula of $\psi_3$.
}\end{defin}

\begin{notation}\label{abbreviation when using aggregation functions}{\rm
(i) When denoting a formula in $PLA^*(\sigma)$ by for example $\varphi(\bar{x})$ then we assume that $\bar{x}$ is a sequence
of different variables and that every free variable in the formula denoted by $\varphi(\bar{x})$ belongs to $\rng(\bar{x})$
(but we do not require that every variable in $\rng(\bar{x})$ actually occurs in the formula).\\
(ii) If all $\chi_1, \ldots, \chi_k$ are the same formula $\chi$,
then we may abbreviate 
\[
F(\varphi_1, \ldots, \varphi_k : \bar{y} : 
\chi_1, \ldots, \chi_k)
\qquad \text{ by } \qquad
F(\varphi_1, \ldots, \varphi_k : \bar{y} : \chi).
\]
}\end{notation}

\begin{defin}\label{definition of literal} {\rm
The $PLA^*(\sigma)$-formulas described in parts~(2) and~(3) of 
Definition~\ref{syntax of PLA*}
are called {\bf \em first-order atomic $\sigma$-formulas}.
A $PLA^*(\sigma)$-formula is called a {\bf \em first-order $\sigma$-literal} if it has the form $\varphi(\bar{x})$
or $\neg\varphi(\bar{x})$, where $\varphi(\bar{x})$ is a first-order atomic formula and
$\neg$ is like in Definition~\ref{special connectives}
(so it corresponds to negation when truth values are restricted to 0 and 1).
}\end{defin}

\begin{defin}\label{semantics of PLA*}{\bf (Semantics of $PLA^*$)} {\rm
For every $\varphi \in PLA^*(\sigma)$ and every sequence of distinct variables $\bar{x}$ such that 
$Fv(\varphi) \subseteq \rng(\bar{x})$ we associate a mapping from pairs $(\mcA, \bar{a})$,
where $\mcA$ is a finite $\sigma$-structure and $\bar{a} \in A^{|\bar{x}|}$, to $[0, 1]$.
The number in $[0, 1]$ to which $(\mcA,\bar{a})$ is mapped is denoted $\mcA(\varphi(\bar{a}))$
and is defined by induction on the complexity of formulas.
\begin{enumerate}
\item If $\varphi(\bar{x})$ is a constant $c$ from $[0, 1]$, then $\mcA(\varphi(\bar{a})) = c$.

\item If $\varphi(\bar{x})$ has the form $x_i = x_j$, then $\mcA(\varphi(\bar{a})) = 1$ if $a_i = a_j$,
and otherwise $\mcA(\varphi(\bar{a})) = 0$.

\item For every $R \in \sigma$, of arity $r$ say, if $\varphi(\bar{x})$ has the form $R(x_{i_1}, \ldots, x_{i_r})$,
then $\mcA(\varphi(\bar{a})) = 1$ if $\mcA \models R(a_{i_1}, \ldots, a_{i_r})$
(where `$\models$' has the usual meaning
of first-order logic), and otherwise $\mcA(\varphi(\bar{a})) = 0$.

\item If $\varphi(\bar{x})$ has the form $\msfC(\varphi_1(\bar{x}), \ldots, \varphi_k(\bar{x}))$,
where $\msfC : [0, 1]^k \to [0, 1]$ is a continuous connective, then
\[
\mcA\big(\varphi(\bar{a})\big) \ = \ 
\msfC\big(\mcA(\varphi_1(\bar{a})), \ldots, \mcA(\varphi_k(\bar{a}))\big).
\]

\item Suppose that $\varphi(\bar{x})$ has the form 
\[
F(\varphi_1(\bar{x}, \bar{y}), \ldots, \varphi_k(\bar{x}, \bar{y}) : \bar{y} : 
\chi_1(\bar{x}, \bar{y}), \ldots, \chi_k(\bar{x}, \bar{y}))
\] 
where $\bar{x}$ and $\bar{y}$ are sequences of distinct variables, $\rng(\bar{x}) \cap \rng(\bar{y}) = \es$, and
$F : \big( [0, 1]^{<\omega} \big)^k \to [0, 1]$ is an aggregation function.
If, for every $i = 1, \ldots, k$, the set 
$\{\bar{b} \in A^{|\bar{y}|} : \mcA(\chi_i(\bar{a}, \bar{b})) = 1\}$ is nonempty, then
let 
\[
\bar{p}_i = 
\big(\mcA\big(\varphi_i(\bar{a}, \bar{b})\big) : \bar{b} \in A^{|\bar{y}|} \text{ and } 
\mcA\big(\chi_i(\bar{a}, \bar{b})\big) = 1\big)
\]
and 
\[
\mcA\big(\varphi(\bar{a})\big) = F(\bar{p}_1, \ldots, \bar{p}_k).
\]
Otherwise let $\mcA\big(\varphi(\bar{a})\big) = 0$.
\end{enumerate}
}\end{defin}

\begin{defin}\label{L-basic subformula} {\bf (Special kinds of formulas)} {\rm 
(i) A formula in $PLA^*(\sigma)$ such that no aggregation function occurs in it is called {\bf \em aggregation-free}.\\
(ii) If $\varphi(\bar{x}) \in PLA^*(\sigma)$ and there is a finite $V \subseteq [0, 1]$
such that for every finite $\sigma$-structure $\mcA$,
 and every $\bar{a} \in A^{|\bar{x}|}$,
$\mcA(\varphi(\bar{a})) \in V$, then we call $\varphi(\bar{x})$ {\bf \em finite valued}.
In the special case when, additionally, $V = \{0, 1\}$
then we call  $\varphi(\bar{x})$ {\bf \em 0/1-valued}.\\
(iii) If $L \subseteq PLA^*(\sigma)$ and every formula in $L$ is 0/1-valued, then 
we say that $L$ is {\bf \em $0/1$-valued}. \\
(iv) Let $L \subseteq PLA^*(\sigma)$ be $0/1$-valued. 
A formula of $PLA^*(\sigma)$ is called an {\bf \em $L$-basic (formula)} if it has the form
$\bigwedge_{i=1}^k \big(\varphi_i(\bar{x}) \to c_i\big)$ where $\varphi_i \in L$ and $c_i \in [0, 1]$ for all $i = 1, \ldots, k$,
and $\lim_{n\to\infty}\mbbP_n\big(\forall \bar{x} \bigvee_{i=1}^k \varphi_i(\bar{x})\big) = 1$.
}\end{defin}

\begin{defin}\label{definition of equivalence} {\rm
Let $\varphi(\bar{x}), \psi(\bar{x}) \in PLA^*(\sigma)$. We say that $\varphi$ and $\psi$
are {\bf \em equivalent} if for every finite $\sigma$-structure $\mcA$ and every $\bar{a} \in A^{|\bar{x}|}$,
$\mcA(\varphi(\bar{a})) = \mcA(\psi(\bar{a}))$.
}\end{defin}

\begin{notation}\label{Notation for consequence in finite}{\rm
(i) For any formula $\varphi(\bar{x}, \bar{y}) \in PLA^*(\sigma)$, finite $\sigma$-structure $\mcA$ and $\bar{a} \in A^{|\bar{x}|}$, let
\[
\varphi(\bar{a}, \mcA) = \{\bar{b} \in A^{|\bar{y}|} : \mcA(\varphi(\bar{a}, \bar{b})) = 1 \}.
\]
(ii) Let $\varphi(\bar{x}), \psi(\bar{x}) \in PLA^*(\sigma)$. 
When writing 
\[
\varphi(\bar{x}) \models \psi(\bar{x})
\] 
we mean that for every finite $\sigma$-structure $\mcA$
and $\bar{a} \in A^{|\bar{x}|}$, if $\mcA(\varphi(\bar{a})) = 1$ then $\mcA(\psi(\bar{a})) = 1$.
(We will only use the notation together with $0/1$-valued formulas.)
}\end{notation}

\begin{exam}\label{example of page rank} {\rm
We exemplify what can be expressed with $PLA^*(\sigma)$, provided that it
contains a binary relation symbol, with the notion of PageRank \cite{BP}.
The PageRank of an internet site can be approximated in ``stages'' as follows
(if we supress the ``damping factor'' for simplicity), where $IN_x$ is the set of sites that link to $x$,
and $OUT_y$ is the set of sites that $y$ link to:
\begin{align*}
&PR_0(x) = 1/N \text{ where $N$ is the number of sites,} \\
&PR_{k+1}(x) = \sum_{y \in IN_x} \frac{PR_k(y)}{|OUT_y|}.
\end{align*}
It is not difficult to prove, by induction on $k$, that for every $k$ the sum of all $PR_k(x)$ as $x$ ranges over all sites
is 1. Hence the sum in the definition of $PR_{k+1}$ is less or equal to 1 (which will matter later).
Let $E \in \sigma$ be a binary relation symbol representing a link.
Define the aggregation function $\mr{length}^{-1} : [0, 1]^{<\omega} \to [0, 1]$ by $\mr{length}^{-1}(\bar{p}) = 1/|\bar{p}|$.
Then $PR_0(x)$ is expressed by the $PLA^*(\sigma)$-formula $\mr{length}^{-1}(x = x : y : \top)$.

Suppose that $PR_k(x)$ is expressed by $\varphi_k(x) \in PLA^*(\sigma)$.
Note that multiplication is a continuous connective from $[0, 1]^2$ to $[0, 1]$ so it can be used in $PLA^*(\sigma)$-formulas.
Then observe that the quantity $|OUT_y|^{-1}$ is expressed by the $PLA^*(\sigma)$-formula
\[
\mr{length}^{-1}\big(y=y : z : E(y, z)\big)
\]
which we denote by $\psi(y)$.
Let $\mr{tsum} : [0, 1]^{<\omega} \to [0, 1]$ be the ``truncated sum'' defined by letting
$\mr{tsum}(\bar{p})$ be the sum of all entries in $\bar{p}$ if the sum is at most 1, and otherwise $\mr{tsum}(\bar{p}) = 1$.
Then $PR_{k+1}(x)$ is expressed by the $PLA^*(\sigma)$-formula
\[
\mr{tsum}\big(x = x \wedge (\varphi_k(y) \cdot \psi(y)) : y : E(y, x)\big).
\]
With $PLA^*(\sigma)$ we can also define all stages of the SimRank \cite{JW} in a simpler way than done
in \cite{KW1} with the sublogic $PLA(\sigma) \subseteq PLA^*(\sigma)$.
}\end{exam}

\begin{rem}\label{FO is expressible in PLA*}  {\rm
Suppose that $\varphi(\bar{x}, \bar{y}) \in FO(\sigma)$,  $\psi(\bar{x}, \bar{y}) \in PLA^*(\sigma)$
and for every finite $\sigma$-structure $\mcA$, $\bar{a} \in A^{|\bar{x}|}$
and $\bar{b} \in A^{|\bar{y}|}$,
$\mcA \models \varphi(\bar{a}, \bar{b})$ if and only if $\mcA\big(\psi(\bar{a}, \bar{b})\big) = 1$.
Then, for all $\bar{a} \in A^{|\bar{x}|}$,
\[
\mcA \models \exists \bar{y} \varphi(\bar{a}, \bar{y}) \text{ if and only if } 
\mcA\big(\max(\psi(\bar{a}, \bar{y}) : \bar{y} : \top)\big) = 1.
\]
Similarly, the quantifier $\forall$ can be expressed in $PLA^*(\sigma)$ by using  the aggregation function min.
By induction on the complexity of first-order formulas it follows that for every $\varphi(\bar{x}) \in FO(\sigma)$
there is a 0/1-valued $\psi(\bar{x}) \in PLA^*(\sigma)$ such that $\varphi$ and $\psi$ are equivalent.
}\end{rem}

\begin{notation}\label{using first-order notation} {\bf (Using $\exists$ and $\forall$ as abbreviations)} {\rm
Due to Remark~\ref{FO is expressible in PLA*}, if $\varphi(\bar{x}, \bar{y}) \in PLA^*(\sigma)$ 
is a 0/1-valued formula
then we will often (in particular in Sections~\ref{the base sequence} and~\ref{properties of closure types})
write
`$\exists \bar{y} \varphi(\bar{x}, \bar{y})$' to mean the same as
`$\max(\varphi(\bar{x}, \bar{y}) : \bar{y} : \top)$', and
`$\forall \bar{y} \varphi(\bar{x}, \bar{y})$' to mean the same as
`$\min(\varphi(\bar{x}, \bar{y}) : \bar{y} : \top)$'.
}\end{notation}

\noindent
The next basic lemma has analogue in first-order logic and is proved straightforwardly by induction on the complexity of 
$PLA^*$-formulas.

\begin{lem}\label{truth values only depend on the relation symbols used}
Suppose that $\sigma' \subseteq \sigma$,
$\varphi(\bar{x}) \in PLA^*(\sigma')$, $\mcA$ is a finite $\sigma$-structure, $\mcA' = \mcA \uhrc \sigma'$,
and $\bar{a} \in A^{|\bar{x}|}$.
Then $\mcA(\varphi(\bar{a})) = \mcA'(\varphi(\bar{a}))$.
\end{lem}

\section{A general approach to asymptotic elimination of aggregation functions}\label{A general approach}

\noindent
In this section we recall the main parts of the general approach to ``asymptotic elimination of aggregation functions''
that is described by the author and Weitkämper in \cite{KW3}.
The set up will be as follows throughout this section:
We assume that $\sigma$ is a finite relational signature and for each 
$n \in \mbbN^+$, $D_n$ is a finite set such that $\lim_{n\to\infty}|D_n| = \infty$.
(We do not assume that $D_n \subseteq D_{n+1}$.)
Furthermore, we assume that (for each $n \in \mbbN^+$) $\mbW_n$ is a nonempty set of 
(not necessarily all)
$\sigma$-structures with domain $D_n$.
We begin by defining the notion of ``asymptotic equivalence (with respect to a sequence of probability distributions)''  between formulas.
After that we consider notions of continuity of {\em aggregation} functions, 
in particular the notion that we simply call ``continuity'' and the weaker notion ``admissibility''.
The aggregation functions arithmetic mean and geometric mean are continuous while maximum and minimum are 
admissible (but not continuous), as stated in Example~\ref{example of continuous aggregation functions}.
Very roughly speaking, the main result of this section, Theorem~\ref{general asymptotic elimination},
tells that under some conditions all continuous aggregation functions can be asymptotically eliminated, 
and under somewhat stronger conditions all admissible aggregation functions can be asymptotically eliminated.
Theorem~\ref{general asymptotic elimination} will be used in 
Section~\ref{Asymptotic elimination of aggregation functions}
to ``asymptotically eliminate'' aggregation functions in the context of working with expansions of a sequence
of base structures with bounded degree.
The work of sections~\ref{Proving convergence in the inductive step} and~\ref{Finding the balance in the inductive step}
will show that the conditions of 
Assumption~\ref{assumptions on the basic logic} below
are satisfied within the context of the kind of base structures considered in this article and consequently
Theorem~\ref{general asymptotic elimination}
can be applied.

\begin{defin}\label{definition of asymptotic probability distribution}{\rm
By a {\bf \em sequence of probability distributions (for $(\mbW_n : n \in \mbbN^+)$)} we mean a sequence
$\mbbP = (\mbbP_n : n \in \mbbN^+)$ such that for every $n$, $\mbbP_n$ is a probability distribution on $\mbW_n$.
}\end{defin}

\begin{defin}\label{definition of asymptotically equivalent formulas} {\rm
Let $\varphi(\bar{x}), \psi(\bar{x}) \in  PLA^*(\sigma)$ and 
let $\mbbP = (\mbbP_n : n \in \mbbN^+)$ be a sequence of probability distributions.
We say that $\varphi(\bar{x})$ and $\psi(\bar{x})$ are {\bf \em asymptotically equivalent (with respect to $\mbbP$)} 
if for all $\varepsilon > 0$
\[
\mbbP_n\Big(\big\{\mcA \in \mbW_n : \text{ there is $\bar{a} \in (B_n)^{|\bar{x}|}$ such that 
$|\mcA(\varphi(\bar{a})) - \mcA(\psi(\bar{a}))| > \varepsilon$}\big\} \Big) \to 0
\]
as $n \to \infty$.
}\end{defin}

\noindent
{\bf \em For the rest of this section we fix a sequence of probability distributions
$\mbbP = (\mbbP_n : n \in \mbbN^+)$. When saying that formulas are asymptotically equivalent, we mean 
``with respect to $\mbbP$''.}

In order to define the notion of continuity that we will use we need the notion of
{\em convergence testing sequence} which generalizes a similar
notion used by Jaeger in \cite{Jae98a}.
The intuition is that a sequence $\bar{p}_n \in [0, 1]^{<\omega}$, $n \in \mbbN^+$ is convergence testing
for parameters 
$c_1, \ldots, c_k \in [0,1]$ and $\alpha_1, \ldots \alpha_k  \in  [0,1]$
if the length of $\bar{p}_n$ tends to infinity as $n\to\infty$ and, as $n\to\infty$,
all entries in $\bar{p}_n$ congregate ever closer to the ``convergence points" in the set $\{c_1, \ldots, c_k\}$,
and the proportion of entries in $\bar{p}$ which are close to $c_i$ converges to $\alpha_i$.

\begin{defin}\label{definition of convergence testing} {\rm 
(i) A sequence $\bar{p}_n \in [0, 1]^{<\omega}$, $n \in \mbbN$,  is called {\bf \em convergence testing for parameters} 
$c_1, \ldots, c_k \in [0,1]$ and $\alpha_1, \ldots \alpha_k  \in  [0,1]$ if the following hold, 
where $p_{n,i}$ denotes the $i$th entry of $\bar{p}_n$:
\begin{enumerate}
\item $|\bar{p}_n| < |\bar{p}_{n+1}|$ for all $n \in \mbbN$.
\item For every disjoint family of open (with respect to the induced topology on $[0, 1]$)
intervals $I_1, \ldots I_k \subseteq [0,1]$ 
 such that $c_i \in I_i$ for each $i$, 
there is an $ N \in \mbbN$ such that $\mathrm{rng}(\bar{p}_n) \subseteq \bigcup\limits_{j=1}^{k} I_j$ for all $n \geq N$, 
and for every $j \in \{1, \ldots, k \}$, 
\[
\lim\limits_{n \rightarrow \infty} \frac{\left| \{ i  : p_{n,i} \in I_j \} \right| }{|\bar{p}_n|} = \alpha_j.
\]
\end{enumerate}   
(ii) More generally, a sequence of $m$-tuples 
$(\bar{p}_{1, n}, \ldots, \bar{p}_{m, n}) \in  \big([0, 1]^{<\omega}\big)^m$, $n \in \mbbN$,  
is called {\bf \em convergence testing for parameters} $c_{i,j} \in [0,1]$ and $\alpha_{i,j} \in [0,1]$, 
where $i \in \{1, \ldots, m\}$, $j \in \{ 1, \ldots, k_i \}$ and $k_1, \ldots k_m \in \mbbN^+$, if  for every fixed 
$i \in \{1, \ldots, m \}$
the sequence $\bar{p}_{i, n}$, $n \in \mbbN$, is convergence testing for $c_{i,1}, \ldots, c_{i, k_i}$, and 
$\alpha_{i,1}, \ldots, \alpha_{i, k_i}$.
}\end{defin}

\begin{defin} \label{definition of ct-continuous} {\rm 
An aggregation function \\ 
$F : \big([0, 1]^{<\omega}\big)^m \to [0, 1]$ is called 
{\bf \em ct-continuous (convergence test continuous)} with respect to the sequence of parameters
$c_{i,j}, \alpha_{i,j} \in [0, 1]$, $i = 1, \ldots, m$, $j = 1, \ldots, k_i$, 
if the following condition holds:
\begin{enumerate}
\item[] For all convergence testing sequences of $m$-tuples
$(\bar{p}_{1, n}, \ldots, \bar{p}_{m, n}) \in  \big([0, 1]^{<\omega}\big)^m$, $n \in \mbbN$,
and $(\bar{q}_{1, n}, \ldots, \bar{q}_{m, n}) \in  \big([0, 1]^{<\omega}\big)^m$, $n \in \mbbN$,
with the same parameters $c_{i,j}, \alpha_{i,j} \in [0, 1]$, 
$\underset{n \rightarrow \infty}{\lim}  
|F(\bar{p}_{1, n}, \ldots, \bar{p}_{m, n}) - F(\bar{q}_{1, n}, \ldots, \bar{q}_{m, n})| = 0$.
\end{enumerate}
}\end{defin}

\begin{defin}\label{definition of strong admissibility} {\rm
Let $F : \big([0, 1]^{<\omega}\big)^m \to [0, 1]$.\\
(i) We call $F$ {\bf \em continuous}, or {\bf \em strongly admissible}, if $F$ is ct-continuous with respect to {\em every} choice
of parameters $c_{i, j}, \alpha_{i, j} \in [0, 1]$, $i = 1, \ldots, m$ and $j = 1, \ldots, k_i$ (for arbitrary $m$ and $k_i$).\\
(ii) We call $F$ {\bf \em admissible} if $F$ is ct-continuous with respect to every choice of parameters
$c_{i, j}, \alpha_{i, j} \in [0, 1]$, $i = 1, \ldots, m$ and $j = 1, \ldots, k_i$ (for arbitrary $m$ and $k_i$)
{\em such that $\alpha_{i, j} > 0$ for all $i$ and $j$}.
}\end{defin}

\begin{exam}\label{example of continuous aggregation functions} {\rm
The aggregation functions am, gm and $\mr{length}^{-\beta}$ are continuous,
which is proved in \cite{KW1} in the case of am and gm, and in the case of $\mr{length}^{-\beta}$ the claim is easy to prove.
The aggregation functions max and min are admissible (which is proved in \cite{KW1}) but not continuous
(which is easy to see).
The aggregation function $\text{noisy-or}((p_1, \ldots, p_n)) = 1 - \prod_{i=1}^n(1 - p_n)$ is not even admissible
(which is not hard to prove).
For more examples of admissible, or even continuous, aggregation functions (of higher arity) see
Example~6.4 and Proposition~6.5 in \cite{KW1}.
}\end{exam}

\noindent
The intuition behind part~(2) of 
Assumption~\ref{assumptions on the basic logic} 
below is that for the sets $L_0, L_1 \subseteq PLA^*(\sigma)$ of $0/1$-valued formulas
and every $\varphi(\bar{x}, \bar{y}) \in L_0$ there is a 
set $L_{\varphi(\bar{x}, \bar{y})} \subseteq L_1$ of formulas expressig 
some ``allowed'' conditions (with respect to $\varphi(\bar{x}, \bar{y})$)
and there are some 
$\varphi'_1(\bar{x}), \ldots, \varphi'_s(\bar{x}) \in L_0$ such that if $\mcA \models \varphi'_i(\bar{a})$
and $\chi(\bar{x}, \bar{y}) \in L_{\varphi(\bar{x}, \bar{y})}$, 
then the proportion $|\varphi(\bar{a}, \mcA) \cap \chi(\bar{a}, \mcA)| / |\chi(\bar{a}, \mcA)|$ is almost surely close to a 
number $\alpha_i$ that depends only on 
$\varphi(\bar{x}, \bar{y}), \chi(\bar{x}, \bar{y}), \varphi_i(\bar{x})$
and the sequence of probability distributions $\mbbP$.
As we allow aggregation functions with arity $m > 1$, part~(2) needs to simultaneously
speak of a sequence  $\varphi_1(\bar{x}, \bar{y}), \ldots, \varphi_m(\bar{x}, \bar{y})\in L_0$.

\begin{assump}\label{assumptions on the basic logic}{\rm
Suppose that $L_0 \subseteq PLA^*(\sigma)$ and $L_1 \subseteq PLA^*(\sigma)$ 
are $0/1$-valued and that the following hold:
\begin{enumerate}
\item For every aggregation-free $\varphi(\bar{x}) \in PLA^*(\sigma)$ there is an
$L_0$-basic formula $\varphi'(\bar{x})$ such that
$\varphi$ and $\varphi'$ are asymptotically equivalent.

\item For every $m \in \mbbN^+$ and $\varphi_j(\bar{x}, \bar{y})\in L_0$, for $j = 1, \ldots, m$,
there are $L_{\varphi_j(\bar{x}, \bar{y})} \subseteq L_1$ for $j = 1, \ldots, m$
such that if $\chi_j(\bar{x}, \bar{y}) \in L_{\varphi_j(\bar{x}, \bar{y})}$
for $j = 1, \ldots, m$, then 
there are $s, t \in \mbbN^+$, $\varphi'_i(\bar{x}) \in L_0$, $\alpha_{i, j} \in [0, 1]$,
for $i = 1, \ldots, s$, $j = 1, \ldots, m$,
and $\chi'_i(\bar{x}) \in L_0$, for $i = 1, \ldots, t$,
such that for every $\varepsilon > 0$ and $n$ there is $\mbY^\varepsilon_n \subseteq \mbW_n$ such that
$\lim_{n\to\infty}\mbbP_n(\mbY^\varepsilon_n) = 1$ and 
for every $\mcA \in \mbY^\varepsilon_n$ the following hold:
\begin{align*}
&(a) \ \mcA \models \forall \bar{x} \bigvee_{i=1}^s \varphi'_i(\bar{x}), \\
&(b) \ \text{if $i \neq j$ then } 
\mcA \models \forall \bar{x} \neg(\varphi'_i(\bar{x}) \wedge \varphi'_j(\bar{x})), \\
&(c) \ \mcA \models \forall \bar{x}
\Big(\Big(\bigvee_{i=1}^m \neg\exists\bar{y}\chi_i(\bar{x}, \bar{y})\Big) \leftrightarrow
\Big(\bigvee_{i=1}^t \chi'_i(\bar{x})\Big)\Big), \text{ and}\\
&(d) \ \text{for all $i = 1, \ldots, s$ and $j = 1, \ldots, m$, if $\bar{a} \in (D_n)^{|\bar{x}|}$, and
$\mcA \models \varphi'_i(\bar{a})$,} \\
&\text{ then } (\alpha_{i, j} - \varepsilon)|\chi_j(\bar{a}, \mcA)| \ \leq \
|\varphi_j(\bar{a}, \mcA) \cap \chi_j(\bar{a}, \mcA)|
\ \leq \ (\alpha_{i, j} + \varepsilon)|\chi_j(\bar{a}, \mcA)|.
\end{align*}
\end{enumerate}
}\end{assump}

\noindent
The next result will be used in 
Section~\ref{Asymptotic elimination of aggregation functions}
and is a simplification of 
Theorem~5.9 in \cite{KW3}. 
See remarks about how it follows from \cite{KW3} after its statement.

\begin{theor}\label{general asymptotic elimination} {\rm \cite{KW3}} 
Suppose that $L_0, L_1 \subseteq PLA^*(\sigma)$ are 0/1-valued and that 
Assumption~\ref{assumptions on the basic logic} holds.
Let $F : \big([0, 1]^{<\omega}\big)^m \to [0, 1]$, 
let $\psi_i(\bar{x}, \bar{y}) \in PLA^*(\sigma)$, for $i = 1, \ldots, m$, and suppose that each 
$\psi_i(\bar{x}, \bar{y})$
is asymptotically equivalent to an $L_0$-basic formula \\
$\bigwedge_{k=1}^{s_i} (\psi_{i, k}(\bar{x}, \bar{y}) \to c_{i, k})$
(so $\psi_{i, k} \in L_0$ for all $i$ and $k$).
Suppose that for $i = 1, \ldots, m$, 
$\chi_i(\bar{x}, \bar{y}) \in \bigcap_{k=1}^{s_i} L_{\psi_{i, k}(\bar{x}, \bar{y})}$.
Let $\varphi(\bar{x})$ denote the $PLA^*(\sigma)$-formula
\[
F\big(\psi_1(\bar{x}, \bar{y}), \ldots, \psi_m(\bar{x}, \bar{y}) : \bar{y} : 
\chi_1(\bar{x}, \bar{y}), \ldots, \chi_m(\bar{x}, \bar{y})\big).
\]
(i) If $F$ is continuous then $\varphi(\bar{x})$ is asymptotically equivalent to an $L_0$-basic formula.\\
(ii) Suppose, in addition, that the following holds 
if $\varphi_j(\bar{x}, \bar{y}) \in L_0$, $\chi_j(\bar{x}, \bar{y}) \in L_1$, $\varphi'_i(\bar{x}) \in L_0$,
$\mbY^\varepsilon_n$ and $\alpha_{i, j}$ are  like in 
part~2 of Asumption~\ref{assumptions on the basic logic}:
If $\alpha_{i, j} = 0$
then, for all sufficiently large $n$, all $\bar{a} \in (D_n)^{|\bar{x}|}$ and all $\mcA \in \mbY^\varepsilon_n$,
if $\mcA \models \varphi'_i(\bar{a})$
then $\varphi_j(\bar{a}, \mcA) \cap \chi_j(\bar{a}, \mcA) = \es$.
Then it follows that if $F$ is admissible then $\varphi(\bar{x})$ is asymptotically equivalent to an $L_0$-basic formula.
\end{theor}

\begin{rem}\label{remark on general asymptotic elimination}{\rm
Part~(i) of Theorem~\ref{general asymptotic elimination}
follows from the definition of continuous aggregation function and Corollary~4.11 and Theorem~5.9 in \cite{KW3}.
Part~(ii) follows by the following observation:
If, in the context of part~(2) of Definition~5.7 in \cite{KW3},
we have $\psi_j(\bar{a}, \mcA) \cap \chi(\bar{a}, \mcA) = \es$ whenever 
$\mcA \models \theta_i(\bar{a})$ and $\alpha_{i, j} = 0$,
then we can omit $c_j$ and $\alpha_{i, j}$ before constructing the 
``sequence of $\bar{y}$-frequency parameters of $\psi$ relative to $\chi$ and $\theta_i$''
(as constructed in \cite[Definition~5.7]{KW3}).
In this way all $\alpha_{i, j}$ that remain are positive.
}\end{rem}

\section{Properties of the base sequence of structures}\label{the base sequence}

\noindent
In this section we will state the conditions on the sequence of base structures 
(Assumption~\ref{properties of the base structures}) that will be assumed in the rest of this study.
These conditions use the notion of {\em $(\sigma, \lambda)$-neighbourhood type} and the 
notion of {\em $(\sigma, \lambda)$-closure type},
where $\sigma$ is a signature and $\lambda \in \mbbN$. 
If we work with a signature $\sigma$, a $(\sigma, \lambda)$-neighbourhood type, 
respectively $(\sigma, \lambda)$-closure type, is a formula
that describes the isomorphism type of the substructure induced
by the $\lambda$-neighbourhood, respectively $\lambda$-closure, of some elements.
The $\lambda$-neighbourhood of some elements is the set of elements within distance $\lambda$
from those elements in a sense that will be explained.
We also prove some basic results about neighbourhood and closure types.
A more detailed study of neighbourhood and closure types will be carried out in 
Section~\ref{properties of closure types}
after we have looked at some examples of sequences of base structures in 
Section~\ref{Examples of sequences of base structures}.
The notions of {\em bounded} and {\em unbounded} neighbourhood and closure type will be critical.
At the end of the section we define the notion of {\em $(\sigma, \lambda)$-basic formula} which, roughly speaking, 
is a $PLA^*(\sigma)$-formula of the form $\bigcup_{i=1}^k (\varphi_i(\bar{x}) \to c_i)$ where
$\varphi_i(\bar{x})$ is a $(\sigma, \lambda)$-basic formula and $c_i \in [0, 1]$ for all $i = 1, \ldots, k$.
The relevance of this notion is that our main results (in Section~\ref{Asymptotic elimination of aggregation functions})
say that, under some conditions, a $PLA^*(\sigma)$-formula is asymptotically equivalent to a
$(\sigma, \lambda)$-basic formula, for some $\lambda \in \mbbN$.

{\bf \em In the rest of the article we assume that $\tau$ and $\sigma$ are finite relational signatures and $\tau \subseteq \sigma$.}
All definitions below make sense if $\tau = \sigma$.
In this and the two following sections we mostly work with the signature $\tau$, 
the signature of the base structures,
but we define some notions (neighbourhood and closure types) 
for the possibly larger signature $\sigma$ since the variant of these notions for $\sigma$
will be used later (starting from Section~\ref{convergence and balance})

We begin by generalizing, in a familiar way, some notions from graph theory to relational structures in general.

\begin{defin}\label{definition of degree}{\rm
Let $\mcB$ be a finite $\tau$-structure.
\begin{enumerate}
\item The {\bf \em Gaifman graph of $\mcB$} is the undirected graph $\mcN$ defined as follows:
\begin{enumerate}
\item The vertex set of $\mcN$ is $B$ (the domain of $\mcB$).
\item Let $a, b \in B$. There is an edge (of $\mcN$) between $a$ and $b$ if and only if $a \neq b$ and
there is $R \in \tau$, of arity $r$ say, and $c_1, \ldots, c_r \in B$ such that $\mcB \models R(c_1, \ldots, c_r)$
and $a, b \in \{c_1, \ldots, c_r\}$.
\end{enumerate}

\item Let $b \in B$. The {\bf \em degree of $b$ (with respect to $\mcB$)}, denoted $\deg_\mcB(b)$, is
the number of neighbours of $b$ in the Gaifman graph of $\mcB$. 

\item The {\bf \em degree of $\mcB$}, denoted $\deg(\mcB)$, 
is the maximum of $\deg_\mcB(b)$ as $b$ ranges over $B$.
\end{enumerate}
}\end{defin}

\begin{defin}\label{definition of distance}{\rm
Let $\mcB$ be a finite $\tau$-structure.
\begin{enumerate}
\item Let $a, b \in B$. 
The {\bf \em distance between $a$ and $b$ (in $\mcB$)}, denoted $\dist_\mcB(a, b)$
or just $\dist(a, b)$, is, by definition, equal to the distance between $a$ and $b$ in the Gaifman graph of $\mcB$
(in other words, it equals the length of the shortest path from $a$ to $b$ in the Gaifman graph of $\mcB$ if 
such a path exists, and otherwise the distance is stipulated to be $\infty$).

\item If $\bar{a} = (a_1, \ldots, a_k) \in B^k$ and $\bar{b} = (b_1, \ldots, b_l) \in B^l$, then
\[
\dist_\mcB(\bar{a}, \bar{b}) = \min \{\dist_\mcB(a_i, b_j) : i = 1, \ldots, k, \ j = 1, \ldots l\}
\]
which may be abbreviated by $\dist(\bar{a}, \bar{b})$ if the structure $\mcB$ is clear from the context.
\item If $\mcA$ is a finite $\sigma$-structure and $\bar{a}$ and $\bar{b}$ are finite sequences of elements from $A$,
then we define $\dist_\mcA(\bar{a}, \bar{b}) = \dist_{\mcA \uhr \tau}(\bar{a}, \bar{b})$.
\end{enumerate}
}\end{defin}

\noindent
{\bf Warning:} Although clear from the definition above, I want to emphasize that even if we work
in a $\sigma$-structure $\mcA$ the distance between elements is always computed in the reduct $\mcA \uhrc \tau$
to the signature $\tau$.

\begin{lem}\label{definability of distance}
For all $\lambda \in \mbbN$ there are $\varphi_\lambda(x, y) \in FO(\tau)$ and a 0/1-valued formula
$\psi_\lambda(x, y) \in PLA^*(\tau)$ such that, for every finite $\tau$-structure $\mcB$ and all $a, b \in B$,
$\dist_\mcB(a, b) \leq \lambda$ if and only $\mcB \models \varphi_\lambda(a, b)$ if and only if 
$\mcA(\psi_\lambda(a, b)) = 1$.
\end{lem}

\noindent
{\bf Proof.}
That the relation `$\dist_\mcB(x, y) \leq \lambda$' is first-order definable is well-known and can be proved
straightforwardly by induction on $\lambda$.
The claim about $PLA^*$-formulas follows from Remark~\ref{FO is expressible in PLA*}.
\hfill $\square$

\begin{notation}\label{remark on definability of distance}{\rm
Due to Lemma~\ref{definability of distance}
we will often use the expression `$\dist(x, y) \leq \lambda$' to denote a 
$PLA^*(\tau)$-formula that expresses (in every finite $\tau$-structure) that ``the distance between $x$ and $y$ 
 is at most $\lambda$''. The expression `$\dist(x, y) > \lambda$' denotes the negation  of such a formula.
For sequences $\bar{x} = (x_1, \ldots, x_k)$ and $\bar{y} = (y_1, \ldots, y_l)$ we let the expression
`$\dist(\bar{x}, \bar{y}) \leq \lambda$' denote the formula
$\bigvee_{i=1}^k \bigvee_{j=1}^l \dist(x_i, y_j) \leq \lambda$, and `$\dist(\bar{x}, \bar{y}) > \lambda$' the negation of it.
}\end{notation}

\begin{defin}\label{definition of neighbourhood}{\rm
Let $\mcB$ be a finite $\tau$-structure, $\lambda \in \mbbN$, and $b_1, \ldots, b_k \in B$.
The {\bf \em $\lambda$-neighbourhood of $b$ (with respect to $\mcB$)}, is the set 
\[
N_\lambda^\mcB(b_1, \ldots, b_k) = \{a \in B : \text{ for some } i \in \{1, \ldots, k\}, \ \dist_\mcB(a, b_i) \leq \lambda \}.
\]
}\end{defin}

\noindent
Let $\lambda \in \mbbN$.
The idea with the next definition is that a (0/1-valued) formula $p(\bar{x})$ is a ``complete $(\sigma, \lambda)$-neighbourhood type''
if the following holds:
If $\mcA_1$ and $\mcA_2$ are finite $\sigma$-structures, $\mcB_1 = \mcA_1 \uhrc \tau$, $\mcB_2 = \mcA_2 \uhrc \tau$,
$\mcA_1 \models p(\bar{a}_1)$ and $\mcA_2 \models p(\bar{a}_2)$, 
then there is an isomorphism $f$ from $\mcA_1 \uhrc N_\lambda^{\mcB_1}(\bar{a}_1)$ to 
$\mcA_2 \uhrc N_\lambda^{\mcB_2}(\bar{a}_2)$ such that $f(\bar{a}_1) = \bar{a}_2$.

\begin{defin}\label{definition of neighbourhood type} {\rm
\begin{enumerate}
\item A {\bf \em $(\sigma, 0)$-neighbourhood type} in the variables $\bar{x} = (x_1, \ldots, x_k)$
is a consistent conjunction of first order $\sigma$-literals 
(see Definition~\ref{definition of literal}) with (only) variables from $\bar{x}$ such that
\begin{enumerate}
\item for every $R \in \tau$ and every choice of $x_{i_1}, \ldots, x_{i_r}$ where $r$ is the arity of $R$,
either $R(x_{i_1}, \ldots, x_{i_r})$ or $\neg R(x_{i_1}, \ldots, x_{i_r})$ is a conjunct, and

\item for all $1 \leq i < j \leq k$, either $x_i = x_j$ or $x_i \neq x_j$ is a conjunct.
\end{enumerate}

\item A {\bf \em complete $(\sigma, 0)$-neighbourhood type} in the variables $\bar{x} = (x_1, \ldots, x_k)$
is a consistent conjunction of $\sigma$-literals with (only) variables from $\bar{x}$ such that
\begin{enumerate}
\item for every $R \in \sigma$ and every choice of $x_{i_1}, \ldots, x_{i_r}$ where $r$ is the arity of $R$,
either $R(x_{i_1}, \ldots, x_{i_r})$ or $\neg R(x_{i_1}, \ldots, x_{i_r})$ is a conjunct, and

\item for all $1 \leq i < j \leq k$, either $x_i = x_j$ or $x_i \neq x_j$ is a conjunct.
\end{enumerate}
Note that the only difference compared with~(1) is that $\tau$ is replaced by $\sigma$ in part~(a).

\item For $\lambda \in \mbbN^+$, a {\bf \em $(\sigma, \lambda)$-neighbourhood type} in the variables 
$\bar{x} = (x_1, \ldots, x_k)$ 
is a consistent formula of the following form
(or an equivalent $PLA^*(\sigma)$-formula), for some $(\sigma, 0)$-neighbourhood type $p(\bar{x}, y_1, \ldots, y_l)$
(in the variables $\bar{x}, y_1, \ldots, y_l$):
\begin{align*}
\exists y_1, \ldots, y_l \bigg(
&\bigwedge_{i=1}^l (\dist(y_i, \bar{x}) \leq \lambda) \  \wedge \
\forall z \Big( (\dist(z , \bar{x}) \leq \lambda) \rightarrow \bigvee_{i = 1}^l z = y_i \Big) \\ 
&\wedge 
p(\bar{x}, y_1, \ldots, y_l) \bigg).
\end{align*}

\item A {\bf \em complete $(\sigma, \lambda)$-neighbourhood type} is defined like a 
$(\sigma, \lambda)$-neighbourhood type, except that we require that $p(\bar{x}, y_1, \ldots, y_l)$ is 
a {\em complete} $(\sigma, 0)$-neighbourhood type.
\end{enumerate}
}\end{defin}

\noindent
Note that since we can let $\sigma = \tau$ we can talk about $(\tau, \lambda)$-neighbourhood types.
Also observe that it follows from the above definition that (for any $\lambda \in \mbbN)$ 
every $(\tau, \lambda)$-neighbourhood type 
is a {\em complete} $(\tau, \lambda)$-neighbourhood type. 
(But if $\tau$ is a proper subset of $\sigma$ then a $(\sigma, \lambda)$-neighbourhood type need not be a
complete $(\sigma, \lambda)$-neighbourhood type.)
Nevertheless we will usually say ``{\em complete} $(\tau, \lambda)$-neighbourhood type'' to emphasize that 
it completely describes the $\tau$-structure of a $\lambda$-neighbourhood.

\begin{lem}\label{neighbourhood types and distance}
Let $p(x_1, \ldots, x_k)$ be a complete $(\tau, \lambda)$-neighbourhood type for some $\lambda \in \mbbN$.
Then for all $i, j \in \{1, \ldots, k\}$ either
\begin{enumerate}
\item[] $p(x_1, \ldots, x_k) \models \dist(x_i, x_j) \leq 2\lambda$, or
\item[] $p(x_1, \ldots, x_k) \models \dist(x_i, x_j) > 2\lambda$.
\end{enumerate}
\end{lem}

\noindent
{\bf Proof.} Immediate from the definition of $(\tau, \lambda)$-neighbourhood type.
\hfill $\square$

\medskip 

\noindent
The next lemma intuitively says that if elements $\bar{a}$ and elements $\bar{b}$ are far apart then the neighbourhood
types of $\bar{a}$ and $\bar{b}$ determine the neighbourhood type of the concatenated sequence $\bar{a}\bar{b}$.

\begin{lem}\label{splitting a neighbourhood type}
Let $p(\bar{x}, \bar{y})$ be a complete $(\tau, \lambda)$-neighbourhood type for some $\lambda \in \mbbN$, and
suppose that $p(\bar{x}, \bar{y}) \models \dist(\bar{x}, \bar{y}) > 2\lambda$.
Let $p_1(\bar{x}) = p \uhrc \bar{x}$ and $p_2(\bar{y}) = p \uhrc \bar{y}$.
Then, for all $n \in \mbbN^+$, for all $\bar{a} \in (B_n)^{|\bar{x}|}$ and $\bar{b} \in (B_n)^{|\bar{y}|}$,
\begin{enumerate}
\item[] $\mcB_n \models p(\bar{a}, \bar{b})$ if and only if $\mcB_n \models p_1(\bar{a}) \wedge p_2(\bar{b})$ and
$\dist(\bar{a}, \bar{b}) > 2\lambda$.
\end{enumerate}
\end{lem}

\noindent
{\bf Proof.}
Note that if $p(\bar{x}, \bar{y}) \models \dist(\bar{x}, \bar{y}) > 2\lambda$
and $\mcB_n \models p(\bar{a}, \bar{b})$ then $N_\lambda^{\mcB_n}(\bar{a}) \cap N_\lambda^{\mcB_n}(\bar{b}) = \es$.
Now the result follows from the definition of $(\tau, \lambda)$-neighbourhood type.
\hfill $\square$

\begin{defin}\label{definition of dimension}{\rm
Let $\lambda \in \mbbN$ and let $p(\bar{x})$ be a complete $(\tau, \lambda)$-neighbourhood type
where $\bar{x} = (x_1, \ldots, x_k)$.
\begin{enumerate}
\item We define a relation $\approx_p$ on $\rng(\bar{x}) = \{x_1, \ldots, x_k\}$ by
\[
x_i \approx_p x_j \ \Longleftrightarrow \ p(\bar{x}) \models \dist(x_i, x_j) \leq 2 \lambda.
\]

\item Let $\sim_p$ be the transitive closure of $\approx_p$ (so $\sim_p$ is an equivalence relation on $\rng(\bar{x})$).
\end{enumerate}
}\end{defin}

\begin{lem}\label{basic property of the equivalence relation}
Let $\bar{x} = (x_1, \ldots, x_k)$ and let $p(\bar{x})$ be a complete $(\tau, \lambda)$-neighbourhood type for some $\lambda \in \mbbN$.
For all $i, j = 1, \ldots, k$,
$x_i \sim_p x_j$ if and only if there is $l \in \mbbN$ such that $p(\bar{x}) \models \dist(x_i, x_j) \leq l$.
\end{lem}

\noindent
{\bf Proof.}
Let $p(\bar{x})$ be a complete $(\tau, \lambda)$-neighbourhood type.
Then there is a finite $\tau$-structure $\mcB$ and $\bar{b} \in B^{|\bar{x}|}$ such that $\mcB \models p(\bar{b})$.
Suppose that $x_i \sim_p x_j$.
It follows straightforwardly from the definition of $\sim_p$ that there is $l$ depending only on $p$ such that $\dist(b_i, b_j) \leq l$.
As we only assumed that $\mcB$ is a finite $\tau$-structure and $\mcB \models p(\bar{b})$ we get $p(\bar{x}) \models \dist(x_i, x_j) \leq l$
(according to Notation~\ref{Notation for consequence in finite}).

Now suppose that $x_i \not\sim_p x_j$. 
Let $\mcB'$ be the substructure of $\mcB$ with domain $N_\lambda^\mcB(\bar{b})$.
Then $\mcB' \models p(\bar{b})$ and $b_i$ and $b_j$ are in different connected components of the Gaifman graph of $\mcB'$.
It follows that, for all $l$, $\dist_{\mcB'}(b_i, b_j) > l$.
So for all $l$, $p(\bar{x}) \not\models \dist(x_i, x_j) \leq l$.
\hfill $\square$

\noindent
{\bf \em For the rest of the article 
(except in Section~\ref{Examples of sequences of base structures} where we give examples) 
we will assume the following:}

\begin{assump}\label{properties of the base structures} {\bf (Properties of the base structures)} {\rm
Let $\mbB = (\mcB_n : n \in \mbbN^+)$ be a sequence of finite $\tau$-structures such that the following hold:
\begin{enumerate}
\item $\lim_{n\to\infty} |B_n| = \infty$.
\item There is $\Delta \in \mbbN$ such that, for every $n$, the degree of $\mcB_n$ is at most $\Delta$.
\item There is a polynomial $P(x)$ such that for all $n$, $|B_n| \leq P(n)$.
\item Suppose that $\lambda \in \mbbN$, $\bar{x} = (x_1, \ldots, x_k)$ is a sequence of distinct variables, 
and that $p(\bar{x})$ is a complete $(\tau, \lambda)$-neighbourhood type such that for all $i, j \in \{1, \ldots, k\}$,
$x_i \sim_p x_j$. Then
\begin{enumerate}
\item either there is $m_p \in \mbbN$ such that for all $n$, $|p(\mcB_n)| \leq m_p$, or

\item there is $f_p : \mbbN \to \mbbR$ such that for all $\alpha \in \mbbR$, $\lim_{n\to\infty} (f_p(n) - \alpha \ln(n)) = \infty$,
and for all sufficiently large $n$, $|p(\mcB_n)| \geq f_p(n)$.
\end{enumerate}
\end{enumerate}
}\end{assump}

\noindent
The role of part~(4)(b) of Assumption~\ref{properties of the base structures} is that it implies that for every
polynomial $P$ and $\alpha > 0$, $P(n) e^{-\alpha f_p(n)}$ tends to zero exponentially fast, as stated by the following lemma:

\begin{lem}\label{consequence of 4b} 
Let $k \in \mbbN^+$, $f : \mbbN \to \mbbR$ and suppose that for all $\alpha \in \mbbR$, $\lim_{n\to\infty} (f(n) - \alpha \ln(n)) = \infty$.
For every $\alpha > 0$, if $n$ is sufficiently large then
$n^k e^{-\alpha f(n)} \leq e^{-\frac{1}{2}\alpha f(n)}$.
\end{lem}

\noindent
{\bf Proof.}
Argue just as in the proof of Lemma~6.5 in \cite{KT} with $f$ in place of `$g_3$' in that proof.
\hfill $\square$

\begin{termin}\label{remark on skipping with respect to}{\rm
Since many definitions that will follow depend on the {\em base sequence} $\mbB = (\mcB_n : n \in \mbbN^+)$
of $\tau$-structures, they should strictly speaking be tagged by ``with respect to $\mbB$''.
But since we have fixed the base sequence $\mbB$ for the rest of the article we omit the phrase ``with respect to $\mbB$''.
}\end{termin}

\begin{defin}\label{definition of cofinally satisfiable}{\rm
We call $\varphi(\bar{x}) \in PLA^*(\tau)$ {\em cofinally satisfiable} if for every $n_0 \in \mbbN$
there is $n \in \mbbN$ such that $n \geq n_0$ and $\mcB_n(\varphi(\bar{b})) = 1$ for some $\bar{b} \in (B_n)^{|\bar{x}|}$.
}\end{defin}

\begin{defin}\label{definition of bounded} {\rm
Let $\varphi(\bar{x}, \bar{y}), \psi(\bar{x}) \in  PLA^*(\tau)$.
\begin{enumerate}
\item For $m \in \mbbN$, we say that $\varphi(\bar{x}, \bar{y})$ is {\bf \em $(\bar{y}, m)$-bounded}
if for all $n$ and all $\bar{a} \in (B_n)^{|\bar{a}|}$, $|\varphi(\bar{a}, \mcB_n)| \leq m$.
In the special case when $\bar{x}$ is empty we may just say that $\varphi(\bar{y})$ is {\em $m$-bounded}.

\item We call $\varphi(\bar{x}, \bar{y})$ {\bf \em $\bar{y}$-bounded} if it is $(\bar{y}, m)$-bounded for some $m \in \mbbN$.
In the special case when $\bar{x}$ is empty we may just say that $\varphi(\bar{y})$ is {\bf \em bounded}.
\end{enumerate}
}\end{defin}

\noindent
We will have to make a distinction between elements with (for some $\lambda \in \mbbN$)
a $\lambda$-neighbourhood that has only a bounded number of isomorphic copies (in every $\mcB_n$)
and elements with a $\lambda$-neighbourhood that has ``arbitrarily many'' isomorphic copies.
This motivates the next definition.

\begin{defin}\label{definition of rare element} {\rm
Let $n \in \mbbN^+$, $\lambda \in \mbbN$, and $a \in B_n$. We say that $a$ is {\bf \em $\lambda$-rare}
if there is a bounded complete $(\tau, \lambda)$-neighbourhood type 
$p(x)$ such that $\mcB_n \models p(a)$. 
}\end{defin}

\begin{lem}\label{definability of m-rareness}
Let $\lambda \in \mbbN$.\\
(i) There are $\varphi(x) \in FO(\tau)$ and a 0/1-valued formula
$\psi(x) \in PLA^*(\tau)$ such that for all $n$ and all $a \in B_n$,
$a$ is $\lambda$-rare if and only if $\mcB_n \models \varphi(a)$ if and only if $\mcB_n(\psi(a)) = 1$.\\
(ii) There is $m \in \mbbN$ such that, for all $n \in \mbbN$, $\mcB_n$ has at most $m$ $\lambda$-rare elements.
\end{lem}

\noindent
{\bf Proof.}
Up to equivalence, there are only finitely many bounded complete $(\tau, \lambda)$-neigh\-bour\-hood types
$p_1(x), \ldots, p_s(x)$ in the variable $x$, so for all $n$ and $a \in B_n$, $a$ is $\lambda$-rare if and only if
$\mcB_n \models p_1(a) \vee \ldots \vee p_s(a)$.
Then recall that every first-order formula is equivalent (in finite structures) to a $PLA^*$-formula.
For the second part, note that since each $p_i(x)$ is bounded there is $m_i$ such that $|p_i(\mcB_n)| \leq m_i$ for all $n$.
\hfill $\square$
\\

\noindent
The $\lambda$-neighbourhood of a sequence of elements need not contain all (or any) of the $\lambda$-rare elements.
Therefore we define the {\em $\lambda$-closure} of a sequence of elements
to be the $\lambda$-neighbourhood of that sequence augmented with all $\lambda$-rare elements.

\begin{defin}\label{definition of closure} {\rm
Let $\lambda \in \mbbN$, let $n \in \mbbN^+$, and let $\bar{b} = (b_1, \ldots, b_k) \in (B_n)^k$.
Then the {\bf \em $\lambda$-closure of $\bar{b}$ (in $\mcB_n$)} is the set
\begin{align*}
&C_\lambda^{\mcB_n}(\bar{b}) = N_\lambda^{\mcB_n}(X)  \ \ \ \text{ where } \\ 
&X = \{b_1, \ldots, b_k\} \cup \{c \in B_n : \text{ $c$ is $\lambda$-rare} \}.
\end{align*}
}\end{defin}

\begin{rem}\label{remark on size of closures}{\rm
It follows from Assumption~\ref{properties of the base structures}
and Lemma~\ref{definability of m-rareness}
that for all $k \in \mbbN^+$ and $\lambda \in \mbbN$ there is a constant $m_{k, \lambda}$ such that 
for all $n$ and $\bar{a} \in (B_n)^k$, 
$|N_\lambda^{\mcB_n}(\bar{a})|, |C_\lambda^{\mcB_n}(\bar{a})| \leq m_{k, \lambda}$.
}\end{rem}

\noindent
In the next definition the idea is that a formula $p(\bar{x})$ is a ``complete $(\sigma, \lambda)$-closure type
if $p(\bar{x})$ expresses that the following holds:
If $\mcA_1$ and $\mcA_2$ are finite $\sigma$-structures, $\mcB_1 = \mcA_1 \uhrc \tau$, $\mcB_2 = \mcA_2 \uhrc \tau$,
$\mcA_1 \models p(\bar{a}_1)$ and $\mcA_2 \models p(\bar{a}_2)$, 
then there is an isomorphism $f$ from $\mcA_1 \uhrc C_\lambda^{\mcB_1}(\bar{a}_1)$ to 
$\mcA_2 \uhrc C_\lambda^{\mcB_2}(\bar{a}_2)$ such that $f(\bar{a}_1) = \bar{a}_2$.

\begin{defin}\label{definition of closure type}  {\rm
Let $\lambda \in \mbbN$.
A {\bf \em $(\sigma, \lambda)$-closure type} in the variables $\bar{x}$ is a formula of the following form, where 
$\varphi_\lambda(y)$ is a formula which expresses that ``$y$ is $\lambda$-rare'' and
$p(\bar{x}, y_1, \ldots, y_k)$ is a $(\sigma, \lambda)$-neighbourhood type:
\begin{align*}
&\exists y_1, \ldots, y_k \bigg( 
\bigwedge_{i=1}^k \varphi_\lambda(y_i) \ \wedge \ 
\forall u \Big(\varphi_\lambda(u) \rightarrow \bigvee_{i=1}^k u = y_i \Big) \ \wedge \  
p(\bar{x}, y_1, \ldots, y_k) \bigg).
\end{align*}
If, in addition, $p$ is a {\em complete} $(\sigma, \lambda)$-neighbourhood type, then the above formula is
called a {\bf \em complete $(\sigma, \lambda)$-closure type}.
}\end{defin}

\begin{rem}\label{remark about empty x in closure types}{\rm
We allow the sequence $\bar{x}$ in the definition above of a $(\sigma, \lambda)$-closure type to be empty.
In this case the $(\sigma, \lambda)$-closure type is a sentence which completely describes the $\tau$-structure
of the $\lambda$-neighbourhood of the $\lambda$-rare elements 
and partially or completely describes the $\sigma$-structure 
of the $\lambda$-neighbourhood of the
$\lambda$-rare elements.
}\end{rem}

\begin{defin}\label{definition of dimension for closure types}{\rm
Let $\lambda \in \mbbN$ and let $p(\bar{x})$ be a complete $(\tau, \lambda)$-closure type
where $\bar{x} = (x_1, \ldots, x_k)$.
Then the relations $\approx_p$ and $\sim_p$
are defined exactly as for complete $(\tau, \lambda)$-neighbourhood types in
Definition~\ref{definition of dimension}.
}\end{defin}

\begin{lem}\label{closure types and cases}
Let $\lambda \in \mbbN$, $\bar{x} = (x_1, \ldots, x_k)$ and let $p(\bar{x})$ be a complete $(\tau, \lambda)$-closure type.\\
(i) For all $i, j = 1, \ldots, k$, either $p(\bar{x}) \models \dist(x_i, x_j) \leq 2\lambda$ or 
$p(\bar{x}) \models \dist(x_i, x_j) > 2\lambda$.\\
(ii) Let $\varphi_\lambda(z)$ be a formula which, in all $\mcB_n$, expresses that ``$z$ is $\lambda$-rare''.
For  $i = 1, \ldots, k$, either $p(\bar{x}) \models \exists z \big(\varphi_\lambda(z) \wedge \dist(z, x_i) \leq 2\lambda\big)$
or $p(\bar{x}) \models \forall z \big(\varphi_\lambda(z) \rightarrow \dist(z, x_i) > 2\lambda\big)$.\\
(iii) For all $i, j = 1, \ldots, k$, 
$x_i \sim_p x_j$ if and only if there is $l \in \mbbN$ such that $p(\bar{x}) \models \dist(x_i, x_j) \leq l$.
\end{lem}

\noindent
{\bf Proof.}
All three parts follows straightforwardly from the definitions of closure types and $\sim_p$ (for closure types),
and Lemmas~\ref{neighbourhood types and distance}
and~\ref{basic property of the equivalence relation}.
\hfill $\square$

\medskip

\noindent
The next definition describes three kinds of ``restrictions'' of $(\sigma, \lambda)$-closure types, 
or $(\sigma, \lambda)$-neighbourhood types,
that will play important technical roles later.
Again, note that the definitions make sense for $\tau$ in place of $\sigma$ since we allow that $\sigma = \tau$.

\begin{defin}\label{restriction of closure type} {\rm
Let, for some $\lambda \in \mbbN$ 
$p(\bar{x})$ be a complete $(\sigma, \lambda)$-neighbourhood type (or a 
complete $(\sigma, \lambda)$-closure type).
\begin{enumerate}
\item If $\bar{y}$ is a subsequence of $\bar{x}$, then the {\bf \em restriction of $p(\bar{x})$ to $\bar{y}$},
denoted $p \uhrc \bar{y}$, is a complete $(\sigma, \lambda)$-neighbourhood type 
(or a complete $(\sigma, \lambda)$-closure type) $q(\bar{y})$ in the variables $\bar{y}$
which is consistent with $p(\bar{x})$. (We say ``the'' restriction of $p(\bar{x})$ to $\bar{y}$ because all such restrictions
are equivalent.)

\item If $\tau \subseteq \sigma' \subseteq \sigma$, then the 
{\bf \em restriction of $p(\bar{x})$ to $\sigma'$},
denoted $p \uhrc \sigma'$, is a complete $(\sigma', \lambda)$-neighbourhood type 
(or a complete $(\sigma', \lambda)$-closure type) $p'(\bar{x})$ in the variables $\bar{x}$
which is consistent with $p(\bar{x})$. (Again we say ``the'' restriction of $p(\bar{x})$ to $\sigma'$ because all such restrictions
are equivalent.)

\item If $\gamma \in \mbbN$ and $\gamma \leq \lambda$, then the {\bf \em restriction of $p(\bar{x})$ to $\gamma$},
denoted $p \uhrc \gamma$,
is a complete $(\sigma, \gamma)$-neighbourhood type 
(or a  complete $(\sigma, \gamma)$-closure type)
type $p'(\bar{x})$ which is
consistent with $p(\bar{x})$. (Yet again we say ``the'' restriction of $p(\bar{x})$ to $\gamma$ because all such restrictions
are equivalent.)
\end{enumerate}
}\end{defin}

\begin{rem}\label{remark on restriction sim}{\rm
Let $p(\bar{x})$ be (for some $\lambda \in \mbbN$) a complete $(\tau, \lambda)$-neighbourhood type or
a complete $(\tau, \lambda)$-closure type, let $\bar{y}$ be a subsequence of $\bar{x}$ and let
$q(\bar{y}) = p \uhrc \bar{y}$ be the restriction of $p$ to $\bar{y}$. 
It follows straightforwardly from the definition of $\sim_p$ and the definition
of neighbourhood type, and closure type, that $\sim_q$ and $\sim_p$ coincide on $\rng(\bar{y})$.
This will be used from time to time without specific explanation or reference.
}\end{rem}

\noindent
We now define two conditions about closure types and neighbourhood types that are stronger than
being unbounded.

\begin{defin}\label{definition of unbounded formula}{\rm  $\text{   }$ 
\begin{enumerate}
\item A formula $\varphi(\bar{x}, \bar{y}) \in PLA^*(\sigma)$ is {\bf \em uniformly $\bar{y}$-unbounded} if 
$\varphi(\bar{x}, \bar{y})$ is cofinally satisfiable and
there is $f_\varphi : \mbbN \to \mbbR$
such that 
\begin{enumerate}
\item for all $\alpha > 0$ $\lim_{n\to\infty} (f_\varphi(n) - \alpha\ln(n)) = \infty$ and 
\item for all $n$ and all $\bar{a} \in (B_n)^{|\bar{x}|}$, 
if $\varphi(\bar{a}, \mcB_n) \neq \es$ then
$|\varphi(\bar{a}, \mcB_n)| \geq f_\varphi(n)$.
\end{enumerate}
In the special case when $\bar{x}$ is empty we say that $\varphi(\bar{y})$ is {\bf \em uniformly unbounded}.

\item Let $p(\bar{x}, \bar{y})$ be a complete
$(\tau, \lambda)$-closure type or a complete $(\tau, \lambda)$-neighbourhood type where $\lambda \in \mbbN$.
We say that $p(\bar{x}, \bar{y})$ is {\bf \em strongly $\bar{y}$-unbounded} if, for every subsequence
$\bar{y}'$ of $\bar{y}$, $p \uhrc \bar{x}\bar{y}'$ is uniformly $\bar{y}'$-unbounded.
In the special case when $\bar{x}$ is empty we say that $\varphi(\bar{y})$ is {\bf \em strongly unbounded}.

\item Let $p(\bar{x}, \bar{y})$ be a complete $(\sigma, \lambda)$-closure type 
or a complete $(\sigma, \lambda)$-neighbourhood type where $\lambda \in \mbbN$.
We say that $p$ is 
{\bf \em uniformly $\bar{y}$-unbounded}, respectively
{\bf \em strongly $\bar{y}$-unbounded}, 
if $p \uhrc \tau$
is uniformly $\bar{y}$-unbounded, respectively
strongly $\bar{y}$-unbounded.
\end{enumerate}
}\end{defin}

\noindent
Later, in Lemma~\ref{not bounded implies uniformly unbounded},
we will see that if $p(\bar{x}, \bar{y})$ is a complete $(\tau, \lambda)$-closure type then either $p$ is $\bar{y}$-bounded
or uniformly $\bar{y}$-unbounded. 
Assumption~\ref{properties of the base structures} stipulates the same thing for $(\tau, \lambda)$-neighbourhood types $p$
in which all free variables are in the same $\sim_p$-class.

\begin{lem}\label{strong unboundedness implies large distance}
Let $\lambda \in \mbbN$ and suppose that $p(\bar{x}, \bar{y})$ is a 
complete $(\tau, \lambda)$-closure type or a  complete $(\tau, \lambda)$-neighbourhood type.
Also, let $\varphi_\lambda(z)$ be a formula which expresses that ``$z$ is $\lambda$-rare'' in every $\mcB_n$.
\\
(i) If, for every subsequence $\bar{y}'$ of $\bar{y}$, $p \uhrc \bar{x}\bar{y}'$ is {\em not} $\bar{y}'$-bounded,
then 
\begin{equation}\label{p implies y far from x and from rare elements}
p(\bar{x}, \bar{y}) \models \dist(\bar{x}, \bar{y}) > 2\lambda \ \wedge \ 
\forall z \big(\varphi_\lambda(z) \rightarrow \dist(z, \bar{y}) > 2\lambda \big).
\end{equation}
(ii) If $p$ is strongly $\bar{y}$-unbounded then~(\ref{p implies y far from x and from rare elements}) holds.
\end{lem}

\noindent
{\bf Proof.}
(i) Suppose that~(\ref{p implies y far from x and from rare elements}) does not hold.
By Lemma~\ref{closure types and cases} we have
\[
p(\bar{x}, \bar{y}) \models \dist(\bar{x}, \bar{y}) \leq 2\lambda \ \text{ or } \
\exists z \big(\varphi_\lambda(z) \wedge \dist(z, \bar{y}) \leq 2\lambda \big).
\]
In the first case
there are $x_i \in \rng(\bar{x})$ and $y_i \in \rng(\bar{y})$ such that $p(\bar{x}, \bar{y}) \models \dist(x_i, y_j) \leq 2\lambda$
and since the $2\lambda$-neighbourhood of an element has cardinality at most $\Delta^{2\lambda}$ it follows that 
$p \uhrc \bar{x}y_i$ is $y_i$-bounded, which contradicts the assumption.
In the second case there is $y_j \in \rng(\bar{y})$ such that
$p(\bar{x}, \bar{y}) \models  \exists z \big(\varphi_\lambda(z) \wedge \dist(z, y_j) \leq 2\lambda \big)$.
By the definition of $\lambda$-rare element there is some $m \in \mbbN$, depending only on 
Assumption~\ref{properties of the base structures}, such that every $\mcB_n$ contains at most $m$ $\lambda$-rare elements,
and the cardinality of the neighbourhood of every such element has cardinality at most $\Delta^{2\lambda}$.
Hence $p \uhrc \bar{x}y_j$ is $y_j$-bounded, which contradicts the assumption.

Part~(ii) follows from part~(i) and the definition of `strongly $\bar{y}$-unbounded'.
\hfill $\square$

\begin{lem}\label{elementary properties of closure types}
Suppose that $p(\bar{x}, \bar{y})$ is a complete $(\tau, \lambda)$-closure type 
or a complete $(\sigma, \lambda)$-neighbourhood type for some $\lambda \in \mbbN$.\\
(i) If $\gamma \leq \lambda$ and $p \uhrc \gamma$ is $\bar{y}$-bounded, then $p$ is $\bar{y}$-bounded.\\
(ii) If $p$ is strongly $\bar{y}$-unbounded and $\bar{z}$ is a nonempty subsequence of $\bar{y}$,
then $p$ is strongly $\bar{z}$-unbounded
and $p \uhrc \bar{z}$ is strongly unbounded.\\
(iii) If $\gamma \leq \lambda$ and $p$ is uniformly (respectively strongly) $\bar{y}$-unbounded,
then $p \uhrc \gamma$ is uniformly (respectively strongly) $\bar{y}$-unbounded.
\end{lem}

\noindent
We will usually reason about $(\tau, \lambda)$-closure types (for some $\lambda$),
because they talk about $\lambda$-rare elements which we cannot ignore in the present context.
But in some arguments we want to ``reduce'' our reasoning about a $(\tau, \lambda)$-closure type to reasoning about a
related $(\tau, \lambda)$-neighbourhood since the later is a simpler concept.
The following lemma gives the connection:

\begin{lem}\label{connection between closure types and neighbourhood types}
Let $p(\bar{x})$ be a cofinally satisfiable complete $(\tau, \lambda)$-closure type for some $\lambda \in \mbbN$,
and let $\varphi_\lambda(z)$ be a formula that expresses (in all $\mcB_n$) that ``$z$ is $\lambda$-rare''.
Then there are a sequence $\bar{z} = (z_1, \ldots, z_k)$ of variables
and a complete $(\tau, \lambda)$-neighbourhood type $p^+(\bar{z}, \bar{x})$ such that 
\begin{enumerate}
\item $p(\bar{x})$ is equivalent to 
$\exists \bar{z} \Big(p^+(\bar{z}, \bar{x}) \wedge \forall u \Big(\varphi_\lambda(u) \rightarrow \bigvee_{i=1}^k u = z_i \Big)\Big)$,
\item $p^+\uhrc \bar{z}$ is bounded and, for all $i = 1, \ldots, k$, $p^+ \uhrc z_i \models \varphi_\lambda(z_i)$, and 
\item if $\bar{x}'$ is a subsequence of $\bar{x}$ then $p \uhrc \bar{x}'$ is equivalent to \\
$\exists \bar{z} \Big(p^+ \uhrc \bar{z}\bar{x}' \wedge 
\forall u \Big(\varphi_\lambda(u) \rightarrow \bigvee_{i=1}^k u = z_i \Big)\Big)$.
\end{enumerate}
\end{lem}

\noindent
{\bf Proof.}
Let $p(\bar{x})$ be a complete $(\tau, \lambda)$-closure type for some $\lambda \in \mbbN$.
By Definition~\ref{definition of closure type} 
of a closure type there is a sequence of variables $\bar{z} = (z_1, \ldots, z_k)$
and a complete $(\tau, \lambda)$-{\em neighbourhood} type $p^+(\bar{z}, \bar{x})$
such that 
\begin{equation}\label{equivalence of p to an existential formula}
p(\bar{x}) \ \text{ is equivalent to }  \
\exists \bar{z}\Big(\bigwedge_{i=1}^k \varphi_\lambda(z_i) \ 
\wedge \ \forall u \Big(\varphi_\lambda(u) \rightarrow \bigvee_{i=1}^k u = z_i\Big) \ 
\wedge \ p^+(\bar{z}, \bar{x}) \Big)
\end{equation}
where $\varphi_\lambda(z)$ is a formula which expresses that ``$z$ is $\lambda$-rare''.
By the proof of 
Lemma~\ref{definability of m-rareness}
we can let $\varphi_\lambda(z)$ be $p_1(z) \vee \ldots \vee p_m(z)$ where 
the sequence $p_i(z)$, $i = 1, \ldots, m$,
enumerates all, up to equivalence, bounded complete $(\tau, \lambda)$-neighbourhood types in the variable $z$.
As $p(\bar{x})$ is cofinally satisfiable (hence consistent) and $p^+(\bar{z}, \bar{x})$ is a 
complete $(\tau, \lambda)$-neighbourhood type 
it follows that for each $j = 1, \ldots, k$ there is $i_j \in \{1, \ldots, m\}$ such that  $p^+(\bar{z}, \bar{x}) \models p_{i_j}(z_j)$
and hence $p^+ \uhrc z_j$ is equivalent to $p_{i_j}(z_j)$, so $p^+ \uhrc z_j \models \varphi_\lambda(z_j)$.
Therefore $p^+(\bar{z})$ is bounded.
This proves part~(2).
From $p^+ \uhrc z_j \models \varphi_\lambda(z_j)$ we have $p^+(\bar{z}, \bar{x}) \models \varphi_\lambda(z_j)$
for all $j = 1, \ldots, k$, so part~(1) now follows from~(\ref{equivalence of p to an existential formula}).

Part~(3) holds since (by the definition of $(\tau, \lambda)$-closure type)
if $\bar{x}'$ is a subsequence of $\bar{x}$ then~(\ref{equivalence of p to an existential formula})
holds if $p$ is replaced by $p \uhrc \bar{x}'$ and $p^+$ is replaced by $p^+ \uhrc \bar{z}\bar{x}'$.
\hfill $\square$
\\

\noindent
The next lemma intuitively says that under some circumstances a $(\tau, \lambda)$-closure type can
be decomposed into two independent parts.

\begin{lem}\label{splitting a closure type}
Let $p(\bar{x}, \bar{y})$ be a cofinally satisfiable complete $(\tau, \lambda)$-closure type for some $\lambda \in \mbbN$
and let $\varphi_\lambda(z)$ express (in every $\mcB_n$) that ``$z$ is $\lambda$-rare''.
Suppose that
\[
p(\bar{x}, \bar{y}) \models \dist(\bar{x}, \bar{y}) > 2\lambda \ \wedge \
\forall z \big(\varphi_\lambda(z) \rightarrow \dist(z, \bar{y}) > 2\lambda \big).
\]
Let $p_1(\bar{x}) = p \uhrc \bar{x}$ and $p_2(\bar{y}) = p \uhrc \bar{y}$.
Then there is a complete $(\tau, \lambda)$-neighbourhood type $p^+_2(\bar{y})$
such that $p_2(\bar{y}) \models p^+_2(\bar{y})$ and, 
for all $n \in \mbbN^+$, for all $\bar{a} \in (B_n)^{|\bar{x}|}$ and $\bar{b} \in (B_n)^{|\bar{y}|}$,
\begin{enumerate}
\item[] $\mcB_n \models p(\bar{a}, \bar{b})$ if and only if $\mcB_n \models p_1(\bar{a}) \wedge p^+_2(\bar{b})$ and
$\dist(\bar{c}\bar{a}, \bar{b}) > 2\lambda$, where $\bar{c}$ enumerates all $\lambda$-rare elements in $\mcB_n$.
\end{enumerate}
(As $p_2 = p \uhrc \bar{y}$ and $p_2 \models p^+_2$ the above also holds if we replave $p^+_2$ by $p_2$.)
We allow $\bar{x}$ to be empty in which case $p_1$ is a sentence which completely describes the $\tau$-structure
of the $\lambda$-neighbourhood of the $\lambda$-rare elements.
\end{lem}

\noindent
{\bf Proof.}
Let $p(\bar{x}, \bar{y})$ be as assumed in the lemma.
By Lemma~\ref{connection between closure types and neighbourhood types},
there is $\bar{z} = (z_1, \ldots, z_k)$ and a complete $(\tau, \lambda)$-{\em neighbourhood} type 
$p^+(\bar{z}, \bar{x}, \bar{y})$ such that 
\begin{align}\label{p and exists x p+ are equivalent}
&\text{$p(\bar{x}, \bar{y})$ is equivalent to }
\exists \bar{z} \Big(p^+(\bar{z}, \bar{x}, \bar{y}) \wedge 
\forall u \Big(\varphi_\lambda(u) \rightarrow \bigvee_{i=1}^k u = z_i \Big)\Big),\\
&\text{for all $i = 1, \ldots, k$, $p^+\uhrc z_i \models \varphi_\lambda(z_i)$, and $p^+\uhrc \bar{z}$ is bounded.} \nonumber
\end{align}
Let $p_1(\bar{x}) = p \uhrc \bar{x}$, $p_2(\bar{y}) = p \uhrc \bar{y}$,
$p^+_1(\bar{z}, \bar{x}) = p^+ \uhrc \bar{z}\bar{x}$,  and
$p^+_2(\bar{y}) = p^+ \uhrc \bar{y}$.
By part~(3) of 
Lemma~\ref{connection between closure types and neighbourhood types},
\begin{equation}\label{p+1 equivalent to something}
p_1(\bar{x}) \text{ is equivalent to } \ 
\exists \bar{z} \Big( p^+_1(\bar{z}, \bar{x}) \wedge  
\forall u \Big(\varphi_\lambda(u) \rightarrow \bigvee_{i=1}^k u = z_i \Big)\Big).
\end{equation}
By definitions of the involved formulas we have $p_2(\bar{y}) \models p^+_2(\bar{y})$.

Suppose that $\mcB_n \models p(\bar{a}, \bar{b})$ and let $\bar{c}$ enumerate all $\lambda$-rare elements in $B_n$.
By the assumptions of the lemma, we have
$\dist(\bar{c}\bar{a}, \bar{b}) > 2\lambda$.
By the choices of $p_1$ and $p_2^+$ and by~(\ref{p and exists x p+ are equivalent})
we also have $\mcB_n \models p_1(\bar{a}) \wedge p_2^+(\bar{b})$.

Now suppose that 
$\mcB_n \models p_1(\bar{a}) \wedge p^+_2(\bar{b})$ and $\dist(\bar{c}\bar{a}, \bar{b}) > 2\lambda$
where $\bar{c}$ enumerates all $\lambda$-rare elements.
By~(\ref{p+1 equivalent to something}) and by reordering $\bar{c}$ if necessary, we get
$\mcB_n \models p^+_1(\bar{c}, \bar{a})$. 
Now Lemma~\ref{splitting a neighbourhood type} implies that
$\mcB_n \models p^+(\bar{c}, \bar{a}, \bar{b})$ so
by~(\ref{p and exists x p+ are equivalent}) we get $\mcB_n \models p(\bar{a}, \bar{b})$.
\hfill $\square$

\medskip
\noindent
Recall the notion of {\em $L$-basic formula} from 
Definition~\ref{L-basic subformula}, where $L \subseteq PLA^*(\sigma)$.
In this article we will work with a special case of this notion, namely the notion of a
{\em $(\sigma, \lambda)$-basic formula} defined below. 
It will play an essential role 
from Section~\ref{convergence and balance} and onwards.

\begin{defin}\label{Definition of basic formula} {\rm
Let $\lambda \in \mbbN$. By a {\bf \em $(\sigma, \lambda)$-basic formula} 
(in the variables $\bar{x}$) we mean a formula of the form
$\bigwedge_{i = 1}^k \big(\varphi_i(\bar{x}) \to c_i\big)$ where, for each $i$, $c_i \in [0, 1]$, $\lambda_i \in \mbbN$,
$\lambda_i \leq \lambda$, and $\varphi_i(\bar{x})$ is a complete $(\sigma, \lambda_i)$-closure type.
We allow the possibility that $\bar{x}$ is empty in which case each $\varphi_i$ is taken to be the formula $\top$ 
(which is another name for the formula `1').
}\end{defin}

\begin{rem}\label{remark on basic formulas}{\rm
It is straightforward to see that if $\varphi(\bar{x})$ is a $(\sigma, \lambda)$-basic formula, then there is
a $(\sigma, \lambda)$-basic formula $\bigwedge_{i = 1}^k \big(\varphi_i(\bar{x}) \to c_i\big)$
which is equivalent to $\varphi(\bar{x})$ and such that $\varphi_1(\bar{x}), \ldots, \varphi_k(\bar{x})$
lists, up to equivalence, all complete $(\sigma, \lambda)$-closure types.
}\end{rem}

\noindent
The following result tells that connectives ``preserve'', up to equivalence, the property of being $(\sigma, \lambda)$-basic,
and it follows that every aggregation-free formula in $PLA^*(\sigma)$ is equivalent to a $(\sigma, 0)$-basic formula.

\begin{lem}\label{connectives and basic formulas}
(i) Suppose that $\varphi_1(\bar{x}), \ldots, \varphi_k(\bar{x})$ are $(\sigma, \lambda)$-basic formulas and that
$\msfC : [0, 1]^k \to [0, 1]$. Then the formula $\msfC(\varphi_1(\bar{x}), \ldots, \varphi_k(\bar{x}))$
is equivalent to a $(\sigma, \lambda)$-basic formula.
In particular, if $\varphi_1(\bar{x}), \ldots, \varphi_k(\bar{x})$ are $(\sigma, 0)$-basic formulas, 
then $\msfC(\varphi_1(\bar{x}), \ldots, \varphi_k(\bar{x}))$
is equivalent to a $(\sigma, 0)$-basic formula.\\
(ii) If $\varphi(\bar{x}) \in PLA^*(\sigma)$ is aggregation-free then it is equivalent 
to a $(\sigma, 0)$-basic formula.
\end{lem}

\noindent
{\bf Proof.}
(i) Let $\varphi(\bar{x})$ denote the formula $\msfC(\varphi_1(\bar{x}), \ldots, \varphi_k(\bar{x}))$ where
each $\varphi_i(\bar{x})$ is $(\sigma, \lambda)$-basic.
Let $q_1(\bar{x}), \ldots, q_m(\bar{x})$ enumerate, up to logical equivalence, all complete $(\sigma, \lambda)$-closure types.

Let $\mcA \in \mbW_n$ for some $n$ and $\bar{a} \in (B_n)^{|\bar{x}|}$.
For each $i$ the value $\mcA(\varphi_i(\bar{a}))$ depends only on which $q_j(\bar{x})$ the sequence $\bar{a}$ satisfies.
So let $c_{i, j} = \mcA(\varphi_i(\bar{a}))$ if $\mcA \models q_j(\bar{a})$.
Then let $d_j = \msfC(c_{1, j}, \ldots, c_{k, j})$ for $j = 1, \ldots, m$.
Now $\varphi(\bar{x})$ is equivalent in to 
the $(\sigma, \lambda)$-basic formula $\bigwedge_{j = 1}^m (q_j(\bar{x}) \to d_j)$.

(ii) Let $\varphi(\bar{x}) \in PLA^*(\sigma)$ be aggregation-free.
The proof proceeds by induction on the number of connectives in $\varphi$.
If the number of connectives is 0 then $\varphi(\bar{x})$ can be a constant from $[0, 1]$, or
it can have the form $R(\bar{x}')$ for some $R \in \sigma$ and subsequence $\bar{x}'$ of $\bar{x}$, or it can have the 
form $u = v$ for some $u, v \in \rng(\bar{x})$. 
We leave it to the reader to verify that in each case $\varphi(\bar{x})$ is equivalent 
to a $(\sigma, 0)$-basic formula.
The inductive step follows from part~(i) of this lemma (in the special case when $\lambda = 0$).
\hfill $\square$

\section{Examples of sequences of base structures}\label{Examples of sequences of base structures}

\noindent
Throughout this article we assume that the sequence of base structures $\mbB = (\mcB_n : n \in \mbbN^+)$
satisfies Assumption~\ref{properties of the base structures}, where each $\mcB_n$ is a $\tau$-structure and
$\tau$ is a finite relational signature.
In this section we give examples of $\mbB$ that satisfy Assumption~\ref{properties of the base structures}.

\begin{exam}\label{example with empty tau} {\bf (Sets without structure)} {\rm
Let $\tau$ be empty and, for each $n \in \mbbN^+$, let $\mcB_n$ be the unique $\tau$-structure with 
domain $B_n = \{1, \ldots, n\}$.
Then all conditions of Assumption~\ref{properties of the base structures} hold and in part~(4) of it we can take $m_p = 0$ and
$f_p(n) = n$.
We also have $N_\lambda^{\mcB_n}(a_1, \ldots, a_k) = C_\lambda^{\mcB_n}(a_1, \ldots, a_k) = 
\{a_1, \ldots, a_k\}$ for all $\lambda \in \mbbN$, $n \in \mbbN^+$ and $a_1, \ldots, a_k \in B_n$.
It will follow (from the rest of the article) that the framework of this article and its results generalize the frameworks and results
in previous work of the author and Weitkämper \cite{KW1, KW2}. 
Since first-order logic (evaluated on finite structures) can,
according to Remark~\ref{FO is expressible in PLA*}, 
be seen as a sublogic of $PLA^*$
it will also follow that the classical result about almost sure elimination of 
(first-order) quantifiers \cite{Gle} and the classical zero-one law \cite{Gle, Fag}
are special cases of results proved here.
}\end{exam}

\begin{exam}\label{example with unary relations} {\bf (Unary relations, denoting rare objects)} {\rm
Let $\tau = \{P_1, \ldots, P_s\}$ where all relation symbols in $\tau$ are unary
and let $m \in \mbbN$ be fixed.
For every $n \in \mbbN^+$, let $\mcB_n$ be the $\tau$-structure with domain $B_n = \{1, \ldots, n\}$
and where, for  all $i = 1, \ldots, s$, $(P_i)^{\mcB_n}$ is a subset of $B_n$ of cardinality at most $m$.
Then, for any two different $a, b \in B_n$ we have $\dist(a, b) = \infty$ and it follows that for
any $k \in \mbbN^+$ and $a_1, \ldots, a_k \in B_n$, $N_\lambda^{\mcB_n}(a_1, \ldots, a_k) = \{a_1, \ldots, a_k\}$.
It is clear that all conditions (1)--(3) of Assumption~\ref{properties of the base structures} hold.
Suppose that $p(x_1, \ldots, x_k)$ is a complete $(\tau, \lambda)$-neighbourhood type for some $\lambda \in \mbbN$.
Then (as one can easily check) $p$ is bounded if and only if, for all $i \in  \{1, \ldots, k\}$ there is $j \in \{1, \ldots, s\}$
such that $p \models P_j(x_i)$.
Similarly, if $p(\bar{x}, y_1, \ldots, y_k)$ is a complete $(\tau,  \lambda)$-neighbourhood type then it is $\bar{y}$-bounded
if and only if, for all $i \in  \{1, \ldots, k\}$ there is $j \in \{1, \ldots, s\}$
such that $p(\bar{x}, y_1, \ldots, y_k) \models P_j(y_i)$.
It follows that an element $a \in B_n$ is rare if and only if $\mcB_n \models P_j(a)$ for some $j \in \{1, \ldots, s\}$.
Hence, for all $\lambda \in \mbbN$ and all $a_1, \ldots, a_k \in (B_n)$, 
\[
C_\lambda^{\mcB_n}(a_1, \ldots, a_k) = \{a_1, \ldots, a_k\} \cup \{b \in B_n : \mcB_n \models P_j(b) \text{ for some $j$}\}.
\]
It is now easy to see that condition~(4) of Assumption~\ref{properties of the base structures} holds if we choose
$m_p = m^s$ and $f_p(n) = n - sm$.
Let $p(\bar{x}, y_1, \ldots, y_k)$ is a complete $(\tau,  \lambda)$-closure type.
Then $p$ is $\bar{y}$-bounded if and only if, for all $i \in  \{1, \ldots, k\}$ there is $j \in \{1, \ldots, s\}$
such that $p \models P_j(y_i)$.
If $p$ is not $\bar{y}$-bounded then, by Lemma~\ref{not bounded implies uniformly unbounded} below,
it is uniformly $\bar{y}$-unbounded.
Moreover, $p$ is strongly $\bar{y}$-unbounded if and only if for all $i = 1, \ldots, k$,
$p \models \bigwedge_{j = 1}^s \neg P_j(x_i)$.
}\end{exam}

\begin{exam}\label{example of paths} {\bf (Lists/paths)} {\rm
Let $\tau = \{E\}$ where $E$ is a binary relation symbol.
For each $n \in \mbbN^+$ let the $\tau$-structure $\mcB_n$ have domain $B_n = \{0, 1,  \ldots, n\}$ and the interpretation
\[
E^{\mcB_n} = \{(i, i+1) : i = 0, 1, \ldots, n-1\}.
\]
So $\mcB_n$ is a ``directed path'' of length $n$.
It is clear that conditions (1)--(3) of
Assumption~\ref{properties of the base structures}
are satisfied.
If $\lambda \in \mbbN$ and $i \in B_n$ then 
\[
N_\lambda^{\mcB_n}(i) = \{j \in B_n : i - \lambda \leq j \leq i + \lambda\}.
\]
If $i, j \in B_n$ then let us call $i$ a {\em predecessor} of $j$, and $j$ a {\em successor} of $i$, if $i < j$.
Let $\lambda \in \mbbN$ and let $p_1(x)$ and $p_2(x_1, x_2)$ be complete $(\tau, \lambda)$-neighbourhood types.
Then $p_1(x)$ is {\em not} bounded if and only if, 
\begin{enumerate}
\item[] $p_1 \models$ ``$x$ has at least $\lambda$ successors and at least $\lambda$ predecessors''.
\end{enumerate}
And $p_2$ is {\em not} bounded if and only if
\begin{enumerate}
\item[] $p_2 \models$ ``for $i=1$ or $i=2$, $x_i$ has at least $\lambda$ successors and at least $\lambda$ predecessors, and 
if $j \neq i$ and $x_j$ has less than $\lambda$ successors or less than $\lambda$ predecessors, 
then $\dist(x_i, x_j) > 2\lambda$''.
\end{enumerate}
An element $i \in B_n$ is $\lambda$-rare if and only if $i$ has fewer than $\lambda$ predecessors or fewer than
$\lambda$ successors.
Suppose that $n > 2\lambda > 0$.
Then the set of $\lambda$-rare elements of $\mcB_n$ is
\[
X = \{0, \ldots, \lambda - 1\} \cup \{n - \lambda + 1, \ldots, n\}
\]
and if $i \in B_n$ then
\begin{align*}
&C_\lambda^{\mcB_n}(i) = N_\lambda^{\mcB_n}(X \cup \{i\}) = \\
&\{0, \ldots, 2\lambda - 1\} \cup \{n - 2\lambda + 1, \ldots, n\} \cup
\{k \in B_n : i - \lambda \leq k \leq i + \lambda\}.
\end{align*}
Let $p(x, y)$ be a complete $(\tau, \lambda)$-closure type. 
Then $p$ is {\em not} $y$-bounded if and only if
\begin{enumerate}
\item[] $p \models$ ``$y$ has at least $\lambda$ successors and at least $\lambda$ successors,
$\dist(x, y) > 2\lambda$,
and for all $z$, if $z$ is $\lambda$-rare then $\dist(y, z) > 2\lambda$''.
\end{enumerate}
Given such observations as above it follows straightforwardly that condition~(4) of Assumption~\ref{properties of the base structures}
is satisfied if $m_p = (2\lambda)^{k+1}$ and $f_p(n) = n - 2\lambda - 2\lambda k$.
}\end{exam}

\begin{exam}\label{example of grids} {\bf (Grids)} {\rm
Let $\tau = \{E\}$ where $E$ is a binary relation symbol.
For every $n \in \mbbN^+$ let the domain of $\mcB_n$ be
\[
B_n = \{(i, j) : i = 0, 1, \ldots, n \text{ and } j = 0, 1, \ldots, n\}
\]
and interpret $E$ in such a way that
$\mcB_n \models E\big((i, j), (k, l)\big)$ if and only if 
$|i - k| + |j - l| = 1$. Thus $\mcB_n$ looks like a two-dimensional grid with side length $n$.
One can verify that all conditions of 
Assumption~\ref{properties of the base structures} hold.
One can also show that for every $\lambda \in \mbbN$ and $b \in B_n$, 
$b$ is $\lambda$-rare if and only if the distance from $b$ to a corner of the grid is less than $\lambda$.

This example is an instance of a class of examples:
For any $d \in \mbbN^+$, we can let $\mcB_n$ be a ``$d$-dimensional grid
with side length $n$''. Then all conditions of Assumption~\ref{properties of the base structures} hold.
}\end{exam}

\begin{exam}\label{example of trees} {\bf (Galton-Watson trees)} {\rm
Let $\tau = \{E, \sqsubset\}$ where $E$ and $\sqsubset$ are  binary relation symbols.
By an {\em ordered rooted tree} we mean a finite $\tau$-structure $\mcT$ such that 
\begin{enumerate}
\item $\mcT \models \forall x, y \big(E(x, y) \rightarrow (x \neq y \wedge \neg E(y, x)) \big)$,

\item there is a unique element $a \in T$, called the {\em root}, such that $\mcT \models \forall x \neg E(x, a)$, 

\item for all $a \in T$, if $a$ is not the root then there is a unique $b \in T$ 
(called the {\em parent} of $a$) such that $\mcT \models E(b, a)$
(and in this case $a$ is called a {\em child} of $b$), 

\item there do {\em not} exist $k \in \mbbN^+$ and $a_1, \ldots, a_k \in T$ such that $\mcT \models E(a_k, a_1)$
and $\mcT \models E(a_i, a_{k+1})$ for all $i = 1, \ldots, k-1$, and

\item for all $a \in T$, $\sqsubset$ is a strict linear order on the children of $a$ (if $a$ has any child), and
if $a \in T$ and $b \in T$ do not have a common parent then $a$ and $b$ are incomparable with respect to $\sqsubset$.
\end{enumerate}
For simplicity we will often just say {\em tree} when actually meaning {\em ordered rooted tree}.
If $\delta$ is a natural number then we say that a tree is {\em $\delta$-bounded } if every vertex of it has at most $\delta$ children.
Fix some $\delta \in \mbbN^+$ and for each $n \in \mbbN^+$ let $\mcB_n$ be a 
$\delta$-bounded tree with exactly $n$ vertices. 
Let $\Delta = \delta + 1$. Then conditions~(1) --~(3) of Assumption~\ref{properties of the base structures}
are clearly satisfied. 
It is possible to choose the sequence $\mbB = (\mcB_n : n \in \mbbN^+)$ so that also condition~(4)
of  Assumption~\ref{properties of the base structures} is satisfied.
In fact, under some additional conditions on the trees, most choices 
(in a probabilistic sense) of a $\delta$-bounded ordered rooted tree with $n$ vertices $\mcB_n$, for $n \in \mbbN^+$,
will result in a sequence of $\tau$-structures that satisfies condition~(4)
of  Assumption~\ref{properties of the base structures}.
We do not show this here because it will follow when we later,
in Example~\ref{Galton-Watson trees satisfy a stronger property},
show that a ``random sequence'' $\mbB = (\mcB_n : n \in \mbbN^+)$ of Galton-Watson trees
will, with high probability, have a property that implies both 
condition~(4) of Assumption~\ref{properties of the base structures} 
and another condition,
Assumption~\ref{relative frequency of tau-closure types}.
For this we will use a result of Janson \cite{Jan21} about the number of embeddings of a fixed tree into a random
Galton-Watson tree.
}\end{exam}

\section{Properties of bounded and unbounded closure types}\label{properties of closure types}

\noindent
In this section we make a more detailed study of bounded and unbounded closure types
(Definition~\ref{definition of bounded}).
The notions and (or) results from this section will be used from Section~\ref{convergence and balance}
and in particular in sections~\ref{Finding the balance in the inductive step}
and~\ref{Asymptotic elimination of aggregation functions}.
We pay such attention to closure types because, in the main results
(in Section~\ref{Asymptotic elimination of aggregation functions}) we prove results
concerning asymptotic elimination of aggregation functions for $PLA^*$-formulas 
that use closure types as conditioning subformulas 
(see Definition~\ref{definition of subformula}).
In particular we will see that if $p(\bar{x}, \bar{y})$ is a complete $(\tau, \lambda)$-closure type that is 
{\em not} $\bar{y}$-bounded, then $p$ is uniformly $\bar{y}$-unbounded
(Definition~\ref{definition of unbounded formula}).
For $p(\bar{x}, \bar{y})$ as above, we define a notion of {\em $\bar{y}$-dimension}, $\dim_{\bar{y}}(p)$,
and prove (among other things) that $p$ is {\em not} $\bar{y}$-bounded if and only if $\dim_{\bar{y}}(p) \geq 1$.
Moreover, we will see that for every complete $(\tau, \lambda)$-closure type $p(\bar{x}, \bar{y})$
the variables $\bar{y}$ can be partititioned into two parts $\bar{u}$ and $\bar{v}$ such that 
the restriction of $p$ to $\bar{x}\bar{u}$ is $\bar{u}$-bounded and $p$ is {\em strongly} $\bar{v}$-unbounded
(Definition~\ref{definition of unbounded formula}).
The notion of strongly $\bar{y}$-unbounded (for some $\bar{y}$) closure types will be essential in
the proofs in Section~\ref{Finding the balance in the inductive step}.

{\bf \em We assume that $\tau \subseteq \sigma$ are finite relational signatures and that
$\mbB = (\mcB_n : n \in \mbbN^+)$ is a sequence of $\tau$-structures that satisfies
Assumption~\ref{properties of the base structures}.}
In this section we work almost exclusively with the $\tau$-structures in the sequence $\mbB$
and with $(\tau, \lambda)$-closure types,
and use definitions and results from Sections~\ref{the base sequence}, 
but the notion of dimension is also defined for $(\sigma, \lambda)$-closure types (for later use).

\begin{lem}\label{transitivity of boundedness}
Let $\lambda \in \mbbN$, $p(\bar{x}, \bar{y}, \bar{z})$ be a complete
$(\tau, \lambda)$-neighbourhood type, and $q(\bar{x}, \bar{y}) = p \uhrc \bar{x}\bar{y}$.
If $q$ is $\bar{y}$-bounded and $p$ is $\bar{z}$-bounded, then $p$ is $\bar{y}\bar{z}$-bounded.
\end{lem}

\noindent
{\bf Proof.}
The assumption that $q$ is $\bar{y}$-bounded and $p$ is $\bar{z}$-bounded means that there are
numbers $m_p, m_q \in \mbbN$ such that for all $n$, all $\bar{a} \in (B_n)^{|\bar{x}|}$ 
and all $\bar{b} \in (B_n)^{|\bar{y}|}$, we have
$|q(\bar{a}, \mcB_n)| \leq m_q$ and $|p(\bar{a}, \bar{b}, \mcB_n)| \leq m_p$.
It follows that for all $n$ and $\bar{a} \in (B_n)^{|\bar{x}|}$, 
$|p(\bar{a}, \mcB_n)| \leq m_q \cdot m_p$ so $p$ is $\bar{y}\bar{z}$-bounded.
\hfill $\square$

\begin{lem}\label{not bounded implies uniformly unbounded}
Suppose that $\lambda \in \mbbN$ and that $p(\bar{x}, \bar{y})$ is a complete $(\tau, \lambda)$-closure type.
If $p$ is not $\bar{y}$-bounded then $p$
is uniformly $\bar{y}$-unbounded, that is, 
$p$ is cofinally satisfiable and there is $f_p : \mbbN \to \mbbR$ such that, for all $\alpha > 0$,
$\lim_{n\to\infty}(f_p(n) - \alpha\ln(n)) = \infty$ and, for all sufficiently large $n \in \mbbN$,
if $\bar{a} \in (B_n)^{|\bar{x}|}$ and $p(\bar{a}, \mcB_n) \neq \es$, 
then $|p(\bar{a}, \mcB_n)| \geq f_p(n)$.
In the particular case when $\bar{x}$ is empty we have that $|p(\mcB_n)| \geq f_p(n)$ for all $n$ such that 
$p(\mcB_n) \neq \es$.
\end{lem}

\noindent
{\bf Proof.}
Let $p(\bar{x}, \bar{y})$ be a complete $(\tau, \lambda)$-closure type that is not $\bar{y}$-bounded.
Let $\bar{u}$ be a maximal subsequence of $\bar{y}$ such that $p \uhrc \bar{x}\bar{u}$ is $\bar{u}$-bounded.
(When saying that $\bar{u}$ is maximal we mean that if $\bar{u}' $ is a subsequence of $\bar{y}$ such that
$\rng(\bar{u}) \subset \rng(\bar{u}')$, then $p \uhrc \bar{x}\bar{u}'$ is not $\bar{u}'$-bounded.)
Since $p(\bar{x}, \bar{y})$ is not $\bar{y}$-bounded it follows that $\bar{u}$ must be a proper subsequence of $\bar{y}$.

Let $\bar{v}$ contain all variables in $\bar{y}$ that do not occur in $\bar{u}$.
Without loss of generality we can assume that $\bar{y} = \bar{u}\bar{v}$.
Suppose towards a contradiction that there is $v_i \in \rng(\bar{v})$ such that
$p(\bar{x}, \bar{u}, \bar{v}) \models \dist(\bar{x}\bar{u}, v_i) \leq 2\lambda$.
Then (as every $\mcB_n$ has degree at most $\Delta$ by
Assumption~\ref{properties of the base structures})
$p \uhrc \bar{x}\bar{u}v_i$ is $v_i$-bounded, and as 
$p \uhrc \bar{x}\bar{u}$ is $\bar{u}$-bounded it follows from
Lemma~\ref{transitivity of boundedness}
that $p \uhrc \bar{x}\bar{u}v_i$ is $\bar{u}v_i$-bounded which contradicts the maximality of $\bar{u}$.
Hence, by Lemma~\ref{closure types and cases}, we have
\begin{equation}\label{p1 implies large distance}
p(\bar{x}, \bar{u}, \bar{v}) \models \dist(\bar{x}\bar{u}, \bar{v}) > 2\lambda.
\end{equation}
Recall Definition~\ref{definition of dimension for closure types} of $\sim_p$.
Choose any $v_i \in \rng(\bar{v})$ and let $\bar{v}''$ enumerate the $\sim_p$-class of $v_i$.
Let $\bar{v}'$ enumerate the rest of the variables in $\bar{v}$. 
Without loss of generality we may assume that $\bar{v} = \bar{v}'\bar{v}''$.
It follows from~(\ref{p1 implies large distance}) and 
Definition~\ref{definition of dimension} of $\sim_p$
that 
\begin{equation}\label{p1 implies large distance, second time}
p(\bar{x}, \bar{u}, \bar{v}', \bar{v}'') \models \dist(\bar{x}\bar{u}\bar{v}', \bar{v}'') > 2\lambda.
\end{equation}

Suppose for a contradiction that, for some $v_i \in \rng(\bar{v}'')$, 
$p(\bar{x}, \bar{y}) \models \exists z \big( \varphi_\lambda(z) \wedge \dist(z, v_i) \leq 2\lambda\big)$
where $\varphi_\lambda(z)$ expresses (in all $\mcB_n$) that $z$ is $\lambda$-rare.
Since, for some $m$, there are at most $\lambda$-rare elements in each $\mcB_n$ and the $2\lambda$-neighbourhood of
an element has cardinality at most $\Delta^{2\lambda} + 1$ it follows that $p\uhrc \bar{x}\bar{u}v_i$ is $\bar{u}v_i$-bounded
which contradicts the maximality of $\bar{u}$.
By Lemma~\ref{closure types and cases} we get
\begin{equation}\label{all y are far from rare elements}
p(\bar{x}, \bar{u}, \bar{v}', \bar{v}'') \models \forall z \big( \varphi_\lambda(z) \rightarrow \dist(z, \bar{v}'') > 2\lambda\big)
\end{equation}

Let $p_1(\bar{x}, \bar{u}, \bar{v}') = p \uhrc \bar{x}\bar{u}\bar{v}'$ and
$p_2(\bar{v}'') = p \uhrc \bar{v}''$.
According to~(\ref{p1 implies large distance, second time}), (\ref{all y are far from rare elements}) 
and Lemma~\ref{splitting a closure type},
there is a complete $(\tau, \lambda)$-{\em neighbourhood type} $p_2^+(\bar{v}'')$ such that
$p_2(\bar{v}'') \models p_2^+(\bar{v}'')$ and for all $n$, $\bar{a} \in (B_n)^{|\bar{x}|}$,
$\bar{b} \in (B_n)^{|\bar{u}|}$, $\bar{c}' \in (B_n)^{|\bar{v}'|}$ and $\bar{c}'' \in (B_n)^{|\bar{v}''|}$,
\begin{align}\label{splitting p}
&\mcB_n \models p(\bar{a}, \bar{b}, \bar{c}', \bar{c}'') \text{ if and only if } \\
&\mcB_n \models p_1(\bar{a}, \bar{b}, \bar{c}') \wedge 
p_2^+(\bar{c}'') \text{ and } \dist(\bar{a}\bar{b}\bar{c}'\bar{d}, \ \bar{c}'') > 2\lambda \nonumber
\end{align}
where $\bar{d}$ enumerates all $\lambda$-rare elements in $B_n$.

If $p_2(\bar{v}'')$ would be bounded then $p \uhrc \bar{x}\bar{u}\bar{v}''$ would be $\bar{u}\bar{v}''$-bounded,
contradicting the maximality of $\bar{u}$.
Hence $p_2(\bar{v}'')$ is not bounded and,
as $p_2(\bar{v}'') \models p_2^+(\bar{v}'')$, it follows that $p_2^+(\bar{v}'')$ is not bounded.
Since $p_2^+$ is a complete $(\tau, \lambda)$-neighbourhood type 
and all of its variables are in the same $\sim_p$-class, hence in the same $\sim_{p_2^+}$-class, it follows from
part~(4) of 
Assumption~\ref{properties of the base structures}
that there is $f_{p_2^+} : \mbbN \to \mbbR$ such that, for all $\alpha > 0$, 
$\lim_{n\to\infty} (f_{p_2^+}(n) - \ln(n)) = \infty$ and
\begin{equation}\label{p2 tends to infinity}
\text{for all sufficiently large $n$}, \  |p_2^+(\mcB_n)| \geq f_{p_2^+}(n).
\end{equation}

From the assumption that $p(\bar{x}, \bar{y})$ not $\bar{y}$-bounded it follows that
$p(\bar{x}, \bar{y})$ is cofinally satisfiable. 
Suppose that for some $n$ and $\bar{a} \in (B_n)^{|\bar{x}|}$ we have  $p(\bar{a}, \mcB_n) \neq \es$.
Then there are $\bar{b} \in (B_n)^{|\bar{u}|}$ and $\bar{c}' \in (B_n)^{|\bar{v}'|}$ 
such that $p(\bar{a}, \bar{b}, \bar{c}', \mcB_n) \neq \es$ and hence
$\mcB_n \models p_1(\bar{a}, \bar{b}, \bar{c}')$.
Recall that by 
Lemma~\ref{definability of m-rareness}
there is $m \in \mbbN$ such that for all $n$, $\mcB_n$ has at most $m$ $\lambda$-rare elements.
By~(\ref{splitting p}), if $\bar{d}$ enumerates all $\lambda$-rare elements in $\mcB_n$,
$\dist(\bar{a}\bar{b}\bar{c}'\bar{d}, \ \bar{c}'') > 2\lambda$ and
$\mcB_n \models p_2^+(\bar{c}'')$, then $\mcB_n \models p(\bar{a}, \bar{b}, \bar{c}', \bar{c}'')$.
Since the degree of $\mcB_n$ is bounded by the fixed number $\Delta$
it follows that there is some $C \in \mbbN$ (depending only on $\Delta$ and $\lambda$) 
such that there are at most $C$ tuples $\bar{c}'' \in p_2^+(\mcB_n)$
such that $\dist(\bar{a}\bar{b}\bar{c}'\bar{d}, \ \bar{c}'') \leq 2\lambda$.
Therefore at most $C$ tuples in $p_2^+(\mcB_n)$ do {\em not} satisfy $p(\bar{a}, \bar{b}, \bar{c}, \bar{d}', \bar{v}'')$.
So if we let $f_p(n) = f_{p_2^+}(n)  - C$ then
it follows from~(\ref{splitting p}) and~(\ref{p2 tends to infinity}) that
$|p(\bar{a}, \bar{b}, \bar{c}', \mcB_n)| \geq f_p(n)$, 
hence $|p(\bar{a}, \mcB_n)| \geq f_p(n)$ for all sufficiently large $n$, and we also have
$\lim_{n\to\infty}(f_p(n) - \alpha\ln(n)) = \infty$ for all $\alpha > 0$.
\hfill $\square$

\medskip

\noindent
The proof of Lemma~\ref{not bounded implies uniformly unbounded} 
shows the following (as~(\ref{p1 implies large distance, second time}) and (\ref{all y are far from rare elements}) implies that 
$p$ in that proof is uniformly $\bar{v}$-unbounded):

\begin{cor}\label{corollary to not bounded implies uniformly unbounded}
Suppose that $\lambda \in \mbbN$ and that $p(\bar{x}, \bar{y})$ is a complete $(\tau, \lambda)$-closure type
such that 
\[
p(\bar{x}, \bar{y}) \models \dist(\bar{x}, \bar{y}) > 2\lambda \ \wedge \ 
\forall z \big(\varphi_\lambda(z) \rightarrow \dist(z, \bar{y}) > 2\lambda\big)
\]
where $\varphi_\lambda(z)$ expresses (in all $\mcB_n$) that $z$ is $\lambda$-rare.
Then $p$ is uniformly $\bar{y}$-unbounded.
\end{cor}

\noindent
The next result tells that a neighbourhood type, or a closure type, can be ``decomposed'' 
into a bounded part and a strongly unbounded part.

\begin{lem}\label{a type can be divided in a bounded and strongly unbounded type}
Let $\lambda \in \mbbN$ and let $p(\bar{x}, \bar{y})$ be a complete $(\tau, \lambda)$-neighbourhood type
or a complete $(\tau, \lambda)$-closure type.
Then there is a subsequence  $\bar{u}$ of $\bar{y}$ such that $p \uhrc \bar{x}\bar{u}$ is $\bar{u}$-bounded
and if $\bar{v}$ is the subsequence of all variables in $\bar{y}$ which do not occur in $\bar{u}$
and $\bar{v}$ is nonempty, 
then $p$ is strongly $\bar{v}$-unbounded.
\end{lem}

\noindent
{\bf Proof.}
Let $\bar{u}$ be a maximal subsequence of $\bar{y}$ such that $p \uhrc \bar{x}\bar{u}$ is $\bar{u}$-bounded.
Let $\bar{v}$ be the sequence of all variables in $\bar{y}$ that do not occur in $\bar{u}$
and suppose that $\bar{v}$ is nonempty (otherwise there is nothing to prove).
Without loss of generality we may assume that $\bar{y} = \bar{u}\bar{v}$.

It now suffices to show that $p$ is strongly $\bar{v}$-unbounded.
Towards a contradiction suppose that $p$ is not strongly $\bar{v}$-unbounded.
Then there is a nonempty subsequence $\bar{v}'$ of $\bar{v}$ 
such that $p \uhrc \bar{x}\bar{u}\bar{v}'$ is not uniformly $\bar{v}'$-bounded.
Lemma~\ref{not bounded implies uniformly unbounded} implies that
$p \uhrc \bar{x}\bar{u}\bar{v}'$ is $\bar{v}'$-bounded.
It now follows from
Lemma~\ref{transitivity of boundedness}
that $p \uhrc \bar{x}\bar{u}\bar{v}'$ is $\bar{u}\bar{v}'$-bounded
which contradicts that $\bar{u}$ is a maximal subsequence of $\bar{y}$ 
such that $p \uhrc \bar{x}\bar{u}$ is $\bar{u}$-bounded.
Hence $p$ is strongly $\bar{v}$-unbounded.
\hfill $\square$
\\

\noindent
Recall Definitions~\ref{definition of dimension} and~\ref{definition of dimension for closure types}
of `$\sim_p$' for neighbourhood types and closure types, respectively.
The next lemma gives an alternative characterization of boundedness (of closure types) which will be used later.

\begin{lem}\label{a property of bounded closure types}
Let $\lambda \in \mbbN$ and let $p(\bar{x}, \bar{y})$ be a complete $(\tau, \lambda)$-closure type
which is cofinally satisfiable.
Then $p(\bar{x}, \bar{y})$ is $\bar{y}$-bounded if and only if the following holds:
\begin{enumerate}
\item[($*$)] If $\bar{v}$ is a subsequence of $\bar{x}\bar{y}$ that contains all variables from one (and only one)
$\sim_p$-class and  $\rng(\bar{v}) \cap \rng(\bar{y}) \neq  \es$,
then $p \uhrc \bar{v}$ is bounded or $\rng(\bar{v}) \cap \rng(\bar{x}) \neq \es$.
\end{enumerate}
\end{lem}

\noindent
{\bf Proof.}
Suppose that ($*$) holds.
Let $\bar{v}_1, \ldots, \bar{v}_k$ be a list of tuples where each $\bar{v}_i$ enumerates all 
variables in one (and only one) $\sim_p$-class such that $\rng(\bar{v}_i) \cap \rng(\bar{y}) \neq \es$.
From part~(2) of Assumption~\ref{properties of the base structures} (the bound $\Delta$ of the degree),
the assumption~($*$), the definitions of $\bar{z}$-bounded formula (for some $\bar{z}$) and $\sim_p$ it 
straightforwardly follows that, 
for all $i = 1, \ldots, k$, $p \uhrc \bar{x}\bar{v}_i$ is $\bar{v}_i$-bounded.
As every variable in $\bar{y}$ belongs to some $\bar{v}_i$ it follows that $p$ is $\bar{y}$-bounded.

Now suppose that $p$ is cofinally satisfiable and that ($*$) does not hold.
Then there is a subsequence $\bar{v}$ of $\bar{x}\bar{y}$ that contains all variables from one (and only one)
$\sim_p$-class, $\rng(\bar{v}) \cap \rng(\bar{y}) \neq  \es$,
$p \uhrc \bar{v}$ is not bounded, and $\rng(\bar{v}) \cap \rng(\bar{x}) = \es$.
Without loss of generality we can assume that $\bar{y} = \bar{u}\bar{v}$.
Since $\rng(\bar{v}) \cap \rng(\bar{x}) = \es$ and $p \uhrc \bar{v}$ is bounded it follows that 
\[
p(\bar{x}, \bar{u}, \bar{v}) \models \dist(\bar{x}\bar{u}, \bar{v}) > 2\lambda \ \wedge \
\forall z \big( \varphi_\lambda(z) \rightarrow \dist(z, \bar{v}) > 2\lambda \big),
\]
where $\varphi_\lambda(z)$ expresses (in all $\mcB_n$) that $z$ is $\lambda$-rare.
By Corollary~\ref{corollary to not bounded implies uniformly unbounded},
$p(\bar{x}, \bar{u}, \bar{v})$ is uniformly $\bar{v}$-unbounded.
So $p$ is cofinally satisfiable and there is $f : \mbbN \to \mbbR$ such that $\lim_{n\to\infty}f(n) = \infty$ and for all large enough $n$,
if $\bar{a} \in (B_n)^{|\bar{x}|}$, $\bar{b} \in (B_n)^{|\bar{u}|}$ 
and $p(\bar{a}, \mcB_n) \neq \es$ then $|p(\bar{a}, \bar{b}, \mcB_n)| \geq f(n)$, hence also
$p(\bar{a}, \mcB_n) \geq f(n)$. Thus $p$ is not $\bar{y}$-bounded.
\hfill $\square$

\medskip

\noindent
The next lemma relates the notion of boundedness of a closure type with 
the notion of closure of a sequence of elements (Definition~\ref{definition of closure}).

\begin{lem}\label{boundedness and closure}
Let $\lambda \in \mbbN$ and suppose that $p(\bar{x}, \bar{y})$ is a $\bar{y}$-bounded complete $(\tau, \lambda)$-closure type
that is cofinally satisfiable.
Then there is $\gamma \in \mbbN$ such that for all $n$ and all $\bar{a} \in (B_n)^{|\bar{x}|}$,
if $\bar{b} \in p(\bar{a}, \mcB_n)$ then $\rng(\bar{b}) \subseteq C_{\gamma}^{\mcB_n}(\bar{a})$.
\end{lem}

\noindent
{\bf Proof.}
Let $p(\bar{x}, \bar{y})$ be as assumed in the lemma and let
$\bar{x} = (x_1, \ldots, x_k)$ and $\bar{y} = (y_1, \ldots, y_l)$.
By Lemma~\ref{a property of bounded closure types}
the following holds:
\begin{enumerate}
\item[($*$)] If $\bar{v}$ is a subsequence of $\bar{x}\bar{y}$ that contains all variables from one (and only one)
$\sim_p$-class and  $\rng(\bar{v}) \cap \rng(\bar{y}) \neq  \es$,
then $p \uhrc \bar{v}$ is bounded or $\rng(\bar{v}) \cap \rng(\bar{x}) \neq \es$.
\end{enumerate}
Suppose that $\bar{a} = (a_1, \ldots, a_k) \in (B_n)^k$, $\bar{b} = (b_1, \ldots, b_l) \in (B_n)^l$,
and $\mcB_n \models p(\bar{a}, \bar{b})$.
Fix any $j \in \{1, \ldots, l\}$.
Let $\bar{v}$ be the subsequence of variables $v$ in $\bar{x}\bar{y}$ such that $v \sim_p y_j$.
In particular, $y_j$ occurs in $\bar{v}$.
Suppose that $\bar{v} = (x_{i_1}, \ldots, x_{i_s}, y_{j_1}, \ldots, y_{j_t})$ and let
$\bar{c} = (a_{i_1}, \ldots, a_{i_s}, b_{j_1}, \ldots, b_{j_t})$

According to~($*$), $p \uhrc \bar{v}$ is bounded or $\rng(\bar{v}) \cap \rng(\bar{x}) \neq \es$.
If $\rng(\bar{v}) \cap \rng(\bar{x}) \neq \es$ then there is $m_j \in \mbbN$ that depends only on $p(\bar{x}, \bar{y})$
and $j$
such that for some $i \in \{1, \ldots, k\}$, $\dist(a_i, b_j) \leq m_j$ and hence
$b_j \in C_{m_j}^{\mcB_n}(\bar{a})$. (One can take $m_j = 2\lambda|\bar{x}\bar{y}|$.)
Otherwise $p \uhrc \bar{v}$ is bounded and it follows that there are
$m_j \in \mbbN$ and a bounded complete $(\tau, m_j)$-neighbourhood type $q(x)$ such that $\mcB_n \models q(b_j)$,
so $b_j$ is $m_j$-rare, and hence $b_j \in  C_{m_j}^{\mcB_n}(\bar{a})$.
It follows that if  $\gamma = \max\{m_j : j = 1, \ldots, l\}$, then
$\rng(\bar{b}) \subseteq C_\gamma^{\mcB_n}(\bar{a})$.
\hfill $\square$

\medskip

\noindent
We now define a notion of dimension (of neighbourhood and closure types) which will play an important role
in Section~\ref{Finding the balance in the inductive step} where certain things are proved by induction on the dimension.

\begin{defin}\label{definition of dimension, part 2} {\rm
(i) Let $p(\bar{x})$ be a complete $(\tau, \lambda)$-neighbourhood type or a 
complete $(\tau, \lambda)$-closure type (for some $\lambda \in \mbbN$).
\begin{enumerate}
\item Suppose that $\bar{y}$ is a nonempty subsequence of $\bar{x}$ (and we allow the special case when $\bar{y} = \bar{x}$).
The {\bf \em $\bar{y}$-dimension of $p$}, denoted $\dim_{\bar{y}}(p)$, 
is the number of equivalence classes $\mb{e}$ of  $\sim_p$ 
(see Definitions~\ref{definition of dimension} and~\ref{definition of dimension for closure types}) 
such that
\begin{enumerate}
\item $\mb{e}$ contains only variables from $\bar{y}$, and 
\item $p$ restricted to the variables in $\mb{e}$ is not bounded.
\end{enumerate}

\item The {\bf \em dimension of $p$}, denoted $\dim(p)$, is defined to be $\dim_{\bar{x}}(p)$.
\end{enumerate}
(ii) Let $p(\bar{x})$ be a $(\sigma, \lambda)$-neighbourhood type or a $(\sigma, \lambda)$-closure type
(for some $\lambda \in \mbbN$). If $\bar{y}$ is a nonempty subsequence of $\bar{x}$, then
$\dim_{\bar{y}}(p) = \dim_{\bar{y}}(p \uhrc \tau)$ and $\dim(p) = \dim(p \uhrc \tau)$, where $p \uhrc \tau$ denotes
the restriction of $p$ to $\tau$ (Definition~\ref{restriction of closure type}).
}\end{defin}

\begin{lem}\label{dimension and unboundedness}
Let $p(\bar{x}, \bar{y})$ be a complete $(\tau, \lambda)$-closure type for some $\lambda \in \mbbN$
and suppose that $p$ is cofinally satisfiable.
Then $\dim_{\bar{y}}(p) \geq 1$ if and only if $p$ is not $\bar{y}$-bounded.
\end{lem}

\noindent
{\bf Proof.}
Suppose that $p(\bar{x}, \bar{y})$ is cofinally satisfiable and not $\bar{y}$-bounded.
By Lemma~\ref{a property of bounded closure types},
there is a subsequence $\bar{v}$ of $\bar{y}$ which enumerates  all variables of one $\sim_p$-class
and $p \uhrc \bar{v}$ is unbounded. 
It now follows from the definition of $\dim_{\bar{y}}(p)$ that $\dim_{\bar{y}}(p) \geq 1$.

Now suppose that $p$ is cofinally satisfiable and $\dim_{\bar{y}}(p) \geq 1$.
This means that there is (at least) one $\sim_p$-class such that if $\bar{v}$ enumerates all variables in it,
then $\rng(\bar{v}) \subseteq \rng(\bar{y})$ and
$p \uhrc \bar{v}$ is unbounded.
Lemma~\ref{a property of bounded closure types} implies that $p$ is not $\bar{y}$-bounded.
\hfill $\square$

\medskip

\noindent
We now give a characterization of strongly unbounded closure types 
(Definition~\ref{definition of unbounded formula})
that will be used later.

\begin{lem}\label{another characterization of strong unboundedness}
Let $\lambda \in \mbbN$, $\bar{x} = (x_1, \ldots, x_k)$, $\bar{y} = (y_1, \ldots, y_l)$ (where $l \geq 1$), 
and let $p(\bar{x}, \bar{y})$ be a complete $(\tau, \lambda)$-closure type.
Then $p(\bar{x}, \bar{y})$ is strongly $\bar{y}$-unbounded
if and only if the following conditions hold:
\begin{enumerate}
\item For all $i = 1, \ldots, k$ and all $j = 1, \ldots, l$, $x_i \not\sim_p y_j$
(which is equivalent to say that $p(\bar{x}, \bar{y}) \models \dist(\bar{x}, \bar{y}) > 2\lambda$).

\item For every subsequence $\bar{u}$ of $\bar{y}$, $p \uhrc \bar{u}$ is not bounded.
\end{enumerate}
\end{lem}

\noindent
{\bf Proof.}
Suppose first that $p(\bar{x}, \bar{y})$ is strongly $\bar{y}$-unbounded.
If $x_i \sim_p y_j$ for some $i$ and $j$
then, letting $m = 2\lambda|\bar{x}\bar{y}|$, we get
$p(\bar{x}, \bar{y}) \models \dist(x_i, y_j) \leq m$ and
therefore $p \uhrc \bar{x}y_j$ is $y_j$-bounded which contradicts that $p$ is strongly $\bar{y}$-unbounded.
Hence condition~(1) holds.
If, for some subsequence $\bar{u}$ of $\bar{y}$, $p \uhrc \bar{u}$ is bounded then
$p \uhrc \bar{x}\bar{u}$ is $\bar{u}$-bounded
which contradicts that $p$ is strongly $\bar{y}$-unbounded.
Hence condition~(2) holds.

Now suppose that conditions~(1) and~(2) hold.
Let $\bar{u}$ be a subsequence of $\bar{y}$.
We need to show that $p \uhrc \bar{x}\bar{u}$ is uniformly $\bar{u}$-unbounded.
From~(1) it follows that $p(\bar{x}, \bar{y}) \models \dist(\bar{x}, \bar{y}) > 2\lambda$
and hence $p  \uhrc \bar{x}\bar{u} \models \dist(\bar{x}, \bar{u}) > 2\lambda$.
From~(2) it follows that $p \uhrc \bar{u}$ is not bounded.
This implies that
$p \uhrc \bar{x}\bar{u} \models \forall z \big( \text{``$z$ is $\lambda$-rare'' } \rightarrow \dist(z, \bar{u}) > 2\lambda \big)$.
Now 
Corollary~\ref{corollary to not bounded implies uniformly unbounded}
implies that $p \uhrc \bar{x}\bar{u}$ is uniformly $\bar{u}$-unbounded.
\hfill $\square$

\medskip

\noindent
The next lemma will be used in Section~\ref{Finding the balance in the inductive step}
in to get a base case (for unbounded closure types) in arguments by induction on the dimension.

\begin{lem}\label{getting a strongly unbounded type of dimension 1}
Suppose that, for some $\lambda \in \mbbN$, $p(\bar{x}, \bar{y})$ is a 
uniformly $\bar{y}$-unbounded complete $(\tau, \lambda)$-closure type.
Then there is a subsequence $\bar{v}$ of $\bar{y}$ such that $p$ is strongly $\bar{v}$-unbounded
and $\dim_{\bar{v}}(p) = 1$.
\end{lem}

\noindent
{\bf Proof.}
By Lemma~\ref{a type can be divided in a bounded and strongly unbounded type},
we may assume that $\bar{y} = \bar{u}\bar{w}$, $p \uhrc \bar{x}\bar{u}$ is $\bar{u}$-bounded
and $p$ is strongly $\bar{w}$-unbounded.
From Lemma~\ref{another characterization of strong unboundedness}
it follows that for all $x_i \in \rng(\bar{x})$, 
$y_j \in \rng(\bar{u})$ and $y_k \in \rng(\bar{w})$,
$x_i \not\sim_p y_k$ and $y_j \not\sim_p y_k$.

Let $y_k \in \rng(\bar{w})$ and let $\bar{v}$ enumerate the $\sim_p$ class of $y_k$, so $\rng(\bar{v}) \subseteq \rng(\bar{w})$.
Then, by Lemma~\ref{elementary properties of closure types},
$p$ is strongly $\bar{v}$-unbounded which implies that $p \uhrc \bar{v}$ is not bounded.
It now follows from the definition of dimension that $\dim_{\bar{v}}(p) = 1$.
\hfill $\square$

\medskip

\noindent
The next lemma generalizes part~(4) of Assumption~\ref{properties of the base structures}
to higher dimensions.

\begin{lem}\label{strongly unbounded neighbourhood types tend to infinity}
Let $p(\bar{y})$ be a strongly unbounded complete $(\tau, \lambda)$-{\rm neighbourhood} type for some $\lambda \in \mbbN$.
Then there is $f_p : \mbbN \to \mbbR$ such that, for all $\alpha > 0$,  $\lim_{n\to\infty} (f_p(n) - \alpha\ln(n)) = \infty$
and for all sufficiently large $n$ we have $|p(\mcB_n)| \geq f_p(n)$.
\end{lem}

\noindent
{\bf Proof.}
Let $p(\bar{y})$ be as assumed in the lemma.
We use induction on $\dim(p)$ which is at least 1, by Lemma~\ref{dimension and unboundedness}.
If $\dim(p) = 1$ then the conclusion follows directly from the definition of $\dim(p)$ and part~(4) of 
Assumption~\ref{properties of the base structures}.

Now suppose that $\dim(p) = k > 1$.
Let $\bar{v}$ enumerate a $\sim_p$-class and let $\bar{u}$ enumerate the rest of the variables in $\bar{y}$.
Without loss of generality, suppose that $\bar{y} = \bar{u}\bar{v}$.
By Lemma~\ref{elementary properties of closure types}
both $p_1(\bar{u}) = p \uhrc \bar{u}$ and $p_2 = p \uhrc \bar{v}$ are strongly unbounded.
Since $\bar{v}$ enumerates a $\sim_p$-class it follows 
(since $\sim_p$ and $\sim_{p_2}$ coincide on $\rng(\bar{v})$)
that it also enumerates a $\sim_{p_2}$-class, so $\dim(p_2) = 1$.
Also $\dim(p_1) = \dim_{\bar{u}}(p) < k$.
By the induction hypothesis, there are $f_1, f_2 : \mbbN \to \mbbN$ such that for $i = 1, 2$
and all $\alpha >  0$, $\lim_{n\to\infty}(f_i(n) - \alpha\ln(n)) = \infty$ and 
\begin{equation}\label{two limits to infinity}
\text{for all large enough $n$, } \ |p_1(\mcB_n)| \geq f_1(n) \ \text{ and } \ |p_2(\mcB_n)| \geq f_2(n).
\end{equation}
By the choice of $\bar{v}$ and definition of $\sim_p$ we have 
$p(\bar{u}, \bar{v}) \models \dist(\bar{u}, \bar{v}) > 2\lambda$.
Let $n$ be large enough so that $|p_1(\mcB_n)| \geq f_1(n) > 0$ and $|p_2(\mcB_n)| \geq f_2(n) > 0$.
Take any $\bar{a} \in p_1(\mcB_n)$. 
By Assumption~\ref{properties of the base structures}
there is $m \in \mbbN$ which is independent of $n$ and $\bar{a}$ and such that 
$|p_2(\mcB_n) \setminus N_{2\lambda}^{\mcB_n}(\bar{a})| \geq f_2(n) - m$.
By Lemma~\ref{splitting a neighbourhood type}, for every $\bar{b} \in p_2(\mcB_n) \setminus N_{2\lambda}^{\mcB_n}(\bar{a})$
we have $\mcB_n \models p(\bar{a}, \bar{b})$, so $|p(\mcB_n)| \geq f_2(n) - m$.
It follows that if $f_p(n) = f_2(n) - m$ and $n$ is sufficiently large then $|p(\mcB_n)| \geq f_p(n)$ and,
for all $\alpha > 0$,
$\lim_{n\to\infty} (f_p(n) - \alpha\ln(n)) = \infty$.
\hfill $\square$

\medskip

\noindent
We conclude this section with three lemmas which give more information about closure types and
will be used in Section~\ref{Finding the balance in the inductive step} or
Section~\ref{Asymptotic elimination of aggregation functions}.

\begin{lem}\label{for strongly y-unbounded the x-part determines if y exists}
Suppose that, for some $\lambda \in \mbbN$, $p(\bar{x}, \bar{y})$ is a 
strongly $\bar{y}$-unbounded complete $(\tau, \lambda)$-closure type.
Let $p_1(\bar{x}) = p \uhrc \bar{x}$. For all sufficiently large $n$ and all $\bar{a} \in (B_n)^{|\bar{x}|}$,
$p(\bar{a}, \mcB_n) \neq \es$ if and only if $\mcB_n \models p_1(\bar{a})$.
\end{lem}

\noindent
{\bf Proof.}
Suppose that, for some $\lambda \in \mbbN$, $p(\bar{x}, \bar{y})$ is a 
strongly $\bar{y}$-unbounded complete $(\tau, \lambda)$-closure type.
Let $p_1(\bar{x}) = p \uhrc \bar{x}$ and $p_2(\bar{y}) = p \uhrc \bar{y}$.
By Lemma~\ref{strong unboundedness implies large distance},
\[
p(\bar{x}, \bar{y}) \models \dist(\bar{x}, \bar{y}) > 2\lambda \ \wedge \
\forall z \big(\text{``$z$ is $\lambda$-rare''}  \rightarrow \dist(z, \bar{y}) > 2\lambda \big).
\]
By Lemma~\ref{splitting a closure type},
there is a complete $(\tau, \lambda)$-{\em neighbourhood type} $p_2^+(\bar{y})$ such that 
$p_2(\bar{y}) \models p_2^+(\bar{y})$ and 
for all $n \in \mbbN^+$, for all $\bar{a} \in (B_n)^{|\bar{x}|}$ and $\bar{b} \in (B_n)^{|\bar{y}|}$,
\begin{align}\label{another equation 1}
&\text{$\mcB_n \models p(\bar{a}, \bar{b})$ if and only if $\mcB_n \models p_1(\bar{a}) \wedge p^+_2(\bar{b})$ and
$\dist(\bar{c}\bar{a}, \bar{b}) > 2\lambda$,} \\ 
&\text{where $\bar{c}$ enumerates all $\lambda$-rare elements in $\mcB_n$.} \nonumber
\end{align}
Since $p(\bar{x}, \bar{y})$ is strongly $\bar{y}$-unbounded it follows that $p_2(\bar{y})$ is strongly unbounded
and as $p_2(\bar{y}) \models p_2^+(\bar{y})$ 
it follows that $p_2^+(\bar{y})$ is strongly unbounded.
Lemma~\ref{strongly unbounded neighbourhood types tend to infinity}
now gives that 
\begin{equation}\label{yet another limit to infinity}
\lim_{n\to\infty} |p_2^+(\mcB_n)| = \infty.
\end{equation}
Since (by Assumption~\ref{properties of the base structures}) the degree of every $\mcB_n$  is bounded by $\Delta$ 
it now follows from~(\ref{yet another limit to infinity}) 
that for all sufficiently large $n$, 
if $\mcB_n \models p_1(\bar{a})$ and $\bar{c}$ enumerates all $\lambda$-rare elements of $\mcB_n$, then
there is $\bar{b} \in p_2^+(\mcB_n)$ such that $\dist(\bar{a}\bar{c}, \bar{b}) > 2\lambda$ and, 
by~(\ref{another equation 1}), we get
$\mcB_n \models p(\bar{a}, \bar{b})$, so $p(\bar{a}, \mcB_n) \neq \es$.
And of course, if $\mcB_n \not\models p_1(\bar{a})$ then $p(\bar{a}, \mcB_n) = \es$.
\hfill $\square$

\begin{lem}\label{an x-type that determines if y exists}
Let $p(\bar{x}, \bar{y})$ be a complete $(\tau, \lambda)$-closure type for some $\lambda \in \mbbN$.
Then there is $\xi \in \mbbN$ such that if $q(\bar{x})$ is a complete $(\tau, \xi)$-closure type, 
then either
\begin{enumerate}
\item for all sufficiently large $n$, if $\mcB_n \models q(\bar{a})$ then $p(\bar{a}, \mcB_n) = \es$, or
\item for all sufficiently large $n$, if $\mcB_n \models q(\bar{a})$ then $p(\bar{a}, \mcB_n) \neq \es$.
\end{enumerate}
Moreover, if $p$ is strongly $\bar{y}$-unbounded then we can let $\xi = \lambda$.
\end{lem}

\noindent
{\bf Proof.}
Let $p(\bar{x}, \bar{y})$ be a complete $(\tau, \lambda)$-closure type.
By Lemma~\ref{a type can be divided in a bounded and strongly unbounded type},
we may assume that $\bar{y} = \bar{u}\bar{v}$, $p \uhrc \bar{x}\bar{u}$ is $\bar{u}$-bounded,
and $p$ is strongly $\bar{v}$-unbounded.
Let $r(\bar{x}, \bar{u}) = p \uhrc \bar{x}\bar{u}$, so $r$ is $\bar{u}$-bounded.
From Lemma~\ref{boundedness and closure}
it follows that there is $\gamma \in \mbbN$ such that (for all $n$) if $\bar{a} \in (B_n)^{|\bar{x}|}$ and 
$\bar{c} \in r(\bar{a}, \mcB_n)$, then $\rng(\bar{c}) \subseteq C_\gamma^{\mcB_n}(\bar{a})$.
Let $\xi = \lambda + \gamma$.
Let $q(\bar{x})$ be a complete $(\tau, \xi)$-closure type.
We need to show that either~(1) or~(2) holds.

By the choice of $\xi$ it follows that either
\begin{enumerate}
\item[(3)] for all $n$, if $\mcB_n \models q(\bar{a})$ then $r(\bar{a}, \mcB_n) = \es$, or
\item[(4)] for all $n$, if $\mcB_n \models q(\bar{a})$ then $r(\bar{a}, \mcB_n) \neq \es$.
\end{enumerate}
If (3) holds then clearly, for all $n$, if $\mcB_n \models q(\bar{a})$ then $p(\bar{a}, \mcB_n) = \es$, so~(1) holds.

Suppose that (4) holds. Also suppose that $\mcB_n \models q(\bar{a})$.
Then (by~(4)) $r(\bar{a}, \mcB_n) \neq \es$ so 
$\mcB_n \models r(\bar{a}, \bar{c})$ for some $\bar{c}$.
Now it follows from
Lemma~\ref{for strongly y-unbounded the x-part determines if y exists}
that if $n$ is sufficiently large then $p(\bar{a}, \bar{c}, \mcB_n) \neq \es$.
Hence~(2) holds.

Regarding the ``moreover'' part: If $p$ is strongly $\bar{y}$-unbounded, then 
let $r(\bar{x}) = p \uhrc \bar{x}$, so $r(\bar{x})$ is a complete $(\tau, \lambda)$-closure type.
By Lemma~\ref{for strongly y-unbounded the x-part determines if y exists}, 
for all sufficiently large $n$, $\mcB_n \models r(\bar{a})$
if and only if $p(\bar{a}, \mcB_n) \neq \es$. 
It follows that if $q(\bar{x})$ is a complete $(\tau, \lambda)$-closure type which is equivalent
to $r(\bar{x})$, then, for all sufficiently large $n$, $\mcB_n \models q(\bar{a})$ implies $p(\bar{a}, \mcB_n) \neq \es$;
and if $q$ is not equivalent to $r$ then $\mcB_n \models q(\bar{a})$ implies $p(\bar{a}, \mcB_n) = \es$.
\hfill $\square$

\begin{lem}\label{all tau-neighbourhood types have roughly the same cardinality}
Let $\lambda \in \mbbN$ and  let $p(\bar{x}, \bar{y})$ be a 
strongly $\bar{y}$-unbounded complete $(\tau, \lambda)$-closure type
such that $\dim_{\bar{y}}(p) = 1$.
Then there is a constant $K \in \mbbN$ such that
the following holds for all $n$:
If $\bar{a}, \bar{a}' \in (B_n)^{|\bar{x}|}$ and both $\bar{a}$ and $\bar{a}'$ satisfy $p \uhrc \bar{x}$ (in $\mcB_n$), then
\[
|p(\bar{a}, \mcB_n)| - K \ \leq \ |p(\bar{a}', \mcB_n)| \ \leq \ |p(\bar{a}, \mcB_n)| + K.
\]
\end{lem}

\noindent
{\bf Proof.}
Let $p(\bar{x}, \bar{y})$ be as assumed in the lemma.
By Lemma~\ref{strong unboundedness implies large distance},
\[
p(\bar{x}, \bar{y}) \models \dist(\bar{x}, \bar{y}) > 2\lambda \ \wedge \
\forall z \big(\text{``$z$ is $\lambda$-rare''}  \rightarrow \dist(z, \bar{y}) > 2\lambda \big).
\]
Let $p_1(\bar{x}) = p \uhrc \bar{x}$ and $p_2(\bar{y}) = p \uhrc \bar{y}$.
By Lemma~\ref{splitting a closure type},
there is a complete $(\tau, \lambda)$-{\em neighbourhood type} $p_2^+(\bar{y})$ such that 
$p_2(\bar{y}) \models p_2^+(\bar{y})$ and 
for all $n \in \mbbN^+$, for all $\bar{a} \in (B_n)^{|\bar{x}|}$ and $\bar{b} \in (B_n)^{|\bar{y}|}$,
\begin{align}\label{another equation 2}
&\text{$\mcB_n \models p(\bar{a}, \bar{b})$ if and only if $\mcB_n \models p_1(\bar{a}) \wedge p^+_2(\bar{b})$ and
$\dist(\bar{c}\bar{a}, \bar{b}) > 2\lambda$,} \\ 
&\text{where $\bar{c}$ enumerates all $\lambda$-rare elements in $\mcB_n$.} \nonumber
\end{align}
Suppose that $\mcB_n \models p_1(\bar{a}) \wedge p_1(\bar{a}')$, 
so by Lemma~\ref{for strongly y-unbounded the x-part determines if y exists}
we have $p(\bar{a}, \mcB_n) \neq \es$ and $p(\bar{a}', \mcB_n) \neq \es$.
Let $\bar{b} \in p(\bar{a}, \mcB_n)$, so in particular $\mcB_n \models p_2(\bar{b})$,
and hence $\mcB_n \models p_2^+(\bar{b})$.
Let $\bar{c}$ enumerate all $\lambda$-rare elements in $\mcB_n$.
If $\dist(\bar{a}'\bar{c}, \ \bar{b}) > 2\lambda$ then,
by~(\ref{another equation 2}), 
$\bar{b} \in p(\bar{a}', \mcB_n)$.

Now suppose that $\dist(\bar{a}'\bar{c}, \ \bar{b}) \leq 2\lambda$.
Then there is $b \in \rng(\bar{b})$ such that $b \in N_{2\lambda}^{\mcB_n}(\bar{a}'\bar{c})$.
Since $\dim_{\bar{y}}(p) = 1$ it follows that $\rng(\bar{b}) \subseteq N_{2\lambda|\bar{y}|}^{\mcB_n}(\bar{a}'\bar{c})$.
By Lemma~\ref{definability of m-rareness}
there is a constant $m \in \mbbN$ such that for all $n$, $\mcB_n$ has at most $m$ $\lambda$-rare elements.
By Remark~\ref{remark on size of closures} 
there is a constant $K \in \mbbN$ such that for all $n$ and $\bar{d} \in (B_n)^{|\bar{x}|+m}$,
$|N_{2\lambda|\bar{y}|}^{\mcB_n}(\bar{d})| \leq K$. 
Hence the number of choices of $\bar{b} \in (B_n)^{|\bar{y}|}$ such that $\dist(\bar{a}'\bar{c}, \bar{b}) \leq 2\lambda$
is at most $K^{|\bar{y}|}$ and we get
$|p(\bar{a}, \mcB_n)| - K^{|\bar{y}|} \ \leq \ |p(\bar{a}', \mcB_n)|$.
Since we can switch the roles of $\bar{a}$ and $\bar{a}'$ we get
\[
|p(\bar{a}, \mcB_n)| - K^{|\bar{y}|} \ \leq \ |p(\bar{a}', \mcB_n)| \ \leq \ |p(\bar{a}, \mcB_n)| + K^{|\bar{y}|}. \qquad \square
\]

\section{Expansions of the base structures and probability distributions}\label{Probability distributions}

\noindent
In this section we define the kind of formalism, $PLA^*$-network, that we will then use to
define a probability distribution on the set of expansions of  a base structure to a larger signature $\sigma \supseteq \tau$.
The section ends with an example that illustrates what kind of distributions can be defined by $PLA^*$-networks

{\bf \em In this and the remaining sections we assume the following: $\tau \subseteq \sigma$ are finite relational signatures.
$\mbB = (\mcB_n : n \in \mbbN^+)$ is a sequence of finite $\tau$-structures that
satisfy Assumption~\ref{properties of the base structures}.
For each $n$, $\mbW_n$ is the set of all $\sigma$-structures that expand $\mcB_n$ 
(i.e. $\mbW_n = \{\mcA : \text{ $\mcA$ is a $\sigma$-structure and $\mcA \uhrc \tau = \mcB_n$}\}$).}

\begin{defin}\label{definition of PLA-network} {\rm
(i) A {\bf \em $PLA^*(\sigma)$-network based on $\tau$} is specified by the following two parts:
\begin{enumerate}
\item A DAG $\mbbG$ with vertex set $\sigma  \setminus  \tau$.

\item To each relation symbol $R \in \sigma \setminus \tau$ a formula 
$\theta_R(\bar{x}) \in L(\mr{par}(R) \cup \tau)$ 
is associated where $|\bar{x}|$ equals the arity of $R$ and
$\mr{par}(R)$ is the set of parents of $R$ in the DAG $\mbbG$.
We call $\theta_R$ the {\bf \em probability formula associated to $R$} by the $PLA^*(\sigma)$-network.
\end{enumerate}
We will denote a $PLA^*(\sigma)$-network by the same symbol (usually $\mbbG$, possibly with a sub or superscript)
as its underlying DAG.\\
(ii) Let $\mbbG$ denote a $PLA^*(\sigma)$-network based on $\tau$, 
let $\tau \subseteq \sigma' \subseteq \sigma$, and suppose that for
every $R \in \sigma'$, $\mr{par}(R) \subseteq \sigma'$. 
Then the $PLA^*(\sigma')$-network specified by the induced subgraph,
with vertex set $\sigma' \setminus \tau$, of the underlying DAG of $\mbbG$ 
and the probability formulas $\theta_R$ for all $R \in \sigma' \setminus \tau$ will be called the
{\bf \em $PLA^*(\sigma')$-subnetwork of $\mbbG$ induced by $\sigma'$}.
}\end{defin}

\begin{rem}\label{remark on when sigma equals tau} {\rm
Note that if $\sigma = \tau$ then the underlying DAG of a $PLA^*(\sigma)$-network based on $\tau$
is empty (as $\sigma \setminus \tau = \es$). 
So in this case ($\sigma = \tau$) there is a unique $PLA^*(\sigma)$-network based on $\tau$.
When $\sigma = \tau$ this unique $PLA^*(\sigma)$-network based on $\tau$ will be considered in the base case of an
inductive argument that will follow in 
sections~\ref{convergence and balance} --~\ref{Finding the balance in the inductive step}.
}\end{rem}

\noindent
{\bf  \em From now on let $\mbbG$ be a $PLA^*(\sigma)$-network based on $\tau$.}

\begin{defin}\label{the probability distribution induced by an PLA-network} {\rm 
(i) If $\sigma = \tau$  then $\mbbP_n$, the {\bf \em probability distribution on $\mbW_n$ induced by $\mbbG$}, 
is the unique probability distribution on (the singleton set) $\mbW_n$.\\
(ii) Now suppose that $\tau$ is a proper subset of $\sigma$ and suppose that
for each $R \in \sigma$,  its arity is denoted by $k_R$ and the
probability formula
corresponding to $R$ is denoted by $\theta_R(\bar{x})$ where $|\bar{x}| = k_R$.
Suppose that the underlying DAG of $\mbbG$ has mp-rank $\rho$.
For each $0 \leq r \leq \rho$ let $\mbbG_r$ be the subnetwork which is
induced by $\sigma_r = \{R \in \sigma : \mr{mp}(R) \leq r\}$ and note that $\mbbG_\rho = \mbbG$.
Also let $\sigma_{-1} = \tau$ and let $\mbbP^{-1}_n$ be the unique probability distribution on $\mbW^{-1}_n = \{\mcB_n\}$.
By induction on $r$ we define, for every $r = 0, 1, \ldots, \rho$, a probability distribution $\mbbP^r_n$ on the set 
\[
\mbW^r_n = \{\mcA : \text{ $\mcA$ is a $\sigma_r$-structure that expands $\mcB_n$}\}
\]
as follows:
For every $\mcA \in \mbW^r_n$, let $\mcA' = \mcA \uhrc \sigma_{r-1}$ and
\[
\mbbP^r_n(\mcA) \ = \ \mbbP^{r-1}_n(\mcA') 
\prod_{R \in \sigma_r \setminus \sigma_{r-1}} \ \prod_{\bar{a} \in R^\mcA} 
\mcA' \big(\theta_R(\bar{a})\big) \ 
\prod_{\bar{a} \in \ (B_n)^{^{k_R}} \ \setminus \ R^\mcA} \big( 1 - \mcA' \big(\theta_R(\bar{a})\big) \big).
\]
Finally we let $\mbbP_n = \mbbP^\rho_n$ and note that $\mbW_n = \mbW^\rho_n$, so 
$\mbbP_n$ is a probability distribution on $\mbW_n$ which we call
{\bf \em the probability distribution on $\mbW_n$ induced by $\mbbG$}.
We also call $\mbbP = (\mbbP_n : n \in \mbbN^+)$
{\bf \em the sequence of probability distributions induced by $\mbbG$}.
}\end{defin}

\begin{rem}\label{remark on the definition of probability distribution} {\rm
From
Lemma~\ref{truth values only depend on the relation symbols used}
it follows straightforwardly that the probability distribution $\mbbP_n$ (on $\mbW_n$)
in Definition~\ref{the probability distribution induced by an PLA-network}
could more simply have been defined as 
\[
\mbbP_n(\mcA) = 
\prod_{R \in \sigma \setminus \tau} \ \prod_{\bar{a} \in R^\mcA} 
\mcA\big(\theta_R(\bar{a})\big) \ 
\prod_{\bar{a} \in \ (B_n)^{^{k_R}} \ \setminus \ R^\mcA} \big( 1 - \mcA\big(\theta_R(\bar{a})\big) \big)
\]
for every $\mcA \in \mbW_n$.
However, the proof of the main results will use the inductive definition with respect to the mp-rank described in
Definition~\ref{the probability distribution induced by an PLA-network},
so therefore we can as well use that formulation as the definition.
}\end{rem}

\begin{notation}\label{notation E and P(formula)}{\rm
(i) For any formula $\varphi(\bar{x}) \in PLA^*(\sigma)$, all $n$, and all $\bar{a} \in (B_n)^{|\bar{x}|}$, 
let 
\[
\mbE_n^{\varphi(\bar{a})} = \{\mcA \in \mbW_n : \mcA(\varphi(\bar{a})) = 1 \}.
\]
(ii) If $\mbbP_n$ is a probability distribution on $\mbW_n$,  $\varphi(\bar{x}) \in PLA^*(\sigma)$,
and $\bar{a} \in (B_n)^{|\bar{x}|}$, then
\[
\mbbP_n(\varphi(\bar{a})) = \mbbP_n\big(\mbE_n^{\varphi(\bar{a})}\big).
\]
}\end{notation}

\noindent
Then following lemma is a basic consequence of 
Definition~\ref{the probability distribution induced by an PLA-network}
of $\mbbP^r_n$ and $\mbbP_n$.

\begin{lem}\label{conditioning on A'}
Let $\rho$, $\sigma_r$, $\mbW^r_n$, $\mbbP^r_n$ and $\mbbP_n$ be as in 
Definition~\ref{the probability distribution induced by an PLA-network}.\\
(i) Let $r \in \{0, \ldots, \rho\}$, $R \in \sigma_r \setminus \sigma_{r-1}$, 
$n \in \mbbN^+$, $\bar{a} \in (B_n)^{k_R}$ (where $k_R$ is the arity of $R$), and $\mcA' \in \mbW^{r-1}_n$.
Then
\[
\mbbP^r_n\big( \{ \mcA \in \mbW^r_n : \mcA \models R(\bar{a}) \} \ | \ 
\{ \mcA \in \mbW_n : \mcA \uhrc \sigma_{r-1} = \mcA' \}\big) \  = \ 
\mcA'(\theta_R(\bar{a})).
\]
(ii) Let $R_1, \ldots, R_t \in \sigma_r \setminus \sigma_{r-1}$
(where we allow that $R_i = R_j$ even if $i \neq j$), 
$n \in \mbbN^+$, $\bar{a}_i \in (B_n)^{k_i}$ for $i = 1, \ldots, t$,
where $k_i$ is the arity of $R_i$, and $\mcA' \in \mbW^{r-1}_n$.
Suppose that 
for $i = 1, \ldots, t$, $\varphi_i(\bar{x})$ is a literal in which $R_i$ occurs, 
and if $i \neq j$ then $\bar{a}_i \neq \bar{a}_j$ or $R_i \neq R_j$.
Using the probability distribution $\mbbP_n$,
the event $\mbE_n^{\varphi_1(\bar{a})}$ is 
independent of the event $\bigcap_{i = 2}^t \mbE_n^{\varphi_i(\bar{a})}$, {\rm conditioned} on the event 
$\{\mcA \in \mbW_n : \mcA \uhrc \sigma_{r-1} = \mcA'\}$.
\end{lem}

\begin{exam}\label{examples of PLA-networks}{\rm
Recall that if $\mbbG$ is a $PLA^*(\sigma)$-network based on $\tau$, then $\mbbG$ consists of a DAG, also
denoted $\mbbG$, with vertex set $\sigma \setminus \tau$, and for every $R \in \sigma \setminus \tau$,
a $PLA^*(\mr{par}(R))$-formula $\theta_R$ where $\mr{par}(R)$ is the set of parents of $R$ in the DAG.
Here we give some examples of what such $\theta_R$ can express (without explicitly writing out the corresponding
$PLA^*(\mr{par}(R))$-formula).

Suppose that $R, Q, E \in \sigma \setminus \tau$ where $R$ and $Q$ are unary and $E$ binary.
Depending on the example we assume that $Q \in \mr{par}(R)$ or $E \in \mr{par}(R)$.
Recall that according to 
Definition~\ref{the probability distribution induced by an PLA-network}
the role of $\theta_R(x)$ is to express the probability that $x$ satisfies $R$.
Then $\theta_R(x)$ can (for example) express any of the following values/probabilities
where $\lambda$ and $k$ are some fixed (but arbitrary) natural numbers:
\begin{enumerate}
\item ``The proportion of elements in the $\lambda$-closure (or $\lambda$-neighbourhood) of $x$
that satisfy $Q$.'' 

\item ``The average of the proportion of elements in the $\lambda$-closure of $y$
that satisfy $Q$, as $y$ ranges over all elements in the domain with the same $(\tau, \lambda)$-closure type
(or $(\tau, \lambda)$-neighbourhood type) as $x$''.

\item ``The $k^{th}$ approximation of the PageRank of $x$ (where links are represented by $E$)
computed (only) on the $\lambda$-closure
(or $\lambda$-neighbourhood) of $x$''. 
(Example~\ref{example of page rank} explains what the PageRank is and shows how its $k^{th}$ approximation 
is expressed by a $PLA^*(\{E\})$-formula.)

\item ``The $k^{th}$ approximation of the PageRank of $x$ (where links are represented by $E$) computed on the whole domain.''
\end{enumerate}

\noindent
Now suppose that $R$ is a binary relation symbol.
Let $\lambda \in \mbbN$, $\beta \in (0, 1)$ and $\alpha_0, \ldots, \alpha_\lambda \in [0, 1]$.
Then $\theta_R(x, y)$ can (for example) express the following value/pro\-ba\-bility:
\begin{enumerate}
\item[(5)] ``$\alpha_d$ if for some $d \in \{0, \ldots, \lambda\}$, the distance between $x$ and $y$ is $d$,
and otherwise $n^{-\beta}$ where $n$ is the size of the domain.''
\end{enumerate}
Note that the above examples make sense for every sequence $\mbB = (\mcB_n : n \in \mbbN^+)$
of base structures that satisfies Assumption~\ref{properties of the base structures}, so in particular they make sense for 
Examples~\ref{example with empty tau} --~\ref{example of trees}.
In the case when $\tau = \es$ (that is, the case of Example~\ref{example with empty tau} above) more examples 
of what $\theta_R$ can express are given in \cite{KW2}.
}\end{exam}

\section{Convergence, balance, a base case, and an induction hypothesis}\label{convergence and balance}

\noindent
We are now ready to define the concepts of {\em convergent pairs of formulas} and {\em balanced triples of formulas}.
Recall Notation~\ref{notation E and P(formula)}.
Informally speaking, a pair $(\varphi(\bar{x}), \psi(\bar{x}))$ is convergent if 
$\mbbP_n(\mbE_n^{\varphi(\bar{a})} \ | \ \mbE_n^{\psi(\bar{a})})$ converges as $n\to\infty$,
and a triple $(\varphi(\bar{x}, \bar{y}), \psi(\bar{x}, \bar{y}), \chi(\bar{x}))$ is balanced if there is 
$\alpha \in [0, 1]$ such that if $n$ is large enough then the probability that $\mcA \in \mbW_n$ has
the following property is high:
if $\mcA(\chi(\bar{a})) = 1$ then
the proportion of $\bar{b}$ such that $\mcA(\varphi(\bar{a}, \bar{b}))) = 1$ among the $\bar{b}$ such that 
$\mcA(\psi(\bar{a}, \bar{b})) = 1$ is close to $\alpha$.
The goal is to prove that certain pairs of closure types are convergent and that certain triples
of closure types are balanced.
In Section~\ref{Asymptotic elimination of aggregation functions}
we will use the results about balanced triples to ``asymptotically eliminate'' continuous (or admissible) aggregation
functions.

In this section we prove a result, 
Lemma~\ref{balanced bounded triples},
about balanced triples of {\em bounded} closure types.
But we also want to consider {\em un}bounded closure types.
The result about balanced triples of unbounded closure types will use a result about
convergent pairs of closure types.
We will use induction on the maximal path rank of the $PLA^*(\sigma)$-network $\mbbG$ (based on $\tau$)
which induces  $\mbbP_n$ to prove the results about convergent pairs and balanced triples
of unbounded closure types.
The base case of the induction
is taken care of in this section by
Lemma~\ref{base case}.
Then the induction hypothesis is stated as 
Assumption~\ref{induction hypothesis}.
The induction step concerning convergent pairs of closure types is carried out in 
Section~\ref{Proving convergence in the inductive step}.
The induction step concerning balanced triples of closure types is carried out in 
Section~\ref{Finding the balance in the inductive step}.

{\em We adopt the assumptions made in 
Section~\ref{Probability distributions}.
So in particular, $\mbW_n$ (for $n \in \mbbN^+$) is the set of expansions to $\sigma$ of the base structure 
$\mcB_n$ and $\mbbP_n$ is the probability distribution on $\mbW_n$ induced by $\mbbG$ which is a 
$PLA(\sigma)$-network based on $\tau$.}

\begin{defin}{\rm
Let $\varphi(\bar{x}), \psi(\bar{x}) \in PLA^*(\sigma)$ and $\alpha \in [0, 1]$.\\
(i) We say that $(\varphi, \psi)$ {\bf \em converges to $\alpha$
(with respect to $(\mbB, \mbbG)$)} if for all $\varepsilon > 0$
and all sufficiently large $n$, if $\bar{a} \in (B_n)^{|\bar{x}|}$ and 
$\mbbP_n\big(\mbE_n^{\psi(\bar{a})}\big) > 0$ then
$
\mbbP_n\big(\mbE_n^{\varphi(\bar{a})} \ \big| \ \mbE_n^{\psi(\bar{a})} \big) 
\in [\alpha - \varepsilon, \alpha + \varepsilon].
$
 We say that $(\varphi, \psi)$ {\bf \em converges (with respect to $(\mbB, \mbbG)$)} if it converges
to $\alpha$ (with respect to $(\mbB, \mbbG)$) for some $\alpha$.\\
(iii) We say that $(\varphi, \psi)$ is {\bf \em eventually constant (with respect to $(\mbB, \mbbG)$)} 
if for some $\alpha$ and all sufficiently large $n$ we have
$
\mbbP_n\big(\mbE_n^{\varphi(\bar{a})} \ \big| \ \mbE_n^{\psi(\bar{a})} \big) = \alpha
$
whenever $\bar{a} \in (B_n)^{|\bar{x}|}$ and 
$\mbbP_n\big(\mbE_n^{\psi(\bar{a})}\big) > 0$.
}\end{defin}

\begin{rem}\label{remark on definition of convergence}{\rm
It is immediate from the definition of convergence of $(\varphi, \psi)$
that if $\varphi(\bar{x}) \wedge \psi(\bar{x})$ is not cofinally satisfiable,
then $(\varphi, \psi)$ is eventually constant with $\alpha = 0$ as the eventually
constant conditional probability, and hence $(\varphi, \psi)$ converges to 0.
So the interesting case is of course when $\varphi(\bar{x}) \wedge \psi(\bar{x})$ is cofinally satisfiable.
}\end{rem}

\begin{defin}{\rm
Let $\varphi(\bar{x}, \bar{y}), \psi(\bar{x}, \bar{y})$ and $\chi(\bar{x})$ be $\sigma$-formulas.\\
(i) Let $\alpha \in [0, 1]$, $\varepsilon > 0$ and let $\mcA$ be a finite $\sigma$-structure.
The triple $(\varphi, \psi, \chi)$ is called {\bf \em $(\alpha, \varepsilon)$-balanced in $\mcA$} if
whenever $\bar{a} \in A^{|\bar{x}|}$ and $\mcA(\chi(\bar{a})) = 1$, then
\[
(\alpha - \varepsilon)|\psi(\bar{a}, \mcA)| \leq
|\varphi(\bar{a}, \mcA) \cap \psi(\bar{a}, \mcA)| \leq (\alpha + \varepsilon)|\psi(\bar{a}, \mcA)|.
\]
(ii) Let $\alpha \in [0, 1]$. 
The triple $(\varphi, \psi, \chi)$ is {\bf \em $\alpha$-balanced (with respect to $(\mbB, \mbbG)$)}
if for all $\varepsilon > 0$, if 
\[
\mbX_n^\varepsilon = 
\big\{\mcA \in \mbW_n : \text{ $(\varphi, \psi, \chi)$ is $(\alpha, \varepsilon)$-balanced in $\mcA$} \big\}
\]
then $\lim_{n\to\infty} \mbbP_n\big(\mbX_n^\varepsilon\big) = 1$.
The triple $(\varphi, \psi, \chi)$ is {\bf \em balanced (with respect to $(\mbB, \mbbG)$)} if, for some $\alpha \in [0, 1]$, it is
$\alpha$-balanced with respect to $(\mbB, \mbbG)$.
If, in addition, $\alpha > 0$ then we call $(\varphi, \psi, \chi)$ {\bf \em positively balanced (with respect to $(\mbB, \mbbG)$)}.
}\end{defin}

\begin{lem}\label{balance follows from inconsistency}
Let $\varphi(\bar{x}, \bar{y}), \psi(\bar{x}, \bar{y}), \chi(\bar{x}) \in PLA^*(\sigma)$.
If $\varphi \wedge \psi \wedge \chi$ is not cofinally satisfiable,
then $(\varphi, \psi, \chi)$ is 0-balanced with respect to $(\mbB, \mbbG)$.
\end{lem}

\noindent
{\bf Proof.} 
Suppose that $\varphi \wedge \psi \wedge \chi$ is not cofinally satisfiable.
It suffices to show that, for all sufficiently large $n$, if $\mcA \in \mbW_n$, $\mcA \models q(\bar{a})$,
and $\mcA \models \psi(\bar{a}, \bar{b})$, then $\mcA \not\models \varphi(\bar{a}, \bar{b})$.
But this is immediate from the assumption that  $\varphi \wedge \psi \wedge \chi$ is not cofinally satisfiable.
\hfill $\square$

\medskip
\noindent
The following lemma is a direct consequence of the definition of balanced triples.

\begin{lem}\label{on stronger conditioning in the balance}
Let $\varphi(\bar{x}, \bar{y}), \psi(\bar{x}, \bar{y}), \chi(\bar{x}) \in PLA^*(\sigma)$
be 0/1-valued formulas
and suppose that $(\varphi, \psi, \chi)$ is $\alpha$-balanced with respect to $(\mbB, \mbbG)$.
If $\theta(\bar{x}) \in PLA^*(\sigma)$ is 0/1-valued and $\theta(\bar{x}) \models \chi(\bar{x}))$, 
then $(\varphi, \psi, \theta)$ is $\alpha$-balanced
with respect to $(\mbB, \mbbG)$.
\end{lem}

\begin{rem}\label{remark on eventual constancy and balance} {\rm
With respect to part~(ii) of 
Theorem~\ref{general asymptotic elimination}
we note the following.
Let $p(\bar{x}, \bar{y})$ be a $(\sigma, \lambda)$-closure type, 
$p_\tau(\bar{x}, \bar{y})$ a $(\tau, \gamma)$-closure type, and
$q(\bar{x})$ a $(\sigma, \xi)$-closure type, for some $\lambda, \gamma, \xi \in \mbbN$.
Observe that if $(p, p_\tau)$ is eventually constant and
$\mbbP_n(\mbE_n^{p(\bar{a}, \bar{b})} \ | \ \mbE_n^{p_\tau(\bar{a}, \bar{b})} \cap \mbE_n^{q(\bar{a})}) = 0$
whenever $n$ is large enough, $\bar{a} \in (B_n)^{|\bar{x}|}$, and $\bar{b} \in (B_n)^{|\bar{y}|}$,
then $p(\bar{a}, \mcA) \cap p_\tau(\bar{a}, \mcA) = \es$ whenever 
$n$ is large enough, $\mcA \in \mbW_n$, $\mbbP_n(\mcA) > 0$, 
$\bar{a} \in (B_n)^{|\bar{x}|}$, $\bar{b} \in (B_n)^{|\bar{y}|}$,
and $\mcA(q(\bar{a})) = 1$.
}\end{rem}

\begin{lem}\label{balanced bounded triples}\label{balanced bounded triples}
Suppose that $\lambda_1, \lambda_2 \in \mbbN$, 
$p_1(\bar{x}, \bar{y})$ is a $(\sigma, \lambda_1)$-closure type, and
$p_2(\bar{x}, \bar{y})$ is a $(\sigma, \lambda_2)$-closure type such that $p_2 \uhrc \tau$ is 
cofinally satisfied and $\bar{y}$-bounded.

Then there is $\xi \in \mbbN$ such that for every complete $(\sigma, \xi)$-closure type $q(\bar{x})$
there is $\alpha \in [0, 1]$ such that for all $n$, all $\mcA \in \mbW_n$, 
and all $\bar{a} \in (B_n)^{|\bar{x}|}$, if $\mcA \models q(\bar{a})$ then
\[
|p_1(\bar{a}, \mcA) \cap p_2(\bar{a}, \mcB_n)| = \alpha |p_2(\bar{a}, \mcB_n)|.
\]
In particular, the triple $(p_1, p_2, q)$ is $\alpha$-balanced with respect to $(\mbB, \mbbG)$.
\end{lem}

\noindent
{\bf Proof.}
Let $p_1(\bar{x}, \bar{y})$ and $p_2(\bar{x}, \bar{y})$ be as assumed in the lemma.
Let $p_\tau = p_2 \uhrc \tau$.
Since $p_\tau$ is cofinally satifiable and $\bar{y}$-bounded it follows from
Lemma~\ref{boundedness and closure}
that there is $\gamma \in \mbbN$ such that for all $n$ and all $\bar{a} \in (B_n)^{|\bar{x}|}$,
if $\bar{b} \in p_\tau(\bar{a}, \mcB_n)$ then $\rng(\bar{b}) \subseteq C_{\gamma}^{\mcB_n}(\bar{a})$.
Let $\xi = \max\{\lambda_1, \lambda_2\} + \gamma$ and let $q(\bar{x})$ be a complete $(\sigma, \xi)$-closure type.

It follows that
if $\mcA \in \mbW_n$, $\mcA' \in \mbW_m$, $\mcA \models q(\bar{a})$, and $\mcA' \models q(\bar{a}')$,
then there is an isomorphism from $\mcA \uhrc C_{\xi}^{\mcB_n}(\bar{a})$ to
$\mcA' \uhrc C_{\xi}^{\mcB_m}(\bar{a}')$ which maps $\bar{a}$ to $\bar{a}'$.
By the choice of $\gamma$, we have $\rng(\bar{b}) \subseteq C_\gamma^{\mcB_n}(\bar{a})$
whenever $\mcA \in \mbW_n$ and $\mcA \models p_\tau(\bar{a}, \bar{b})$ 
(so in particular if $\mcA \models p_2(\bar{a}, \bar{b})$).

Hence, if $\mcA \in \mbW_n$ and $\mcA \models q(\bar{a})$, then $q$ alone determines 
\begin{enumerate}
\item the number, say $s$, of $\bar{b}$ such that $\mcB_n \models p_2(\bar{a}, \bar{b})$, and
\item the number, say $t$, of $\bar{b}$ such that $\mcB_n \models p_1(\bar{a}, \bar{b}) \wedge p_2(\bar{a}, \bar{b})$.
\end{enumerate}
If $s = 0$ let $\alpha = 0$, and otherwise let $\alpha = t/s$.
Then we have
\[
|p_1(\bar{a}, \mcA) \cap p_2(\bar{a}, \mcB_n)| = \alpha |p_2(\bar{a}, \mcB_n)|.
\]
\hfill $\square$

\begin{lem}\label{base case} {\bf (Base case of the induction)}
Suppose that $\sigma = \tau$ (so $\mbW_n = \{\mcB_n\}$ for all $n$).\\
(i) Suppose that $\lambda, \gamma \in \mbbN$, $p_1(\bar{x})$ is a complete $(\tau, \lambda)$-closure type and
$p_2(\bar{x})$ is a complete $(\tau, \lambda + \gamma)$-closure type.
Then $(p_1, p_2)$ converges with respect to $(\mbB, \mbbG)$.\\
(ii) A triple $(p_1, p_2, q)$ is balanced with respect to $(\mbB, \mbbG)$ if the following hold:
\begin{enumerate}
\item[(a)] For some $\lambda \in \mbbN$, $p_1(\bar{x}, \bar{y})$ is a complete $(\tau, \lambda)$-closure type,
\item[(b)] for some $\gamma \in \mbbN$, $p_2(\bar{x}, \bar{y})$ is a complete $(\tau, \lambda + \gamma)$-closure
type, and
\item[(c)] $q(\bar{x})$ is a complete $(\tau, \lambda)$-closure type.
\end{enumerate}
\end{lem}

\noindent
{\bf Proof.}
(i) Let $p_1(\bar{x})$  and $p_2(\bar{x})$ be as assumed.
The either $p_2(\bar{x}) \models p_1(\bar{x})$ or $p_1(\bar{x}) \wedge p_2(\bar{x})$ is not satisfiable.
In the first case it follows immediately from the definition of convergence that $(p_1, p_2)$ 
converges to 1 with respect to $(\mbB, \mbbG)$.
In the second case it follows immediately  from the definition of convergence that $(p_1, p_2)$ converges to 0
with respect to  $(\mbB, \mbbG)$.

(ii) Suppose that $p_1(\bar{x}, \bar{y})$ is a complete $(\tau, \lambda)$-closure type, 
$p_2(\bar{x}, \bar{y})$ is a complete $(\tau, \lambda + \gamma)$-closure
type and $q(\bar{x})$ is a complete $(\tau, \lambda)$-closure type.
If $p_2(\bar{x}, \bar{y}) \models p_1(\bar{x}, \bar{y}) \wedge q(\bar{x})$ 
is not cofinally satisfiable, then
Lemma~\ref{balance follows from inconsistency}
implies that $(p_1, p_2, q)$ is 0-balanced with respect to $(\mbB, \mbbG)$.

Now suppose that
$p_2(\bar{x}, \bar{y}) \wedge p_1(\bar{x}, \bar{y}) \wedge q(\bar{x})$ is cofinally satisfiable.
Then $p_2(\bar{x}, \bar{y}) \models  p_1(\bar{x}, \bar{y}) \wedge q(\bar{x})$ because otherwise 
we would be in the previous case. Now it follows immediately from the definition of balanced triples
that, for all $n \in \mbbN^+$ and $\varepsilon > 0$, $(p_1, p_2, q)$ is $(1, \varepsilon)$-balanced in $\mcB_n$.
As we assume that $\sigma = \tau$ we have $\mbW_n = \{\mcB_n\}$ so it follows immediately fro
the definition of balanced triples that $(p_1, p_2, q)$ is 1-balanced with respect to $(\mbB, \mbbG)$.
\hfill $\square$

\medskip
\noindent
As has been said above we will prove results about convergent pairs and balanced triples (of formulas)
by induction on the mp-rank of the $PLA^*(\sigma)$-network $\mbbG$. 
For this purpose we introduce some notation that will be used from now on.

\begin{defin}\label{definition of sigma'}{\rm
(i) Let $\rho$ be the mp-rank of $\mbbG$ (that is, the mp-rank of the underlying DAG of $\mbbG$),
where we make the convention that if $\sigma = \tau$, so that the underlying DAG has empty vertex set,
then the mp-rank of $\mbbG$ is $-1$.\\
(ii) Let 
$
\sigma' = \tau \cup \{ R \in \sigma \setminus \tau : \ \mr{mp}(R) < \rho \}.
$\\
(iii) For all $n \in \mbbN^+$ let $\mbW'_n$ be the set of all $\sigma'$-structures with domain $B_n$ that expand $\mcB_n$.\\
(iv) Let $\mbbG'$ denote the $PLA^*(\sigma')$-subnetwork of $\mbbG$ induced by $\sigma'$, and for each $n \in \mbbN^+$
let $\mbbP'_n$ be the probability distribution on $\mbW'_n$ which is induced by $\mbbG'$.
}\end{defin}

\noindent
Recall the definitions of {\em $(\sigma, \lambda)$-closure types} and 
{\em $(\sigma, \lambda)$-basic formulas} from Section~\ref{the base sequence}
as they will be used in the induction hypothesis below.
Observe that if $\sigma = \tau$ then $\sigma' = \tau$.
The following assumption will be the induction hypothesis that we will use
in Sections~\ref{Proving convergence in the inductive step}
and~\ref{Finding the balance in the inductive step}
The assumption below holds if $\sigma = \tau$ (the base case) as 
commented on in some more detail in 
Remark~\ref{remark on induction hypothesis and base case}.

\begin{assump}\label{induction hypothesis}{\bf (Induction hypothesis)} {\rm 
We assume that there are $\kappa, \kappa' \in \mbbN$ such that the following hold:
\begin{enumerate}
\item For every $R \in \sigma \setminus \tau$ there is a $(\sigma', \kappa)$-basic formula 
$\chi_R(\bar{x}) \in PLA^+(\sigma')$ such that 
$\chi_R(\bar{x})$ and $\theta_R(\bar{x})$ are asymptotically equivalent with respect to 
$\mbbP' = (\mbbP'_n : n \in \mbW'_n)$, 
where $\theta_R \in PLA^+(\sigma')$ is the aggregation formula corresponding to $R$ in $\mbbG$.

\item For all $\lambda, \gamma \in \mbbN$ such that $\gamma \geq \kappa'$, every complete
$(\sigma', \lambda)$-closure type $p'(\bar{x})$
and every complete $(\tau, \lambda + \gamma)$-closure type $p_\tau(\bar{x})$, 
$(p', p_\tau)$ converges with respect to $(\mbB, \mbbG')$.

\item A triple $(p', p_\tau, q')$ is balanced with respect to $(\mbB, \mbbG')$ if the following hold:
\begin{enumerate}
\item For some $\lambda \in \mbbN$, $p'(\bar{x}, \bar{y})$ is a complete $(\sigma', \lambda)$-closure type,

\item  for some $\gamma \geq \kappa'$, 
$p_\tau(\bar{x}, \bar{y})$ is a complete $(\tau, \lambda + \gamma)$-closure type, 

\item $q'(\bar{x})$ is a complete $(\sigma', \lambda)$-closure type, and

\item $p_\tau$ is strongly $\bar{y}$-unbounded and $\dim_{\bar{y}}(p_\tau) = 1$.
\end{enumerate}
\end{enumerate}
}\end{assump}

\begin{rem}\label{remark on induction hypothesis and base case}{\rm
Suppose in this remark that $\sigma = \tau$.
Then part~(1) of Assumption~\ref{induction hypothesis} holds automatically, as $\sigma \setminus \tau = \es$.
From $\sigma = \tau$ we also get $\sigma' = \tau$.
Note that in Lemma~\ref{base case} $\lambda$ and $\gamma$ denote arbitrary natural numbers.
It follows that if we let $\kappa = \kappa' = 0$ then parts~(2) and~(3) of 
Assumption~\ref{induction hypothesis} hold.
}\end{rem}

\noindent
The following basic lemma will sometimes be used without reference.

\begin{lem}\label{same probability in W as in W'}
Let $\varphi(\bar{x}), \psi(\bar{x}), \chi(\bar{x}) \in PLA^*(\sigma')$.\\
(i) If $\bar{a} \in (B_n)^{|\bar{x}|}$ (for any $n$)
then 
\[
\mbbP_n\big( \{\mcA \in \mbW_n : \mcA(\varphi(\bar{a})) = 1 \} \big) = 
\mbbP'_n\big( \{ \mcA \in \mbW'_n : \mcA(\varphi(\bar{a})) = 1 \} \big).
\]
(ii) $(\varphi, \psi)$ converges to $\alpha$ with respect to $(\mbB, \mbbG')$ if and only if $(\varphi, \psi)$
converges to $\alpha$ with respect to $(\mbB, \mbbG)$.\\
(iii) $(\varphi, \psi, \chi)$ is $\alpha$-balanced with respect to $(\mbB, \mbbG')$ if and only if $(\varphi, \psi, \chi)$
is $\alpha$-balanced with respect to $(\mbB, \mbbG)$.
\end{lem}

\noindent
{\bf Proof.} 
Part~(i) follows from 
Lemma~\ref{truth values only depend on the relation symbols used}
and Definition~\ref{the probability distribution induced by an PLA-network} 
of $\mbbP_n$ (and $\mbbP'_n$).
Part~(ii) follows from part~(i).
Part~(iii) follows from 
Lemma~\ref{truth values only depend on the relation symbols used}
and part~(i).
\hfill $\square$

\section{Proving convergence in the inductive step}\label{Proving convergence in the inductive step}

\noindent
We adopt all assumptions made in Section~\ref{Probability distributions}.
We also use the notation from Definition~\ref{definition of sigma'}
and we assume that the induction hypothesis stated in 
Assumption~\ref{induction hypothesis}
holds.
In this section we will prove that if we define $\kappa_1 := \kappa' + \kappa$, then 
part~(2) of Assumption~\ref{induction hypothesis}
holds if we replace $\sigma', \kappa'$ and $\mbbG'$ by $\sigma, \kappa_1$ and $\mbbG$, respectively,
and with this the induction step is completed for part~(2) of Assumption~\ref{induction hypothesis}.
We arrive at this conclusion (stated by
Proposition~\ref{convergence of conditional probabilities, part 3}
and Remark~\ref{the new kappa'})
by first proving a weaker version of 
Proposition~\ref{convergence of conditional probabilities, part 3}
and then strengthening it, in a couple of steps.

\begin{lem}\label{convergence of conditional probabilities, part 1}
Let $\lambda \geq \kappa$ and $\gamma \geq \kappa'$, 
let $p_\tau(\bar{x})$ be a complete $(\tau, \lambda + \gamma)$-closure type,
let $p'(\bar{x})$ be a complete $(\sigma', \lambda)$-closure type and 
let $p(\bar{x})$ be a $(\sigma, 0)$-closure type. 
Then $(p \wedge p', p_\tau)$  converges with respect to $(\mbB, \mbbG)$.
\end{lem}

\noindent
{\bf Proof.}
By Remark~\ref{remark on definition of convergence}, we may assume that
$p \wedge p'  \wedge p_\tau$ is cofinally satisfiable.
By Assumption~\ref{induction hypothesis},
$(p', p_\tau)$ converges to some $\beta$ with respect to $(\mbB, \mbbG')$.
It follows from 
Lemma~\ref{same probability in W as in W'}
that $(p', p_\tau)$ converges to $\beta$ with respect to $(\mbB, \mbbG)$.
This means that for all $\varepsilon > 0$,
all sufficiently large $n$ and all $\bar{a} \in (B_n)^{|\bar{x}|}$,
if $\mbbP_n\big(\mbE_n^{p_\tau(\bar{a})}\big) > 0$, then
$\mbbP_n\big(\mbE_n^{p'(\bar{a})}\ | \ \mbE_n^{p_\tau(\bar{a})}\big) \in 
[\beta - \varepsilon, \beta + \varepsilon]$.

To show that $(p \wedge p', p_\tau)$ converges 
we must show that there is $\alpha$ such that for every $\varepsilon > 0$ and all sufficiently large $n$
and $\bar{a} \in (B_n)^{|\bar{x}|}$, 
if $\mbbP_n\big(\mbE_n^{p_\tau(\bar{a})}\big) > 0$,
then $\mbbP_n\big(\mbE_n^{p(\bar{a})} \cap \mbE_n^{p'(\bar{a})} \ \big| \ \mbE_n^{p_\tau(\bar{a})}\big) 
\in [\alpha - \varepsilon, \alpha + \varepsilon]$.
Since $\mbE_n^{p(\bar{a})} \cap \mbE_n^{p'(\bar{a})} \subseteq \mbE_n^{p(\bar{a})}$ it
follows that if $\beta = 0$, then the above holds if $\alpha = 0$.
So now suppose that $\beta > 0$.

We have 
\[
\mbbP_n\big(\mbE_n^{p(\bar{a})} \cap \mbE_n^{p'(\bar{a})} \ \big| \ \mbE_n^{p_\tau(\bar{a})}\big) =
\mbbP_n\big(\mbE_n^{p(\bar{a})} \ \big| \ \mbE_n^{p'(\bar{a})} \cap \mbE_n^{p_\tau(\bar{a})}\big)
\cdot 
\mbbP_n\big(\mbE_n^{p'(\bar{a})}\ | \ \mbE_n^{p_\tau(\bar{a})}\big)
\]
so it suffices to prove that there is $\alpha$ such that, for all $\varepsilon > 0$, all sufficiently large $n$
and all $\bar{a} \in (B_n)^{|\bar{x}|}$, 
if $\mbbP_n\big(\mbE_n^{p'(\bar{a})} \cap \mbE_n^{p_\tau(\bar{a})}\big) > 0$, then
\begin{equation}\label{balance of p over p' and p-tau}
\mbbP_n\big(\mbE_n^{p(\bar{a})} \ \big| \ \mbE_n^{p'(\bar{a})} \cap \mbE_n^{p_\tau(\bar{a})} \big) 
\in [\alpha - \varepsilon, \alpha + \varepsilon].
\end{equation}
For every $R \in \sigma \setminus \sigma'$,
let $\theta_R$ be the probability formula of $\mbbG$ associated to $R$ and,
according to 
Assumption~\ref{induction hypothesis}~(1),
let $\chi_R$ be a $(\sigma', \kappa)$-basic formula which is asymptotically equivalent to $\theta_R$.
Let $\varepsilon' > 0$.
For every $R \in \sigma \setminus \sigma'$, if $R$ has arity $r$ let 
\[
\mbX_{n, \varepsilon'}^R = \{ \mcA \in \mbW_n : \text{ for all }\bar{a} \in (B_n)^r,  \ 
|\mcA(\theta_R(\bar{a})) - \mcA(\chi_R(\bar{a}))| < \varepsilon' \}.
\]
Since $\theta_R$ and $\chi_R$ are asymptotically equivalent 
for every $R \in \sigma \setminus \sigma'$ (and $\sigma$ is finite) it follows that for every $\varepsilon' > 0$
\begin{equation}\label{convergence to 1 of the Xs}
\lim_{n\to\infty} \mbbP_n\Big( \bigcap_{R \in \sigma \setminus \sigma'} \mbX_{n, \varepsilon'}^R \Big) = 1.
\end{equation}
Note that if $\mcB_n \models p_\tau(\bar{a})$ then $\mcA \models p_\tau(\bar{a})$ for all $\mcA \in \mbW_n$,
and if $\mcB_n \not\models p_\tau(\bar{a})$ then $\mcA \not\models p_\tau(\bar{a})$ for all $\mcA \in \mbW_n$.
Therefore we have $\mbbP_n\big(\mbE_n^{p_\tau(\bar{a})}\big) = 1$ or 
$\mbbP_n\big(\mbE_n^{p_\tau(\bar{a})}\big) = 0$.
It follows that if $\mbbP_n\big(\mbE_n^{p_\tau(\bar{a})}\big) > 0$ then
$\mbbP_n\big(\mbE_n^{p_\tau(\bar{a})}\big) = 1$, and as
$\mbbP_n\big(\mbE_n^{p'(\bar{a})}\ | \ \mbE_n^{p_\tau(\bar{a})}\big) \in [\beta - \varepsilon', \beta + \varepsilon']$
whenever $n$ is large enough,
it follows (from the definition of conditional probability) that 
$\mbbP_n\big(\mbE_n^{p'(\bar{a})} \cap \mbE_n^{p_\tau(\bar{a})}\big) \in [\beta - \varepsilon', \beta + \varepsilon']$
if $n$ is large enough.
This together with~(\ref{convergence to 1 of the Xs}) 
and the assumption that $\beta > 0$ 
implies that if $\mbbP_n\big(\mbE_n^{p_\tau(\bar{a})}\big) > 0$, then
\[
\mbbP_n\Big(\bigcap_{R \in \sigma \setminus \sigma'} \mbX_{n, \varepsilon'}^R \ \big| \ 
\mbE_n^{p'(\bar{a})} \cap \mbE_n^{p_\tau(\bar{a})} \Big) \ \geq \ 1 - \varepsilon'
\]
whenever $n$ is sufficiently large.
Therefore it suffices to show that there is $\alpha$ such that for every  $\varepsilon > 0$,
if $\varepsilon' > 0$ is sufficently small, then for all sufficiently large $n$ and $\bar{a} \in (B_n)^{|\bar{x}|}$,
then 
\[
\mbbP_n\big(\mbE_n^{p(\bar{a})} \ \big| \ \mbE_n^{p'(\bar{a})} \cap \mbE_n^{p_\tau(\bar{a})}
\cap \bigcap_{R \in \sigma \setminus \sigma'} \mbX_{n, \varepsilon'}^R \big) 
\in [\alpha - \varepsilon, \alpha + \varepsilon].
\]
By Lemma~\ref{basic fact about conditional probabilities},
it suffices to prove that there is $\alpha$ and $\varepsilon' > 0$ such that for all $\varepsilon > 0$,
all sufficiently large $n$ and all $\bar{a} \in (B_n)^{|\bar{x}|}$,
if $\mcA' \in \mbW'_n$, $\mcA' \models p_\tau(\bar{a}) \wedge p'(\bar{a})$, and
$|\mcA'(\theta_R(\bar{a})) - \mcA'(\chi_R(\bar{a}))| < \varepsilon'$ for all $R \in \sigma \setminus \sigma'$,
then
\[
\mbbP_n\big(\mbE_n^{p(\bar{a})} \ | \  \{\mcA \in \mbW_n : \mcA \uhrc \sigma' = \mcA'\}\big) 
\in [\alpha - \varepsilon, \alpha + \varepsilon].
\]
Let $\varepsilon' > 0$.
Suppose that  $\mcA' \in \mbW'_n$, $\bar{a} \in (B_n)^{|\bar{x}|}$, 
$\mcA' \models  p_\tau(\bar{a}) \wedge p'(\bar{a})$, and,
if $k_R$ is the arity of $R \in \sigma \setminus \sigma'$, then
$|\mcA'(\theta_R(a_1, \ldots, a_{k_R})) - \mcA'(\chi_R(a_1, \ldots, a_{k_R}))| < \varepsilon'$ 
for all $a_1, \ldots, a_{k_R} \in B_n$.
Suppose that $\bar{x} = (x_1, \ldots, x_m)$ and $\bar{a} = (a_1, \ldots, a_m) \in (B_n)^m$.
Lemma~\ref{conditioning on A'}
implies that
\begin{align*}
&\mbbP_n\big(\mbE_n^{p(\bar{a})} \ | \ \{\mcA \in \mbW_n : \mcA \uhrc \sigma' = \mcA'\}\big) \ = \\
&\prod_{\substack{R \in \sigma \setminus \sigma' \text{ and} \\ p(\bar{x}) \models R(x_{i_1}, \ldots, x_{i_{k_R}})}}
\mcA'(\theta_R(a_{i_1}, \ldots, a_{i_{k_R}}))
\prod_{\substack{R \in \sigma \setminus \sigma' \text{ and} \\ p(\bar{x}) \models \neg R(x_{i_1}, \ldots, x_{i_{k_R}})}}
\big(1 - \mcA'(\theta_R(a_{i_1}, \ldots, a_{i_{k_R}}))\big).
\end{align*}
Recall that, by the assumptions on $\mcA'$, 
for all $R \in \sigma \setminus \sigma'$ and $a_{i_1}, \ldots, a_{i_{k_R}} \in \rng(\bar{a})$ we have
\[
\big| \mcA'(\theta_R(a_{i_1}, \ldots, a_{i_{k_R}}))  - \mcA'(\chi_R(a_{i_1}, \ldots, a_{i_{k_R}})) \big| \leq \varepsilon'.
\]
It follows that
if $\varepsilon'  > 0$ is chosen sufficiently small and
\begin{equation*}
\alpha = \prod_{\substack{R \in \sigma \setminus \sigma' \\ p(\bar{x}) \models R(x_{i_1}, \ldots, x_{i_{k_R}})}}
\mcA'(\chi_R(a_{i_1}, \ldots, a_{i_{k_R}}))
\prod_{\substack{R \in \sigma \setminus \sigma' \\ \neg p(\bar{x}) \models R(x_{i_1}, \ldots, x_{i_{k_R}})}}
\big(1 - \mcA'(\chi_R(a_{i_1}, \ldots, a_{i_{k_R}}))\big)
\end{equation*}
then
\begin{equation*}
\Big| \mbbP_n\big(\mbE_n^{p(\bar{a})} \ | \ \{\mcA \in \mbW_n : \mcA \uhrc \sigma' = \mcA'\}\big)
\ - \ \alpha \Big| \ \leq \ \varepsilon.
\end{equation*}

Since $\chi_R$ is a $(\sigma', \kappa)$-basic formula if follows that the value
$\mcA'\big(\chi_R(a_{i_1}, \ldots, a_{i_{k_R}})\big)$ 
depends only on the complete $(\sigma', \kappa)$-closure type which $(a_{i_1}, \ldots, a_{i_{k_R}})$ satisfies in $\mcA'$.
But this complete $(\sigma', \kappa)$-closure type is determined by $p'(\bar{x})$ because 
$p'(\bar{x})$ is a complete $(\sigma', \lambda)$-closure type where 
(by assumption) $\lambda \geq \kappa$ and $\mcA' \models p'(\bar{a})$.
It follows that $\alpha$ depends only on $p$, $p'$ and $p_\tau$.
\hfill $\square$
\\

\noindent
The following corollary will be used in a proof in 
Section~\ref{Finding the balance in the inductive step}.

\begin{cor}\label{corollary about convergence of (p, p' and p-tau)}
For every $R \in \sigma \setminus \sigma'$
let $\theta_R$ be the probability formula of $\mbbG$ associated to $R$ and,
according to 
Assumption~\ref{induction hypothesis}~(1),
let $\chi_R$ be a $(\sigma', \kappa)$-basic formula which is asymptotically equivalent to $\theta_R$.
Also, for every $\varepsilon' > 0$, let
\[
\mbX_{n, \varepsilon'}^R = \{ \mcA \in \mbW_n : \text{ for all }\bar{a} \in (B_n)^r,  \ 
|\mcA(\theta_R(\bar{a})) - \mcA(\chi_R(\bar{a}))| < \varepsilon' \}.
\]
If $p(\bar{x}), p'(\bar{x})$ and $p_\tau(\bar{x})$ are as in
Lemma~\ref{convergence of conditional probabilities, part 1},
then there is $\alpha$ such that, for every  $\varepsilon > 0$,
if $\varepsilon' > 0$ is sufficently small, then for all sufficiently large $n$ and $\bar{a} \in (B_n)^{|\bar{x}|}$,
\[
\mbbP_n\big(\mbE_n^{p(\bar{a})} \ \big| \ \mbE_n^{p'(\bar{a})} \cap \mbE_n^{p_\tau(\bar{a})}
\cap \bigcap_{R \in \sigma \setminus \sigma'} \mbX_{n, \varepsilon'}^R \big) 
\in [\alpha - \varepsilon, \alpha + \varepsilon].
\]
\end{cor}

\noindent
{\bf Proof.}
This was proved in the proof of 
Lemma~\ref{convergence of conditional probabilities, part 1}.
\hfill $\square$

\begin{rem}\label{remark about basic probability formulas} {\rm
In this remark let us make the following stronger assumptions than conditions~(1) and~(2)
of Assumption~\ref{induction hypothesis}:
\begin{enumerate}
\item[(I)] For every $R \in \sigma \setminus \sigma'$ the corresponding probability formula $\theta_R$ is
a $(\sigma', \kappa)$-basic formula.
\item[(II)] For $\lambda, \gamma \in \mbbN$ such that $\gamma \geq \kappa'$, 
every complete $(\sigma', \lambda)$-closure type $p'(\bar{x})$,
and every complete $(\tau, \lambda + \gamma)$-closure type $p_\tau(\bar{x})$ that is consistent with $p'$,
$(p', p_\tau)$ is eventually constant (with respect to $(\mbB, \mbbG')$).
\end{enumerate}
Then the proof of Lemma~\ref{convergence of conditional probabilities, part 1}
works out if we let $\chi_R$ be the {\em same} formula as $\theta_R$ for all $R \in \sigma \setminus \sigma'$.
Moreover, in the beginning of the proof we can, by assumption (II), suppose that
$\mbbP_n\big(\mbE_n^{p'(\bar{a})}\ | \ \mbE_n^{p_\tau(\bar{a})}\big) = \beta$,
and towards the end of the proof we get the stronger conclusion that
$\mbbP_n\big(\mbE_n^{p(\bar{a})} \ | \ \{\mcA \in \mbW_n : \mcA \uhrc \sigma' = \mcA'\}\big) = \alpha$.
It follows that $(p, p' \wedge p_\tau)$ and $(p \wedge p', p_\tau)$ are eventually constant (with respect to $(\mbB, \mbbG)$).
}\end{rem}

\begin{lem}\label{convergence of conditional probabilities, part 2}
Let $\gamma \geq \kappa'$, let $p_\tau(\bar{x})$ be a complete $(\tau, \kappa + \gamma)$-closure type
and let $p(\bar{x})$ be a $(\sigma, 0)$-closure type.
Then $(p, p_\tau)$ converges with respect to $(\mbB, \mbbG)$.
\end{lem}

\noindent
{\bf Proof.}
We may assume that $p \wedge p_\tau$ is cofinally satisfiable since otherwise the lemma follows immediately.
Let $p'_1(\bar{x}), \ldots, p'_s(\bar{x})$ enumerate all, up to equivalence,
complete $(\sigma', \kappa)$-closure types $p'(\bar{x})$ such that $p_\tau \wedge p' \wedge p$
is cofinally satisfiable. 
We also assume that if $i \neq j$ then $p'_i$ is not equivalent to $p'_j$.
By Lemma~\ref{convergence of conditional probabilities, part 1}, 
for all $i$, $(p \wedge p'_i, p_\tau)$ converges to some $\alpha_i$.
For all sufficiently large $n$ and all $\bar{a} \in (B_n)^{|\bar{x}|}$ such that $\mcB_n \models p_\tau(\bar{a})$,
$\mbE_n^{p(\bar{a})}$ is the disjoint union of 
$\mbE_n^{p'_1(\bar{a}) \wedge p(\bar{a})}, \ldots, \mbE_n^{p'_s(\bar{a}) \wedge p(\bar{a})}$
so 
\[\mbbP_n\big( \mbE_n^{p(\bar{a})} \ \big| \ \mbE_n^{p_\tau(\bar{a})} \big) = 
\sum_{i=1}^s \mbbP_n\big( \mbE_n^{p'_i(\bar{a}) \wedge p(\bar{a})} \ \big| \ \mbE_n^{p_\tau(\bar{a})} \big)
\]
where for all $\varepsilon > 0$ and all large enough $n$,
\[
\mbbP_n\big( \mbE_n^{p'_i(\bar{a}) \wedge p(\bar{a})} \ \big| \ \mbE_n^{p_\tau(\bar{a})} \big) \in
[\alpha_i - \varepsilon, \alpha_i + \varepsilon]
\]
whenever $\mbbP_n\big(\mbE_n^{p_\tau(\bar{a})}\big) > 0$.
Therefore $(p, p_\tau)$ converges to $(\alpha_1 + \ldots + \alpha_s)$.
\hfill $\square$

\begin{prop}\label{convergence of conditional probabilities, part 3}
Let $\lambda \in \mbbN$, let $\gamma \geq \kappa'$, let 
$p(\bar{x})$ be a complete $(\sigma, \lambda)$-closure type, and let
$p_\tau(\bar{x})$ be a complete $(\tau, \lambda + \kappa + \gamma)$-closure type.
Then $(p, p_\tau)$ converges with respect to $(\mbB, \mbbG)$.
\end{prop}

\noindent
{\bf Proof.}
Let $p(\bar{x})$ and $p_\tau(\bar{x})$ be as in the lemma.
We may also assume that $p \wedge p_\tau$ is cofinally satisfiable because otherwise the conclusion is trivial.
Then there is a cofinally satisfiable complete $(\tau, \kappa + \gamma)$-closure type
$p^*_\tau(\bar{x}, \bar{y})$ be  such that,
for all $n$, if $\mcB_n \models p^*_\tau(\bar{a}, \bar{b})$, then $\mcB_n \models p_\tau(\bar{a})$
and $\rng(\bar{b}) = C^{\mcB_n}_\lambda(\bar{a}) \setminus \rng(\bar{a})$.
Let $p^*(\bar{x}, \bar{y})$ be a complete $(\sigma, 0)$-closure type such that
$p^*_\tau \wedge p^* \wedge p$ is consistent.
It follows that (for all $n$) if $\mcA \in \mbW_n$ and $\mcA  \models p^*(\bar{a}, \bar{b}) \wedge p_\tau(\bar{a})$,
then $\mcA \models p(\bar{a})$.
It also follows that if $\mcA \models p(\bar{a}) \wedge p_\tau(\bar{a})$ then there is $\bar{b} \in (B_n)^{|\bar{y}|}$
such that $\mcA \models p^*(\bar{a}, \bar{b}) \wedge p^*_\tau(\bar{a}, \bar{b})$.

By Lemma~\ref{convergence of conditional probabilities, part 2},
$(p^*, p^*_\tau)$ converges to some $\alpha$.
Let $n$ be arbitrary and suppose that $\mcA \in \mbW_n$ and 
$\mcA \models p^*(\bar{a}, \bar{b}) \wedge p^*_\tau(\bar{a}, \bar{b})$.
Let $\bar{b}_1, \ldots, \bar{b}_s$ be a maximal sequence (with respect to length) of
different permutations of $\bar{b}$ such that
$\mcA \models p^*(\bar{a}, \bar{b}_i) \wedge p^*_\tau(\bar{a}, \bar{b}_i)$ for all $i$,
and if $i \neq j$, $\mcA', \mcA'' \in \mbW_n$,
$\mcA' \models p^*(\bar{a}, \bar{b}_i) \wedge p^*_\tau(\bar{a}, \bar{b}_i)$ and
$\mcA'' \models p^*(\bar{a}, \bar{b}_j) \wedge p^*_\tau(\bar{a}, \bar{b}_j)$, then
$\mcA' \uhrc \rng(\bar{a}\bar{b}) \neq \mcA'' \uhrc \rng(\bar{a}\bar{b})$.
Note that $s$ depends only on $p^*$ and $p^*_\tau$.
Also observe that if $i \neq j$ then $\mbE_n^{p^*(\bar{a}, \bar{b}_i)}$ is disjoint from
$\mbE_n^{p^*(\bar{a}, \bar{b}_j)}$.
Since $(p^*, p^*_\tau)$ converges to $\alpha$ it follows that $(p, p_\tau)$ converges to $s\alpha$.
\hfill $\square$

\begin{rem}\label{the new kappa'}  {\rm
Proposition~\ref{convergence of conditional probabilities, part 3}
shows that if we define $\kappa_1 := \kappa' + \kappa$, then 
part~(2) of Assumption~\ref{induction hypothesis}
holds if we replace $\sigma', \kappa'$ and $\mbbG'$ by $\sigma, \kappa_1$ and $\mbbG$, respectively.
Thus the induction step is completed for part~(2) of Assumption~\ref{induction hypothesis}.
}\end{rem}

\begin{rem}\label{remark on eventually constant pairs in the general case} {\rm
Suppose that conditions~(I) and~(II) of 
Remark~\ref{remark about basic probability formulas} hold.
It follows from that remark and the proofs of
Lemma~\ref{convergence of conditional probabilities, part 2}
and Proposition~\ref{convergence of conditional probabilities, part 3}
that if
$\lambda \in \mbbN$, $\gamma \geq \kappa'$, 
$p(\bar{x})$ is a complete $(\sigma, \lambda)$-closure type, and 
$p_\tau(\bar{x})$ is a complete $(\tau, \lambda + \kappa + \gamma)$-closure type,
then $(p, p_\tau)$ is eventually constant with respect to $(\mbB, \mbbG)$.
}\end{rem}

\section{Proof of balance in the inductive step}\label{Finding the balance in the inductive step}

\noindent
In this section we prove that if $\kappa'$, $\sigma'$ and $\mbbG'$ from Assumption~\ref{induction hypothesis}
are replaced by $\kappa' + \kappa$, $\sigma$ and $\mbbG$, respectively,
then part~(3) of Assumption~\ref{induction hypothesis} still holds.
This follows from Lemma~\ref{balanced triples, strongly unbounded, dimension 1},
as pointed out by Remark~\ref{the new kappa' for balance}.
However, we continue to prove more general results about balanced triples
because we need these results to prove more general results about asymptotic elimination of aggregation functions
(than can be done with Lemma~\ref{balanced triples, strongly unbounded, dimension 1} alone),
and these will used to complete 
(in Section~\ref{Asymptotic elimination of aggregation functions}) 
the inductive step for part~(1) of
Assumption~\ref{induction hypothesis}.
Just as in the previous section we adopt all assumptions made in Section~\ref{Probability distributions}.
We also use the notation from Definition~\ref{definition of sigma'}
and we assume that the induction hypothesis stated in 
Assumption~\ref{induction hypothesis}
holds. It follows that we can use all results from 
Section~\ref{Proving convergence in the inductive step}
(but actually we only use
Corollary~\ref{corollary about convergence of (p, p' and p-tau)}
in this section).

\begin{defin} {\rm
If $\mcA' \in \mbW'_n$ then we let $\mbW_n^{\mcA'} = \{\mcA \in \mbW_n : \mcA \uhrc \sigma' = \mcA'\}$.
}\end{defin}

\noindent
For every $R \in \sigma \setminus \sigma'$
let $\theta_R$ denote the probability formula of $\mbbG$ associated to $R$ 
(so $\theta_R \in PLA^*(\sigma')$) and,
according to 
Assumption~\ref{induction hypothesis}~(1),-
let $\chi_R$ be a $(\sigma', \kappa)$-basic formula
(recall Definition~\ref{Definition of basic formula}) which is asymptotically equivalent to $\theta_R$.
The following lemma is the technical fundament on which the following results of this section rest.

\begin{lem}\label{simplest form of number of realizations, 1 dimension}
Let us assume the following:
\begin{enumerate}
\item $p(\bar{x}, \bar{y})$ is a $(\sigma, 0)$-closure type,
\item $p'(\bar{x}, \bar{y})$ is a complete $(\sigma', \lambda)$-closure type where $\lambda \geq \kappa$,
\item $p_\tau(\bar{x}, \bar{y})$ is a complete $(\tau, \lambda + \gamma)$-closure type where $\gamma \geq \kappa'$,
\item $\dim_{\bar{y}}(p_\tau) = 1$, and
\item whenever $R \in \sigma \setminus \sigma'$, $r$ is the arity of $R$, and
$p(\bar{x}, \bar{y}) \models R(\bar{z})$ or $p(\bar{x}, \bar{y}) \models \neg R(\bar{z})$
where $\bar{z}$ is a subsequence of $\bar{x}\bar{y}$ of length $r$,
then $\bar{z}$ contains at least one variable from $\bar{y}$.
\end{enumerate}
Furthermore, suppose that for some $n$, $\bar{a} \in (B_n)^{|\bar{x}|}$, 
$\bar{b}_1, \ldots, \bar{b}_{n_0} \in (B_n)^{|\bar{y}|}$ are different tuples,
$\mcA' \in \mbW'_n$,
$\mcA' \models p'(\bar{a}, \bar{b}_i) \wedge p_\tau(\bar{a}, \bar{b}_i)$ for all $i = 1, \ldots, n_0$, 
and, for every $R \in \sigma \setminus \sigma'$, if $r$ is the arity of $R$ and $\bar{c} \in (B_n)^r$ then
$|\mcA'(\theta_R(\bar{c})) - \mcA'(\chi_R(\bar{c}))| < \varepsilon'$.

Then there is $\gamma \in [0, 1]$, depending only on $p'$ and $p$, 
such that for every $\varepsilon > 0$ there is $c > 0$, depending only on $\varepsilon$,
such that if $\varepsilon' > 0$ is small enough and $n$ and $n_0$ are 
large enough, then
\begin{align*}
\mbbP_n\big( \big\{ \mcA \in \mbW_n : 
(\gamma - \varepsilon)n_0 \leq 
|\{ i : \mcA \models p(\bar{a}, \bar{b}_i) \}| \leq 
(\gamma + \varepsilon)n_0 \big\} \ \big| \ 
\mbW_n^{\mcA'} \big) \\ 
\geq \ 1 - e^{-c \gamma n_0}.
\end{align*}
\end{lem}

\noindent
{\bf Proof.}
We adopt the assumptions of the lemma.
Let $P = \{\bar{b}_1, \ldots, \bar{b}_{n_0}\}$.
First we prove the following claim.

\medskip

\noindent
{\bf Claim.}
{\em There are $k_0, k_1 \in \mbbN^+$, 
depending only on the sequence $(\mcB_n : n \in \mbbN^+)$ of base structures,
and a partition $P_1, \ldots, P_{k_0}$ of $P$ such that, for all $k = 1, \ldots, k_0$,
 $|P_k| \geq n_0/k_1$ and if $\bar{b}_i, \bar{b}_j \in P_k$ and $i \neq j$, 
then $\rng(\bar{b}_i) \cap \rng(\bar{b}_j) = \es$.}

\medskip

\noindent
{\bf Proof of the claim.}
By assumption, all $\bar{b}_i$ have the same length.
By Assumption~\ref{properties of the base structures},
for every $n$, the degree of $\mcB_n$
is at most $\Delta$. It follows that there is $t_0 \in \mbbN^+$,
depending only on the sequence $(\mcB_n : n \in \mbbN^+)$ of base structures, 
such that for all $c \in B_n$, $\big|N_{2(\lambda + \gamma)|\bar{b}_1|}^{\mcB_n}(c)\big| \leq t_0$.
Without loss of generality we may assume that $t_0 \geq 2$.
It follows that for every $\bar{b}_i$ we have $|N_{2(\lambda + \gamma)|\bar{b}_1|}^{\mcB_n}(\bar{b}_i)| \leq t_0|\bar{b}_1|$.
Since $p$ is strongly $\bar{y}$-unbounded and $\dim_{\bar{y}}(p) = 1$ it follows
from Lemma~\ref{another characterization of strong unboundedness}
and Definition~\ref{definition of dimension, part 2} that,
for each $i$, the distance (in $\mcB_n$) between any two elements in $\rng(\bar{b}_i)$ is at most 
$2(\lambda + \gamma)|\bar{b}_1|$.
It follows that if $\bar{b}_i$ and $\bar{b}_j$ have a common element, 
then $\rng(\bar{b}_j) \subseteq N_{2(\lambda + \gamma)|\bar{b}_1|}^{\mcB_n}(\bar{b}_i)$, 
so for each $\bar{b}_i$ there are at most $t := (t_0|\bar{b}_1|)^{|\bar{b}_1|}$ choices of $\bar{b}_j$ that has
a common element with $\bar{b}_i$.

Let $P_1 \subseteq P$ be maximal (with respect to `$\subseteq$')
such that if $\bar{b}_i, \bar{b}_j \in P_1$ and $i \neq j$, then $\rng(\bar{b}_i) \cap \rng(\bar{b}_j) = \es$.
Then $|P_1| \geq |P|/t > |P|/(2t) = n_0/(2t)$.

Now suppose that $P_1, \ldots, P_k$ are disjoint subsets of $P$ such that,
for all $l = 1, \ldots, k$, $|P_l| \geq n_0/(2t)^{k+1}$ and if 
$\bar{b}_i, \bar{b}_j \in P_l$ and $i \neq j$, then $\rng(\bar{b}_i) \cap \rng(\bar{b}_j) = \es$.
If $P_1 \cup \ldots \cup P_k = P$ then we are done and let $k_0$ as in the claim be equal to $k$
and let $k_1$ be equal to $(2t)^{k+1}$.

So suppose that $P_1 \cup \ldots \cup P_k \neq P$.
Let $P'_{k+1} \subseteq P \setminus (P_1 \cup \ldots \cup P_k)$ be maximal such that
if $\bar{b}_i, \bar{b}_j \in P_k$ and $i \neq j$, then $\rng(\bar{b}_i) \cap \rng(\bar{b}_j) = \es$.
If $|P'_{k+1}| \geq n_0/(2t)^{k+2}$ then let $P_{k+1} = P'_{k+1}$.

Now suppose that $|P'_{k+1}| < n_0/(2t)^{k+2}$.
At most $t|P'_{k+1}|$ tuples in $P_k$ have a common element with some tuple in $P'_{k+1}$.
Also $t|P'_{k+1}| < \frac{n_0}{2(2t)^{k+1}}$.
Since $|P_k| \geq n_0/(2t)^{k+1}$ we can choose $Q \subseteq P_k$ such that 
$|Q| \geq n_0/(2t)^{k+2}$, $|P_k \setminus Q| \geq n_0/(2t)^{k+2}$ 
and all pairs of different tuples in $P'_{k+1} \cup Q$ do not have a common element.
Now let $P_k^* = P_k \setminus Q$ and $P_{k+1}^* = P'_{k+1} \cup Q$.
Then $|P_k^*| \geq n_0 / (2t)^{k+2}$ and $|P_{k+1}^*| \geq n_0 / (2t)^{k+2}$.
Now let $P_k \supseteq P_k^*$ be maximal such that $P_k \subseteq P \setminus (P_1 \cup \ldots \cup P_{k-1} \cup P_{k+1})$
and every pair of different tuples in $P_k$ have no common element.
Finally, let $P_{k+1} \supseteq P_{k+1}^*$ be maximal such that 
$P_{k+1} \subseteq P \setminus (P_1 \cup \ldots \cup P_k)$
and every pair of different tuples in $P_{k+1}$ have no common element.

It remains to show that for some $k$, depending only on 
$(\mcB_n : n \in \mbbN^+)$, we will have $P = P_1 \cup \ldots \cup P_k$.
Suppose that $P \neq P_1 \cup \ldots \cup P_k$, so we can choose
$\bar{b} \in P \setminus (P_1 \cup \ldots \cup P_k)$.
Let also $\bar{b}_l \in P_l$ for $l = 1, \ldots, k$.
Recall that by the choice of $t$ we can have $\rng(\bar{b}_l) \cap \rng(\bar{b}_{l'}) \neq \es$ 
for at most $t$ of the indices $l' = 1, \ldots, k$.
Hence we can choose a subsequence $\bar{b}_{l_i} \in P_{l_i}$, for $i = 1, \ldots, \lfloor \frac{k}{t} \rfloor$,
such that $\rng(\bar{b}_{l_i}) \cap \rng(\bar{b}_{l_j}) = \es$ if $i \neq j$.
By the maximality of each $P_l$ it follows that $\rng(\bar{b}) \cap \rng(\bar{b}_{l_i}) \neq \es$ for all 
$i = 1, \ldots, \lfloor \frac{k}{t} \rfloor$.
From the choice of $t$ it follows that
$\lfloor \frac{k}{t} \rfloor \leq t$ and hence $k \leq t^2$.
Thus there is $k \leq t^2$ such that after $k$ iterations we have $P = P_1 \cup \ldots \cup P_k$.
Hence we can let $k_1 = (2t)^{t^2+1}$.
\hfill $\square$

\medskip

\noindent
By the claim, there are $k_0$ and $k_1$ that depend only on
$\mbB = (\mcB_n : n \in \mbbN^+)$ and a
partition $P_1, \ldots, P_{k_0}$ of $P$ such that, for all $k = 1, \ldots, k_0$,
$|P_k| \geq n_0/k_1$ and if $\bar{b}_i, \bar{b}_j \in P_k$ and $i \neq j$, 
then $\rng(\bar{b}_i) \cap \rng(\bar{b}_j) = \es$.

Fix any $l \in \{1, \ldots, {k_0}\}$.
Recall the assumption that whenever $R \in \sigma \setminus \sigma'$, $k_R$ is the arity of $R$, and
$p(\bar{x}, \bar{y}) \models R(\bar{z})$ or $p(\bar{x}, \bar{y}) \models \neg R(\bar{z})$
where $\bar{z}$ is a subsequence of $\bar{x}\bar{y}$ of length $k_R$,
then $\bar{z}$ contains at least one variable from $\bar{y}$.
This property of $p$ together with 
Lemma~\ref{conditioning on A'}
implies that if $\bar{b}_1, \ldots, \bar{b}_m \in P_l$ are distinct sequences
then, conditioned on $\mbW_n^{\mcA'}$, for every $i = 1, \ldots, m$ the event
$\mbE_n^{p(\bar{a}, \bar{b}_i)}$ is independent from all events
$\mbE_n^{p(\bar{a}, \bar{b}_j)}$ where  $j \in \{1, \ldots, m\} \setminus \{i\}$.
By assumption, for each $\bar{b}_i \in P_l$ we have $\mcA' \models p'(\bar{a}, \bar{b}_i) \wedge p_\tau(\bar{a}, \bar{b}_i)$,
so $\mbW_n^{\mcA'} \subseteq \mbE_n^{p'(\bar{a}, \bar{b}_i)} \cap \mbE_n^{p_\tau(\bar{a}, \bar{b}_i)}$.
For every $R \in \sigma \setminus \sigma'$, of arity $r$ say, and $\varepsilon' > 0$, let
\[
\mbX_{n, \varepsilon'}^R = \{ \mcA \in \mbW_n : \text{ for all }\bar{a} \in (B_n)^r,  \ 
|\mcA(\theta_R(\bar{a})) - \mcA(\chi_R(\bar{a}))| < \varepsilon' \}.
\]
By the assumption on $\mcA'$ in the present lemma
and by Lemma~\ref{truth values only depend on the relation symbols used}
we have $\mbW_n^{\mcA'} \subseteq \bigcap_{R \in \sigma \setminus \sigma'} \mbX_{n, \varepsilon'}^R$.
Now Corollary~\ref{corollary about convergence of (p, p' and p-tau)}
implies that there is $\gamma$ (not depending on $\mcA'$, $n$ or $n_0$)
such that, for all  $i = 1, \ldots, n_0$ and $\varepsilon > 0$,
if $\varepsilon' > 0$ is sufficently small, then for all sufficiently large $n$ 
\[
\mbbP_n\big(\mbE_n^{p(\bar{a}, \bar{b}_i)} \ | \ \mbW_n^{\mcA'}\big) \in [\gamma - \varepsilon, \gamma + \varepsilon].
\]
Recall that $n_0 = |P|$ and $|P_l| \geq n_0/k_1$.
Let
\[
\mbY_l^\varepsilon = \big\{ \mcA \in \mbW_n^{\mcA'} : 
(1 - \varepsilon)(\gamma - \varepsilon)n_0 \leq 
|\{ \bar{b}_i \in P_l : \mcA \models p(\bar{a}, \bar{b}_i) \}| 
\leq (1 + \varepsilon)(\gamma + \varepsilon)n_0 \big\}.
\]
By Corollary~\ref{independent bernoulli trials, second version}, 
there is $c_0 > 0$, depending only on $\varepsilon$, such that if $n_0$ and $n$
are large enough, then 
\[
\mbbP_n\big( \mbY_l^\varepsilon \ \big| \ \mbW_n^{\mcA'} \big) 
\geq \ 1 - e^{-c_0 \gamma |P_l|} \geq 1 - e^{-c_0 \gamma n_0/k_1}.
\]
It follows that for some $c > 0$ 
\[
\mbbP_n\bigg( \bigcap_{l = 1}^{k_0} \mbY_l^\varepsilon \ \bigg| \ 
\mbW_n^{\mcA'} \bigg) 
\geq \ 1 - e^{-c \gamma n_0}
\]
if $n_0$ and $n$ are large enough.
Since  $P_1, \ldots, P_{k_0}$ is a partition of $P$ the conclusion of the lemma now follows.
\hfill $\square$

\medskip
\noindent
In the next lemma we combine Lemma~\ref{simplest form of number of realizations, 1 dimension}
with part~(3) of Assumption~\ref{induction hypothesis} to prove the first results about balanced triples.

\begin{lem}\label{balanced triples, strongly unbounded, dimension 1, 0-closure}
Suppose that the following hold:
\begin{enumerate}
\item $p(\bar{x}, \bar{y})$ is a $(\sigma, 0)$-closure type,
\item $p'(\bar{x}, \bar{y})$ is a complete $(\sigma', \gamma)$-closure type where $\gamma \geq \kappa$,
\item $p_\tau(\bar{x}, \bar{y})$ is a complete $(\tau, \xi)$-closure type where $\xi \geq \gamma + \kappa'$, and
\item $p_\tau$ is strongly $\bar{y}$-unbounded and $\dim_{\bar{y}}(p_\tau) = 1$.
\end{enumerate}
Let $q(\bar{x}) = p \uhrc \bar{x}$ and $q'(\bar{x}) = p' \uhrc \bar{x}$.
Then $(p \wedge p', p_\tau, q \wedge q')$ is balanced with respect to $(\mbB, \mbbG)$.
\end{lem}

\noindent
{\bf Proof.}
If $p \wedge p' \wedge p_\tau$ is not cofinally satisfiable then it follows from
Lemma~\ref{balance follows from inconsistency}
that $(p \wedge p', p_\tau, q \wedge q')$ is 0-balanced with respect to $(\mbB, \mbbG)$.
So we assume that $p \wedge p' \wedge p_\tau$ is cofinally satisfiable.
This implies that $p' \models p_\tau$.
For every $R \in \sigma \setminus \sigma'$
let $\theta_R$ be the probability formula of $\mbbG$ associated to $R$ and,
according to 
Assumption~\ref{induction hypothesis}~(1),
let $\chi_R$ be a $(\sigma', \kappa)$-basic formula which is asymptotically equivalent to $\theta_R$.
Also, for every $\varepsilon' > 0$, if $r$ is the arity of $R$ let
\[
\mbX_{n, \varepsilon'}^R = \{ \mcA \in \mbW_n : \text{ for all }\bar{a} \in (B_n)^r,  \ 
|\mcA(\theta_R(\bar{a})) - \mcA(\chi_R(\bar{a}))| < \varepsilon' \}.
\]
Then, by the definition of asymptotic equivalence and 
Lemma~\ref{same probability in W as in W'}, 
for every $\varepsilon' > 0$,
\[
\lim_{n\to\infty}\mbbP_n\Big(\bigcap_{R \in \sigma \setminus \sigma'} \mbX_{n, \varepsilon'}^R \Big) = 1.
\]
By Assumption~\ref{induction hypothesis},
there is $\beta$ such that
$(p', p_\tau, q')$ is $\beta$-balanced with respect to $(\mbB, \mbbG')$.
By Lemma~\ref{same probability in W as in W'},
it is also $\beta$-balanced with respect to $(\mbB, \mbbG)$.
This means that, for every $\varepsilon > 0$, if 
\[
\mbY_n^\varepsilon = 
\big\{ \mcA \in \mbW_n : \text{ $(p', p_\tau, q')$ is $(\beta, \varepsilon)$-balanced in $\mcA$} \big\},
\]
then
$\lim_{n\to\infty}\mbbP_n\big(\mbY_n^\varepsilon\big) = 1$.

We need to show that there is $\alpha \in [0, 1]$ such that, for every $\varepsilon > 0$, if
\[
\mbZ_n^{\alpha, \varepsilon} = \{ \mcA \in \mbW_n : 
\text{ $(p \wedge p', p_\tau, q \wedge q')$ is $(\alpha, \varepsilon)$-balanced in $\mcA$}\},
\]
then $\lim_{n\to\infty} \mbbP_n(\mbZ_n^{\alpha, \varepsilon}) = 1$.
If $\beta = 0$ then it is straightforward to see that the above holds if $\alpha = 0$.
So now we assume that $\beta > 0$.

Since $\lim_{n\to\infty}\mbbP_n\Big( \mbY_n^\varepsilon \cap 
\bigcap_{R \in \sigma \setminus \sigma'} \mbX_{n, \varepsilon'}^R  \Big) = 1$ for all $\varepsilon, \varepsilon' > 0$
it suffices to show that there is $\alpha$ such that for all $\varepsilon > 0$ 
there is $\varepsilon' > 0$ such that
\begin{equation}\label{equation about X-epsilon'}
\lim_{n\to\infty} \mbbP_n\big( \mbZ_n^{\alpha, \varepsilon} \ \big| \ 
\mbY_n^{\varepsilon'} \cap 
\bigcap_{R \in \sigma \setminus \sigma'} \mbX_{n, \varepsilon'}^R \big) \ = \ 1.
\end{equation}
For each $\varepsilon' > 0$, let 
\[
\mbU_n^{\varepsilon'} = \Big\{ \mcA \uhrc \sigma' :  
\mcA \in \mbY_n^{\varepsilon'} \cap \bigcap_{R \in \sigma \setminus \sigma'} \mbX_{n, \varepsilon'}^R \Big\}
\]
so $\mbU_n^{\varepsilon'} \subseteq \mbW'_n$ and
$\mbY_n^{\varepsilon'} \cap \bigcap_{R \in \sigma \setminus \sigma'} \mbX_{n, \varepsilon'}^R$ is the 
disjoint union of all $\mbW_n^{\mcA'}$ as $\mcA'$ ranges over $\mbU_n^{\varepsilon'}$.

By Definition~\ref{definition of unbounded formula}
of strongly $\bar{y}$-unbounded closure types 
there is $f_{p_\tau} : \mbbN \to \mbbR$ such that,
for all $K > 0$,  $\lim_{n\to\infty}(f_{p_\tau}(n) - K\ln(n)) = \infty$,
and for all $n$ and $\bar{a} \in (B_n)^{|\bar{x}|}$,
if $p_\tau(\bar{a}, \mcB_n) \neq \es$ then $|p_\tau(\bar{a}, \mcB_n)| \geq f_{p_\tau}(n)$.
By Lemma~\ref{basic fact about conditional probabilities}, 
to prove~(\ref{equation about X-epsilon'}) it suffices to show that there are $\alpha$ and $d > 0$ such that for 
all $\varepsilon > 0$, there is $\varepsilon' > 0$ such that if $n$ sufficently large,
and $\mcA' \in \mbU_n^{\varepsilon'}$,
then 
\begin{equation}\label{the lower limit conditioned on A'}
\mbbP_n\big(\mbZ_n^{\alpha, \varepsilon} \ \big| \  \mbW_n^{\mcA'} \big)  \ \geq \ 
1 - e^{-d \, f_{p_\tau}(n)}.
\end{equation}
From the definition of a $(\sigma, 0)$-closure type it follows that 
$p(\bar{x}, \bar{y})$ is a conjunction of $\sigma$-literals with variables from $\rng(\bar{x}\bar{y})$.
Let $p^*(\bar{x}, \bar{y})$ be the conjunction of all literals $\varphi(\bar{z})$
of the form $R(\bar{z})$ or $\neg R(\bar{z})$,
where $R \in \sigma \setminus \sigma'$, such that
$p(\bar{x}, \bar{y}) \models \varphi(\bar{z})$
and $\rng(\bar{z}) \cap \rng(\bar{y}) \neq \es$.
Then, for all $n$, $\mcA \in \mbW_n$, $\bar{a} \in (B_n)^{|\bar{x}|}$, and $\bar{b} \in (B_n)^{|\bar{y}|}$:
\begin{equation}\label{the special property of p*}
\text{if $\mcA \models p'(\bar{a}, \bar{b}) \wedge q(\bar{a})$, then 
$\mcA \models p(\bar{a}, \bar{b})$ if and only if $\mcA \models p^*(\bar{a}, \bar{b})$.} 
\end{equation}

Fix $n$, $\varepsilon' > 0$ and $\mcA' \in \mbU_n^{\varepsilon'}$. 
We aim at proving that there are $\alpha$ and $d > 0$ (which are independent of $\mcA'$) such that, for all $\varepsilon > 0$, 
(\ref{the lower limit conditioned on A'}) holds
if $\varepsilon'$ is small enough and $n$ large enough.
For every $\bar{a} \in (B_n)^{|\bar{x}|}$, let
\[
B'_{\bar{a}} = \{ \bar{b} \in (B_n)^{|\bar{y}|} : 
\mcA' \models p'(\bar{a}, \bar{b}) \}.
\]
Since $\mcA' \in \mbU_n^{\varepsilon'}$ it follows that $(p', p_\tau, q')$ is $(\beta, \varepsilon')$-balanced in $\mcA'$.
So
\begin{equation}\label{bounds on B'}
\text{ if $\mcA' \models q'(\bar{a})$, then } 
(\beta - \varepsilon')|p_\tau(\bar{a}, \mcB_n)| \ \leq \ |B'_{\bar{a}}| \ \leq \
(\beta + \varepsilon')|p_\tau(\bar{a}, \mcB_n)|.
\end{equation}
For all $\bar{a} \in (B_n)^{|\bar{x}|}$ and $\mcA \in \mbW_n^{\mcA'}$ define
\[
B_{\bar{a}, \mcA} = \{ \bar{b} \in B'_{\bar{a}} : \mcA \models p^*(\bar{a}, \bar{b}) \wedge p'(\bar{a}, \bar{b}) \}.
\]
By Lemma~\ref{simplest form of number of realizations, 1 dimension}
there are $\gamma \in [0, 1]$ and $c > 0$, where $c$ depends only on $\varepsilon'$,
such that, conditioned on $\mcA \in \mbW_n^{\mcA'}$ and $\mcA' \models q'(\bar{a})$, the probability that 
\begin{equation}\label{bounds on B}
(\gamma - \varepsilon')|B'_{\bar{a}}| \leq |B_{\bar{a}, \mcA}| \leq (\gamma + \varepsilon')|B'_{\bar{a}}|
\end{equation}
is at least $1 - e^{-c \gamma |B'_{\bar{a}}|} \geq 1 - e^{-c \gamma |p_\tau(\bar{a}, \mcB_n)|}
\geq 1 -  e^{-c \gamma f_{p_\tau}(n)}$ for all sufficiently large $n$
(because $|B'_{\bar{a}}| \geq (\beta - \varepsilon')|p_\tau(\bar{a}, \mcB_n)| \geq 
(\beta - \varepsilon')f_{p_\tau}(n) \to \infty$).
It follows from~(\ref{bounds on B'}) and~(\ref{bounds on B}) that,
conditioned on $\mcA \in \mbW_n^{\mcA'}$ and $\mcA' \models q'(\bar{a})$, the probability that 
\[
(\beta - \varepsilon')(\gamma - \varepsilon') |p_\tau(\bar{a}, \mcB_n)| \ \leq \ 
|B_{\bar{a}, \mcA}| \ \leq \ 
(\beta + \varepsilon')(\gamma + \varepsilon') |p_\tau(\bar{a}, \mcB_n)|
\]
is at least $1 - e^{-c \gamma f_{p_\tau}(n)}$.
Let $\alpha = \beta \gamma$.
It follows that, conditioned on $\mcA \in \mbW_n^{\mcA'}$ and $\mcA' \models q'(\bar{a})$, 
the probability that
\begin{equation}\label{final inequalities about B}
(\alpha - 3\varepsilon') |p_\tau(\bar{a}, \mcB_n)| \ \leq \ 
|B_{\bar{a}, \mcA}| \ \leq \ 
(\alpha + 3\varepsilon') |p_\tau(\bar{a}, \mcB_n)|
\end{equation}
is at least $1 - e^{-c \gamma f_{p_\tau}(n)}$.

By Assumption~\ref{properties of the base structures}, 
there is a polynomial $P(x)$ such that for all $n$, $|B_n| \leq P(n)$.
Then $\big|(B_n)^{|\bar{x}|}\big| \leq P(n)^{|\bar{x}|}$ where the right side is a function in $n$
which can be expressed by a polynomial.
From the assumption about $f_{p_\tau}$ and
Lemma~\ref{consequence of 4b} 
it follows that there is $d > 0$ such that, for all sufficently large $n$,
$P(n)^{|x|} e^{-c \gamma f_{p_\tau}(n)} \leq e^{-d f_{p_\tau}(n)}$.
It follows that, 
conditioned on $\mcA \in \mbW_n^{\mcA'}$, the probability that~(\ref{final inequalities about B})
holds for {\em every} $\bar{a} \in (B_n)^{|\bar{x}|}$ such that $\mcA' \models q'(\bar{a})$
is at least $1 - e^{-d  f_{p_\tau}(n)}$.

Observe that if $\mcA \in \mbW_n^{\mcA'}$ and $\mcA \models q'(\bar{a}) \wedge q(\bar{a})$ (hence $\mcA' \models q'(\bar{a})$),
then it follows from~(\ref{the special property of p*}) that
\[
B_{\bar{a}, \mcA} = \{\bar{b} \in (B_n)^{|\bar{y}|} : 
\mcA \models p(\bar{a}, \bar{b}) \wedge p'(\bar{a}, \bar{b}) \}.
\]
Therefore we have proved that~(\ref{the lower limit conditioned on A'}) holds if $\varepsilon = 3\varepsilon'$
and $n$ is sufficiently large. This completes the proof of the lemma.
\hfill $\square$

\bigskip
\noindent
Our next result generalizes the previous lemma to the case when
closure types ``speak about'' larger closures.

\begin{lem}\label{balanced triples, strongly unbounded, dimension 1}
Suppose that the following hold:
\begin{enumerate}
\item $p(\bar{x}, \bar{y})$ is a complete $(\sigma, \lambda)$-closure type,
\item $p_\tau(\bar{x}, \bar{y})$ is a complete $(\tau, \lambda + \gamma)$-closure type where 
$\gamma \geq \max(\lambda, \kappa + \kappa')$, and
\item $p_\tau$ is strongly $\bar{y}$-unbounded and $\dim_{\bar{y}}(p_\tau) = 1$.
\end{enumerate}
If $q(\bar{x})$ is a complete $(\sigma, \lambda)$-closure type
then $(p, p_\tau, q)$ is balanced with respect to $(\mbB, \mbbG)$.
\end{lem}

\noindent
{\bf Proof.}
We assume that $p \wedge p_\tau \wedge q$ is cofinally satisfiable since otherwise the conclusion of the lemma follows from
Lemma~\ref{balance follows from inconsistency}.
Hence $q(\bar{x})$ is equivalent to $p \uhrc \bar{x}$ and we may as well assume that $q(\bar{x}) = p \uhrc \bar{x}$.

From the definition of closure types it follows that there is a 
$(\tau, \lambda + \gamma)$-neighbourhood type $p'_\tau(\bar{u}, \bar{x}, \bar{y})$ and a 
$(\tau, \gamma)$-neighbourhood type $p''_\tau(\bar{u}, \bar{x}, \bar{v}, \bar{y}, \bar{w})$ such that
the following equivalences hold in $\mcB_n$ for all $n$:
\begin{align*}
&p_\tau(\bar{x}, \bar{y}) \Longleftrightarrow \\
&\exists \bar{u} \big( \text{``$\bar{u}$ enumerates all $(\lambda + \gamma)$-rare elements''} \wedge 
p'_\tau(\bar{u}, \bar{x}, \bar{y}) \big) \\ 
&\Longleftrightarrow  \\
&\exists \bar{u}, \bar{v}, \bar{w} \big(
\text{``$\bar{u}$ enumerates all $(\lambda + \gamma)$-rare elements''} \wedge \\
&\text{``$\bar{v}$ enumerates all elements in the $\lambda$-neighbourhood of $\bar{u}\bar{x}$''} \wedge \\
&\text{``$\bar{w}$ enumerates all elements in the $\lambda$-neighbourhood of $\bar{y}$''} \wedge \\
&p''_\tau(\bar{u},\bar{x}, \bar{v}, \bar{y}, \bar{w})\big). 
\end{align*}
Since we assume that $\gamma \geq \max(\lambda, \kappa + \kappa')$ it follows
(from the definition of closure type) that 
there is a $(\tau, \gamma)$-{\em closure} type $p^+_\tau(\bar{u}, \bar{x}, \bar{v}, \bar{y}, \bar{w})$
such that (in all $\mcB_n$):
\begin{align}\label{the formula that p-+-tau is equivalent to}
&p^+_\tau(\bar{u}, \bar{x}, \bar{v}, \bar{y}, \bar{w}) \Longleftrightarrow \\
&p''_\tau(\bar{u},\bar{x}, \bar{v}, \bar{y}, \bar{w}) \wedge \nonumber\\
&\text{``$\bar{u}$ contains all $\gamma$-rare elements''} \wedge \nonumber\\
&\text{``$\bar{v}$ enumerates all elements in the $\lambda$-neighbourhood of $\bar{u}\bar{x}$''} \wedge  \nonumber \\
&\text{``$\bar{w}$ enumerates all elements in the $\lambda$-neighbourhood of $\bar{y}$''}.  \nonumber
\end{align}
Thus we get the following equivalence (in all $\mcB_n$):
\begin{align}\label{equivalent condition to p-tau}
&p_\tau(\bar{x}, \bar{y}) \Longleftrightarrow \\
&\exists \bar{u}, \bar{v}, \bar{w} \big( \text{``$\bar{u}$ enumerates all $(\lambda + \gamma)$-rare elements''} \wedge 
\nonumber \\
&p^+_\tau(\bar{u}, \bar{x}, \bar{v}, \bar{y}, \bar{w}) \big). \nonumber
\end{align}
Now we want to show that $p^+_\tau(\bar{u}, \bar{x}, \bar{v}, \bar{y}, \bar{w}) $ is 
strongly $\bar{y}\bar{w}$-unbounded and $\dim_{\bar{y}\bar{w}}(p^+_\tau)$, because then we will be able to
use Lemma~\ref{balanced triples, strongly unbounded, dimension 1, 0-closure}
to a triple of closure types that involves $p^+_\tau$.
For this we use the characterization of strong unboundedness in 
Lemma~\ref{another characterization of strong unboundedness}.
Recall the assumption that $p_\tau$ is strongly $\bar{y}$-unbounded with $\bar{y}$-dimension 1.
It follows from Lemma~\ref{another characterization of strong unboundedness}
and the definition of dimension that $y_i \sim_{p_\tau} y_j$ for all
$y_i, y_j \in \rng(\bar{y})$ and $x_i \not\sim_{p_\tau} y_j$ for all $x_i \in \rng(\bar{x})$
and $y_j \in \rng(\bar{y})$. 
From~(\ref{the formula that p-+-tau is equivalent to})
it follows that
for every $w_j \in \rng(\bar{w})$ there is $y_i \in \rng(\bar{y})$ such that 
$y_i \sim_{p^+_\tau} w_j$. 
Let $y_i, y_j \in \rng(\bar{y})$. 
So $y_i \sim_{p_\tau} y_j$ and from~(\ref{the formula that p-+-tau is equivalent to}) 
it follows that there 
are $w_k, w_l \in \rng(\bar{w})$ such that $y_i \sim_{p^+_\tau} w_k \sim_{p^+_\tau} w_l \sim_{p^+_\tau} y_j$
and, as $p^+_ \tau$ is an equivalence relation,  $y_i \sim_{p^+_\tau} y_j$.
Now it follows that $w_i \sim_{p^+_\tau} w_j$ for all $w_i, w_j \in \rng(\bar{w})$.

If we would have $x_i \sim_{p^+_\tau} y_j$ for some $x_i \in \rng(\bar{x})$ and $y_j \in \rng(\bar{y})$
then it would follow from the definitions of $\sim_{p^+_\tau}$ and $\sim_{p_\tau}$ that 
$x_i \sim_{p_\tau} y_j$ which contradicts 
(by the use of Lemma~\ref{another characterization of strong unboundedness})
that $p_\tau$ is strongly $\bar{y}$-unbounded.
So $x_i \not\sim_{p_\tau} y_j$.
Using what we have proved and that $p''_\tau$, and hence $p^+_\tau$, implies that ``every $u_i \in rng(\bar{u})$
is $(\lambda + \gamma)$-rare it follows in a straightforward maner that
$z \not\sim_{p^+_\tau} z'$ for all $z \in \rng(\bar{u}\bar{x}\bar{v})$ and $z' \in \rng(\bar{y}\bar{w})$.

If there is a subsequence $\bar{z}$ of $\bar{y}\bar{w}$ such that $p^+_\tau \uhrc \bar{z}$ is bounded
then, as $y_i \sim_{p^+_\tau} w_j$ for all $y_i \in \rng(\bar{y})$ and $w_j \in \rng(\bar{w})$,
it follows that $p^+_\tau \uhrc \bar{y}\bar{w}$ is bounded and 
from~(\ref{equivalent condition to p-tau}) and~(\ref{the formula that p-+-tau is equivalent to})  
we get that $p_\tau \uhrc \bar{y}$ is bounded, contradicting
that $p_\tau$ is strongly $\bar{y}$-unbounded.
By Lemma~\ref{another characterization of strong unboundedness} we conclude that
$p^+_\tau$ is strongly $\bar{y}\bar{w}$-unbounded. 
From what we have proved about $\sim_{p^+_\tau}$ it also follows that $\dim_{\bar{y}\bar{w}}(p^+_\tau) = 1$.

Since we assume that $p \wedge p_\tau$ is cofinally satisfiable it follows, using~(\ref{equivalent condition to p-tau}),
that there is a complete $(\sigma, 0)$-closure type $p^+(\bar{u}, \bar{x}, \bar{v}, \bar{y}, \bar{w})$
such that $p^+ \wedge p \wedge p_\tau$ is cofinally satisfiable and (in all $\mcB_n$)
\begin{align}\label{p and p-tau is equivalent to p+ and p+tau}
&p(\bar{x}, \bar{y}) \wedge p_\tau(\bar{x}, \bar{y}) \Longleftrightarrow \\
&\exists \bar{u}, \bar{v}, \bar{w} \big( \text{``$\bar{u}$ enumerates all $(\lambda + \gamma)$-rare elements''} \wedge 
\nonumber \\
&p^+(\bar{u}, \bar{x}, \bar{v}, \bar{y}, \bar{w}) \wedge p^+_\tau(\bar{u}, \bar{x}, \bar{v}, \bar{y}, \bar{w}) \big). \nonumber
\end{align}
Let $q^+_\tau(\bar{u}\bar{x}\bar{v}) = p^+_\tau \uhrc \bar{u}\bar{x}\bar{v}$.
Let $p'_i(\bar{u}, \bar{x}, \bar{v}, \bar{y}, \bar{w})$, $i = 1, \ldots, s$, enumerate, up to equivalence, all 
complete $(\sigma', \kappa)$-closure types in the variables $\bar{u}, \bar{x}, \bar{v}, \bar{y}, \bar{w}$.
Let $q'_i(\bar{u}\bar{x}\bar{v}) = p'_i \uhrc \bar{u}\bar{x}\bar{v}$.
It follows that if 
$\mcA \models p^+(\bar{c}, \bar{a}, \bar{d}, \bar{b}, \bar{e}) \wedge 
p^+_\tau(\bar{c}, \bar{a}, \bar{d}, \bar{b}, \bar{e})$ then there is a unique $i$ such that
$\mcA \models p'_i(\bar{c}, \bar{a}, \bar{d}, \bar{b}, \bar{e})$.

Lemma~\ref{balanced triples, strongly unbounded, dimension 1, 0-closure}
implies that for each $i$ there is $\alpha_i$ such that
$(p^+ \wedge p_i, p^+_\tau, q^+ \wedge q'_i)$ is $\alpha_i$-balanced with respect to $(\mbB, \mbbG)$.
For all $\varepsilon > 0$ let $\mbX_n^\varepsilon$
be the set of all $\mcA \in \mbW_n$ such that for all $i$, 
$(p^+ \wedge p_i, p^+_\tau, q^+ \wedge q'_i)$ is $(\alpha_i, \varepsilon)$-balanced in $\mcA$.
Then $\lim_{n\to\infty}\mbbP_n\big(\mbX_n^\varepsilon\big) = 1$ for all $\varepsilon > 0$.

From~(\ref{the formula that p-+-tau is equivalent to}) it follows that for all $n \in \mbbN^+$
(and tuples from $B_n$ of appropriate length) 
if $\mcB_n \models p^+_\tau(\bar{c}, \bar{a}, \bar{d}, \bar{b}, \bar{e})$ then
the following conditions are equivalent:
\begin{enumerate}
\item[\textbullet] $\mcB_n \models p^+_\tau(\bar{c}, \bar{a}, \bar{d}, \bar{b}, \bar{e}')$.
\item[\textbullet] $\bar{e}'$ is a permutation of $\bar{e}$ and there is an isomorphism from
$\mcB_n \uhrc N_\gamma^{\mcB_n}(\bar{c} \bar{a} \bar{d} \bar{b} \bar{e})$
to $\mcB_n \uhrc N_\gamma^{\mcB_n}(\bar{c} \bar{a} \bar{d} \bar{b} \bar{e}')$
that maps $\bar{c} \bar{a} \bar{d} \bar{b} \bar{e}$ to $\bar{c} \bar{a} \bar{d} \bar{b} \bar{e}'$.
\end{enumerate}
It follows that there is $\xi_\tau \in \mbbN^+$, depending only on $p^+_\tau$, such that 
$|p^+_\tau(\bar{c}, \bar{a}, \bar{d}, \bar{b}, \mcB_n)|$ is either 0 or $\xi_\tau$.
Since $p^+_\tau \in PLA^*(\tau)$ the same holds for every $\mcA \in \mbW_n$.

Moreover, for all $i = 1, \ldots, s$, if $\mcA \in \mbW_n$ (for any $n$) and
$\mcA \models p^+(\bar{c}, \bar{a}, \bar{d}, \bar{b}, \bar{e}) \wedge p_i(\bar{c}, \bar{a}, \bar{d}, \bar{b}, \bar{e})
\wedge p^+_\tau(\bar{c}, \bar{a}, \bar{d}, \bar{b}, \bar{e})$, then the following are equivalent
\begin{enumerate}
\item[\textbullet] $\mcA \models p^+(\bar{c}, \bar{a}, \bar{d}, \bar{b}, \bar{e}') \wedge 
p_i(\bar{c}, \bar{a}, \bar{d}, \bar{b}, \bar{e}')
\wedge p^+_\tau(\bar{c}, \bar{a}, \bar{d}, \bar{b}, \bar{e}')$

\item[\textbullet] $\bar{e}'$ is a permuation of $\bar{e}$, there is an isomorphism 
$f$ from $\mcA \uhrc \rng(\bar{c} \bar{a} \bar{d} \bar{b} \bar{e})$ to 
$\mcA \uhrc \rng(\bar{c} \bar{a} \bar{d} \bar{b} \bar{e}')$ that maps $\bar{c} \bar{a} \bar{d} \bar{b} \bar{e}$
to $\bar{c} \bar{a} \bar{d} \bar{b} \bar{e}'$, there is an isomorphism $g$
from $(\mcA \uhrc \sigma') \uhrc N_\lambda^{\mcB_n}(\bar{c} \bar{a} \bar{d} \bar{b} \bar{e})$
to $(\mcA \uhrc \sigma') \uhrc N_\lambda^{\mcB_n}(\bar{c} \bar{a} \bar{d} \bar{b} \bar{e}')$
that extends $f$, and there is an isomorphism from
$\mcB_n \uhrc N_{\lambda + \gamma}^{\mcB_n}(\bar{c} \bar{a} \bar{d} \bar{b} \bar{e})$ to 
$\mcB_n \uhrc N_{\lambda + \gamma}^{\mcB_n}(\bar{c} \bar{a} \bar{d} \bar{b} \bar{e}')$ that extends $g$.
\end{enumerate}
It follows that there is $\xi_i \in \mbbN^+$, depending only on $p^+$, $p_i$ and $p^+_\tau$, such that
\[
\big|p^+(\bar{c}, \bar{a}, \bar{d}, \bar{b}, \mcA) \cap p_i(\bar{c}, \bar{a}, \bar{d}, \bar{b}, \mcA)
\cap p^+_\tau(\bar{c}, \bar{a}, \bar{d}, \bar{b}, \mcA)\big| 
\]
is either 0 or $\xi_i$.

After all this preparation we are ready to prove that $(p, p_\tau, q)$ is balanced.
Suppose that $\mcA \in \mbX_n^\varepsilon$, $\mcA \models q(\bar{a})$ and $p_\tau(\bar{a}, \mcA) \neq \es$.
Since $q(\bar{x}) = p \uhrc \bar{x}$ there are $\bar{c} \in (B_n)^{|\bar{u}|}$ and $\bar{d} \in (B_n)^{|\bar{v}|}$
such that $\mcA \models q^+(\bar{c}, \bar{a}, \bar{d})$, $\bar{c}$ enumerates all $(\lambda + \gamma)$-rare
elements and $\bar{d}$ enumerates $N_\gamma^{\mcB_n}(\bar{c}\bar{a})$.
Now we have 
\begin{align*}
\big|p_\tau(\bar{a}, \mcA)\big| = \frac{\big|p^+_\tau(\bar{c}, \bar{a}, \bar{d}, \mcA)\big|}{\xi_\tau}.
\end{align*}
By~(\ref{p and p-tau is equivalent to p+ and p+tau}) and the choice of $p'_i$ we have
\begin{align*}
\big|p(\bar{a}, \mcA) \cap p_\tau(\bar{a}, \mcA)\big| = 
\sum_{i=1}^s 
\frac{\big| p^+(\bar{c}, \bar{a}, \bar{d}, \mcA) \cap 
p'_i(\bar{c}, \bar{a}, \bar{d}, \mcA) \cap
p^+_\tau(\bar{c}, \bar{a}, \bar{d}, \mcA)\big|}{\xi_i}.
\end{align*}
Hence
\[
\frac{|p(\bar{a}, \mcA) \cap p_\tau(\bar{a}, \mcA)\big|}{|p_\tau(\bar{a}, \mcA)\big|} = 
\sum_{i=1}^s 
\frac{\xi_\tau  \big|p^+(\bar{c}, \bar{a}, \bar{d}, \mcA) \cap 
p'_i(\bar{c}, \bar{a}, \bar{d}, \mcA) \cap
p^+_\tau(\bar{c}, \bar{a}, \bar{d}, \mcA)\big|}
{\xi_i  \big|p^+_\tau(\bar{c}, \bar{a}, \bar{d}, \mcA)\big|}.
\]
Since $\mcA \in \mbX_n^\varepsilon$ it follows that if $p^+_\tau(\bar{c}, \bar{a}, \bar{d}, \mcA) \neq \es$,
then
\[
\alpha_i - \varepsilon \leq \frac{\big|p^+(\bar{c}, \bar{a}, \bar{d}, \mcA) \cap 
p'_i(\bar{c}, \bar{a}, \bar{d}, \mcA) \cap
p^+_\tau(\bar{c}, \bar{a}, \bar{d}, \mcA)\big|}{\big| p^+_\tau(\bar{c}, \bar{a}, \bar{d}, \mcA) \big|}
\leq \alpha_i + \varepsilon.
\]
Thus we get
\[
\sum_{i=1}^s  (\alpha_i - \varepsilon)\frac{\xi_\tau}{\xi_i} \ \leq \ 
\frac{|p(\bar{a}, \mcA) \cap p_\tau(\bar{a}, \mcA)\big|}{|p_\tau(\bar{a}, \mcA)\big|} \ \leq \ 
\sum_{i=1}^s  (\alpha_i + \varepsilon)\frac{\xi_\tau}{\xi_i}.
\]
As $\varepsilon > 0$ can be taken as small as we like it follows that if 
$\alpha = \sum_{i=1}^s  \frac{\alpha_i \xi_\tau}{\xi_i}$, then
$(p, p_\tau, q)$ is $\alpha$-balanced with respect to $(\mbB, \mbbG)$.
\hfill $\square$

\begin{rem}\label{the new kappa' for balance} {\rm
Lemma~\ref{balanced triples, strongly unbounded, dimension 1}
shows that if we let $\kappa_1 = \kappa' + \kappa$ then 
part~(3) of Assumption~\ref{induction hypothesis} is satisfied if
we replace $\sigma', \kappa'$ and $\mbbG'$ by $\sigma, \kappa_1$ and $\mbbG$, respectively.
Hence the induction step for part~(3) of Assumption~\ref{induction hypothesis} is completed.
}\end{rem}

\noindent
In spite of Remark~\ref{the new kappa' for balance} 
we continue to prove more results about balanced triples, for more general closure types,
because we need these results to prove more general results about asymptotic elimination of aggregation functions
(than can be done with Lemma~\ref{balanced triples, strongly unbounded, dimension 1}),
and these will used to complete 
(in Section~\ref{Asymptotic elimination of aggregation functions}) 
the inductive step for part~(1) of
Assumption~\ref{induction hypothesis}.
First we generalize Lemma~\ref{balanced triples, strongly unbounded, dimension 1}
to the case when $p_\tau(\bar{x}, \bar{y})$ is {\em uniformly} $\bar{y}$-unbounded.

\begin{lem}\label{balanced triples, uniformly unbounded, dimension 1}
Suppose that the following hold:
\begin{enumerate}
\item $p(\bar{x}, \bar{y})$ is a complete $(\sigma, \lambda)$-closure type,
\item $p_\tau(\bar{x}, \bar{y})$ is a complete $(\tau, \lambda + \gamma)$-closure type where 
$\gamma \geq \max(\lambda, \kappa + \kappa')$, and
\item $p_\tau$ is uniformly $\bar{y}$-unbounded and $\dim_{\bar{y}}(p_\tau) = 1$.
\end{enumerate}
There is $\xi \in \mbbN^+$ such that if 
$q^*(\bar{x})$ is a complete $(\sigma, \xi)$-closure type, 
then $(p, p_\tau, q^*)$ is balanced with respect to $(\mbB, \mbbG)$.
\end{lem}

\noindent
{\bf Proof.}
As usual we assume that $p \wedge p_\tau$ is cofinally satisfiable since otherwise we get the conclusion from 
Lemma~\ref{balance follows from inconsistency}.
By Lemma~\ref{a type can be divided in a bounded and strongly unbounded type}
we can assume that $\bar{y} = \bar{u}\bar{v}$ where $r_\tau(\bar{x}, \bar{u}) = p_\tau \uhrc \bar{x}\bar{u}$
is $\bar{u}$-bounded and $p_\tau$ is {\em strongly} $\bar{v}$-unbounded.
Since $r_\tau$ is $\bar{u}$-bounded there is $m \in \mbbN$ such that for all $n$ and all $\bar{a} \in (B_n)^{|\bar{x}|}$,
$|r_\tau(\bar{a}, \mcB_n)| \leq m$.

Let $r(\bar{x}, \bar{u}) = p \uhrc \bar{x}\bar{u}$.
According to
Lemma~\ref{boundedness and closure}
there is $\xi \in \mbbN^+$ such that if
$\mcB_n \models r_\tau(\bar{a}, \bar{b})$ then $\rng(\bar{b}) \subseteq C_\xi^{\mcB_n}(\bar{a})$.
Let $q^*(\bar{x})$ be a complete $(\sigma, \xi)$-closure type.
Then there are $s \leq t  \leq m$, depending only on $q^*$, 
such that (for every $n$) if $\mcA \in \mbW_n$ and $\mcA \models q^*(\bar{a})$,
then $|r_\tau(\bar{a}, \mcA)| = t$ and $|r(\bar{a}, \mcA) \cap r_\tau(\bar{a}, \mcA)| = s$.

Now suppose that $\mcA \in \mbW_n$, $\mcA \models q^*(\bar{a})$,
$\mcA \models r_\tau(\bar{a}, \bar{b}_i)$, for $i = 1, \ldots, t$, 
and $\mcA \models r(\bar{a}, \bar{b}_i) \wedge r_\tau(\bar{a}, \bar{b}_i)$, 
for $i = 1, \ldots, s$.
Recall that $p_\tau$ is strongly $\bar{v}$-unbounded.
As $\dim_{\bar{y}}(p) = 1$ it follows that $\dim_{\bar{v}}(p) = 1$.
Let $q(\bar{x}, \bar{u}) = p \uhrc \bar{x}\bar{u}$.
Lemma~\ref{balanced triples, strongly unbounded, dimension 1}
now implies that $(p, p_\tau, q)$ is $\alpha$-balanced with respect to $(\mbB, \mbbG)$ for some $\alpha$.
So for any $\varepsilon > 0$, if
\[
\mbX_n^\varepsilon = \big\{\mcA \in \mbW_n : \text{ $(p, p_\tau, q)$ is $(\alpha, \varepsilon)$-balanced
in $\mcA$} \big\}
\]
then $\lim_{n\to\infty}\mbbP_n\big(\mbX_n^\varepsilon\big) = 1$.

Now suppose that $\mcA \in \mbX_n^\varepsilon$, $\mcA \models q^*(\bar{a})$,
$\mcA \models r_\tau(\bar{a}, \bar{b}_i)$, for $i = 1, \ldots, t$, 
and $\mcA \models r(\bar{a}, \bar{b}_i) \wedge r_\tau(\bar{a}, \bar{b}_i)$, 
for $i = 1, \ldots, s$.
Since $\mcA \in \mbX_n^\varepsilon$ it follows that for all $i = 1, \ldots, s$,
\[
(\alpha - \varepsilon)|p_\tau(\bar{a}, \bar{b}_i, \mcA)| \leq 
|p(\bar{a}, \bar{b}_i, \mcA) \cap p_\tau(\bar{a}, \bar{b}_i, \mcA)| 
\leq (\alpha + \varepsilon)|p_\tau(\bar{a}, \bar{b}_i, \mcA)|.
\]
By the choice of the $\bar{b}_i$ we also have 
$|p(\bar{a}, \bar{b}_i, \mcA) \cap p_\tau(\bar{a}, \bar{b}_i, \mcA)| = 0$ for all $i = s+1, \ldots, t$.
As $p_\tau$ is strongly $\bar{v}$-unbounded it follows from
Lemma~\ref{all tau-neighbourhood types have roughly the same cardinality},
that there is a constant $C > 0$, depending only on $\mbB$, such that
for all $i, j \in \{1, \ldots, t\}$,
\begin{equation}\label{almost same cardinality of the p-taus, first application}
|p_\tau(\bar{a}, \bar{b}_i, \mcA)| - C \ \leq \ |p_\tau(\bar{a}, \bar{b}_j, \mcA)| \ \leq \ 
|p_\tau(\bar{a}, \bar{b}_i, \mcA)| + C.
\end{equation}

For numbers $c$ and $d$ we will use the notation `$c \pm d$' to denote some unspecified number in the interval 
$[c - d, c + d]$. In fact, in each case when the notation is used in this proof 
it can be taken to denote the average of a sequence of numbers
that are within distance $d$ from $c$.
Using~(\ref{almost same cardinality of the p-taus, first application}) we get
\begin{equation}\label{approximation of size of p-tau, first time}
|p_\tau(\bar{a}, \mcA)| = \sum_{i=1}^t |p_\tau(\bar{a}, \bar{b}_i, \mcA)| = 
t \cdot \big( |p_\tau(\bar{a}, \bar{b}_1, \mcA)| \pm C \big).
\end{equation}
Since $\mcA \in  \mbX_n^\varepsilon$ we also get
\begin{align}\label{approximation of size of p, first time}
&|p(\bar{a}, \mcA) \cap p_\tau(\bar{a}, \mcA)| = 
\sum_{i=1}^s |p(\bar{a}, \bar{b}_i, \mcA) \cap p_\tau(\bar{a}, \bar{b}_i,\mcA)| = \\
&\sum_{i=1}^s (\alpha \pm \varepsilon)|p_\tau(\bar{a}, \bar{b}_i, \mcA)| = 
s \cdot (\alpha \pm \varepsilon)\big(|p_\tau(\bar{a}, \bar{b}_1, \mcA)| \pm C\big).
\nonumber
\end{align}
It follows from~(\ref{approximation of size of p-tau, first time}) and~(\ref{approximation of size of p, first time}) that
if $|p_\tau(\bar{a}, \mcB_n)| > 0$, then
\[
\frac{|p(\bar{a}, \mcA)  \cap p_\tau(\bar{a}, \mcB_n)|}{|p_\tau(\bar{a}, \mcB_n)|} = 
\frac{s \cdot (\alpha \pm \varepsilon)}{t}.
\]
Since $\varepsilon > 0$ can be chosen arbitrarily small it follows that
$(p, p_\tau, q^*)$ is $s\alpha/t$-balanced with respect to $(\mbB, \mbbG)$.
\hfill $\square$

\medskip
\noindent
Next, we generalize Lemma~\ref{balanced triples, uniformly unbounded, dimension 1}
to the case when $p_\tau(\bar{x}, \bar{y})$ may have arbitrary (positive) $\bar{y}$-dimension.

\begin{lem}\label{balanced triples, uniformly unbounded}
Suppose that the following hold:
\begin{enumerate}
\item $p(\bar{x}, \bar{y})$ is a complete $(\sigma, \lambda)$-closure type,
\item $p_\tau(\bar{x}, \bar{y})$ is a complete $(\tau, \lambda + \gamma)$-closure type where 
$\gamma \geq \max(\lambda, \kappa + \kappa')$, and
\item $p_\tau$ is uniformly $\bar{y}$-unbounded.
\end{enumerate}
There is $\xi \in \mbbN^+$ such that if 
$q^*(\bar{x})$ is a complete $(\sigma, \xi)$-closure type, 
then $(p, p_\tau, q^*)$ is balanced with respect to $(\mbB, \mbbG)$.
Moreover, if $p_\tau$ strongly $\bar{y}$-unbounded then we can let $\xi = \lambda$.
\end{lem}

\noindent
{\bf Proof.}
We use induction on $\dim_{\bar{y}}(p_\tau)$. 
If $\dim_{\bar{y}}(p_\tau) = 1$, then the conclusion is given by
Lemma~\ref{balanced triples, uniformly unbounded, dimension 1}.
So suppose that $\dim_{\bar{y}}(p_\tau) = k+1$ where $k \geq 1$.

It follows from 
Definition~\ref{definition of dimension, part 2} and
Lemma~\ref{getting a strongly unbounded type of dimension 1} 
(and Remark~\ref{remark on restriction sim}),
that we can assume that
\begin{align*}
&\bar{y} = \bar{u}\bar{v}, \\
&\dim_{\bar{u}}(p_\tau \uhrc \bar{x}\bar{u}) = k, \\
&\dim_{\bar{v}}(p_\tau) = 1,   \text{ and } \\
&\text{$p_\tau$ is strongly $\bar{v}$-unbounded.}
\end{align*}
Let $r(\bar{x}, \bar{u}) = p \uhrc \bar{x}\bar{u}$ and  $r_\tau(\bar{x}, \bar{u}) = p_\tau \uhrc \bar{x}\bar{u}$.
so $\dim_{\bar{u}}(r_\tau) = k  \geq 1$.
By Lemma~\ref{dimension and unboundedness},
 $r_\tau$ is not $\bar{u}$-bounded and therefore,
by Lemma~\ref{not bounded implies uniformly unbounded},
it is uniformly $\bar{u}$-unbounded.
By the induction hypothesis,
there are $\xi \in \mbbN$ such that if $q^*(\bar{x})$ is a complete $(\sigma, \xi)$-closure type,
then $(r, r_\tau, q^*)$ is balanced with respect to $(\mbB, \mbbG)$.
Note that if $p_\tau$ is strongly $\bar{y}$-unbounded, 
then $r_\tau$ is strongly $\bar{u}$-unbounded, so the induction hypothesis 
says that we can let $\xi = \lambda$.

Let $q^*(\bar{x})$ be a complete $(\sigma, \xi)$-closure type, so there is $\alpha$ such that 
$(r, r_\tau, q^*)$ is $\alpha$-balanced with respect to $(\mbB, \mbbG)$.
We want to show that $(p, p_\tau, q^*)$ is balanced.
As $(r, r_\tau, q^*)$ is $\alpha$-balanced it follows  that, for all $\varepsilon > 0$, if 
\[
\mbX_n^\varepsilon = \big\{ \mcA \in \mbW_n : (r, r_\tau, q) \text{ is $(\alpha, \varepsilon)$-balanced in } \mcA \big\}
\]
then $\lim_{n\to\infty} \mbbP_n\big(\mbX_n^\varepsilon\big) = 1$.
By Lemma~\ref{balanced triples, strongly unbounded, dimension 1},
there is $\beta$ such that
$(p, p_\tau, r)$ is $\beta$-balanced with respect to $(\mbB, \mbbG)$.
So for all $\varepsilon > 0$, if
\[
\mbY_n^\varepsilon = \big\{ \mcA \in \mbW_n : (p, p_\tau, r) \text{ is $(\beta, \varepsilon)$-balanced in } \mcA \big\}
\]
then $\lim_{n\to\infty} \mbbP_n\big(\mbY_n^\varepsilon\big) = 1$.
We also get  $\lim_{n\to\infty} \mbbP_n\big(\mbX_n^\varepsilon \cap \mbY_n^\varepsilon\big) = 1$.

It now suffices to prove that for every $\delta > 0$ there is $\varepsilon > 0$ such that 
$(p, p_\tau, q^*)$ is $(\alpha\beta, \delta)$-balanced in every
$\mcA \in  \mbX_n^\varepsilon \cap \mbY_n^\varepsilon$.
Suppose that $\mcA \in \mbX_n^\varepsilon \cap \mbY_n^\varepsilon$ and $\mcA \models q^*(\bar{a})$.
It is enough to show that for any $\delta > 0$, if $\varepsilon > 0$ is sufficiently small 
(where $\varepsilon$ depends only on $\delta$) and $|p_\tau(\bar{a}, \mcA)| > 0$, then
\begin{equation}\label{p divided by p-tau}
\alpha\beta - \delta \ \leq \ \frac{|p(\bar{a}, \mcA) \cap p_\tau(\bar{a}, \mcA)|}{|p_\tau(\bar{a}, \mcA)|} 
\ \leq \ \alpha\beta + \delta.
\end{equation}

By Lemma~\ref{all tau-neighbourhood types have roughly the same cardinality},
there is a constant $C > 0$, depending only on $\mbB$ and $|\bar{x}\bar{y}|$, 
such that if $\bar{b}, \bar{b}' \in r_\tau(\bar{a}, \mcA)$,
then 
\begin{equation}\label{almost same cardinality of the r-taus}
|p_\tau(\bar{a}, \bar{b}, \mcA)| - C \ \leq \ |p_\tau(\bar{a}, \bar{b}', \mcA)| \ \leq \ 
|p_\tau(\bar{a}, \bar{b}, \mcA| + C.
\end{equation}
Below we use the notation `$c \pm d$' in the way explained in the proof of
Lemma~\ref{balanced triples, uniformly unbounded, dimension 1}.
Using~(\ref{almost same cardinality of the r-taus}) and that 
$\mcA \in  \mbX_n^\varepsilon \cap \mbY_n^\varepsilon$ we get, for an arbitrary 
choice of $\bar{b}_0 \in r_\tau(\bar{a}, \mcA)$,
\begin{equation}\label{approximation of size of p-tau}
|p_\tau(\bar{a}, \mcA)| = \sum_{\bar{b} \in r_\tau(\bar{a}, \mcA)} |p_\tau(\bar{a}, \bar{b}, \mcA)| = 
|r_\tau(\bar{a}, \mcA)| \cdot \big( |p_\tau(\bar{a}, \bar{b}_0, \mcA)| \pm C \big)
\end{equation}
and 
\begin{align}\label{approximation of size of p}
&|p(\bar{a}, \mcA) \cap p_\tau(\bar{a}, \mcA)| = \\
&\sum_{\bar{b} \in r(\bar{a}, \mcA)} |p(\bar{a}, \bar{b}, \mcA) \cap p_\tau(\bar{a}, \bar{b}, \mcA)| = 
\sum_{\bar{b} \in r(\bar{a}, \mcA)} (\beta \pm \varepsilon)|p_\tau(\bar{a}, \bar{b}, \mcA)| = \nonumber \\
&|r(\bar{a}, \mcA)|(\beta \pm \varepsilon)\big(|p_\tau(\bar{a}, \bar{b}_0, \mcA)| \pm C\big) = \nonumber \\
&(\alpha \pm \varepsilon)|r_\tau(\bar{a}, \mcA)| (\beta \pm \varepsilon)\big(|p_\tau(\bar{a}, \bar{b}_0, \mcA)| \pm C\big).
\nonumber
\end{align}
It follows from~(\ref{approximation of size of p-tau}) and~(\ref{approximation of size of p}) that
\[
\frac{|p(\bar{a}, \mcA) \cap p_\tau(\bar{a}, \mcA)|}{|p_\tau(\bar{a}, \mcA)|} = (\alpha \pm \varepsilon)(\beta \pm \varepsilon)
\]
so if $\varepsilon > 0$ is small enough we get~(\ref{p divided by p-tau}) and the proof is finished.
\hfill $\square$

\medskip
\noindent
Finally, we generalize 
Lemma~\ref{balanced triples, uniformly unbounded}
to the case when $p(\bar{x}, \bar{y})$ is a, {\em not} necessarily complete, $(\sigma, \lambda)$-closure type.

\begin{prop}\label{incomplete balanced triples, uniformly unbounded}
Suppose that the following hold:
\begin{enumerate}
\item $p(\bar{x}, \bar{y})$ is a, not necessarily complete, $(\sigma, \lambda)$-closure type,
\item $p_\tau(\bar{x}, \bar{y})$ is a complete $(\tau, \lambda + \gamma)$-closure type where 
$\gamma \geq \max(\lambda, \kappa + \kappa')$, and
\item $p_\tau$ is uniformly $\bar{y}$-unbounded.
\end{enumerate}
There is $\xi \in \mbbN^+$ (depending only on $p_\tau$) such that if 
$q(\bar{x})$ is a complete $(\sigma, \xi)$-closure type, 
then $(p, p_\tau, q)$ is balanced with respect to $(\mbB, \mbbG)$.
Moreover, if $p_\tau$ strongly $\bar{y}$-unbounded then we can let $\xi = \lambda$.
\end{prop}

\noindent
{\bf Proof.}
In this proof we omit saying ``with respect to $(\mbB, \mbbG)$'' since this will always be the case.
Let $p$ and $p_\tau$ be as assumed in the lemma.
Then there are $s \in \mbbN^+$ and
non-equivalent complete $(\sigma, \lambda)$-closure types $p_i(\bar{x}, \bar{y})$, for $i = 1, \ldots, s$,
such that $p$ is equivalent to $\bigvee_{i=1}^s p_i$.
By Lemma~\ref{balanced triples, uniformly unbounded}
there are $\xi_i$, $i = 1, \ldots, s$ such that, for each $i$, if $q(\bar{x})$ is a complete $(\sigma, \xi_i)$-closure type,
then $(p_i, p_\tau, q)$ is balanced, say $\alpha_{q, i}$-balanced.
Let $\xi = \max\{\xi_1, \ldots, \xi_s\}$.
By Lemma~\ref{on stronger conditioning in the balance},
for every complete $(\sigma, \xi)$-closure type $q(\bar{x})$,
$(p_i, p_\tau, q)$ is $\alpha_{q, i}$-balanced for all $i$.
It now follows straightforwardly from the definition of balanced triples that
$(p, p_\tau, q)$ is $(\alpha_{q, 1} + \ldots + \alpha_{q, s})$-balanced for every  
complete $(\sigma, \xi)$-closure type $q(\bar{x})$.
The ``moreover'' part follows directly from the given argument and
Lemma~\ref{balanced triples, uniformly unbounded}.
\hfill $\square$

\section{Asymptotic elimination of aggregation functions}\label{Asymptotic elimination of aggregation functions}

\noindent
In this section we complete the induction step for part~(1) of 
Assumption~\ref{induction hypothesis}
and we prove the main results of this investigation, which concern asymptotic elimination of aggregation functions and
the asymptotic distribution of values of formulas.
We prove two versions of these results.

In the first version, treated in Section~\ref{Results without further assumptions},
we assume that $\kappa = \kappa' = 0$ 
(where these numbers come from Assumption~\ref{induction hypothesis}), which follows if
only strongly unbounded closure types are used as conditioning formulas in aggregations
(as will be explained).
Under this assumption the inductive step of part~(1) of Assumption~\ref{induction hypothesis}
is completed by 
Proposition~\ref{asymptotic elimination for strongly unbounded aggregations}
and Remark~\ref{remark on kappa without further assumptions}.
Then Proposition~\ref{asymptotic elimination for strongly unbounded aggregations} 
is used to prove our first main result about
asymptotic elimination of aggregation functions and
about the asymptotic distribution of values of formulas, 
Theorem~\ref{main result about strongly unbounded aggregations}
and its corollary. Roughly speaking, the assumptions in 
Theorem~\ref{main result about strongly unbounded aggregations} are, 
besides Assumption~\ref{properties of the base structures},
that all formulas associated to the $PLA^*(\sigma)$-network, and the formula that expresses the query,
uses only continuous aggregation functions and use only strongly unbounded closure types as 
conditioning formulas in aggregations.

In Section~\ref{Results with an additional assumption}
we prove another version of the results in Section~\ref{Results without further assumptions}
under an additional (but reasonable I think) assumption on the sequence of base structures $\mbB = (\mcB_n : n \in \mbbN^+)$.
The additional assumption is stated as Assumption~\ref{relative frequency of tau-closure types}
and it is satisfied by all examples of sequences of base structures described in 
Section~\ref{Examples of sequences of base structures}
(see Examples~\ref{easy examples satisfy a stronger property} and~\ref{Galton-Watson trees satisfy a stronger property} below).
With stronger restrictions on the sequence of base structures we can now prove results about
asymptotic elimination of aggregation functions and about the asymptotic distribution of values with 
{\em less} restrictive assumptions on the formulas to which these results apply.
Theorem~\ref{asymptotic elimination under stronger conditions} and 
Remark~\ref{remark on kappa with an additional assumption}
complete the inductive step of part~(1) of Assumption~\ref{induction hypothesis}
without any additional assumption on the numbers $\kappa$ and $\kappa'$ and with a more liberal
(than in Proposition~\ref{asymptotic elimination for strongly unbounded aggregations})
assumption on conditioning subformulas that appear in aggregations.
Then Theorem~\ref{asymptotic elimination under stronger conditions}
is used to prove our second version of the main result,
Theorem~\ref{main result about aggregations in general with extra assumption} and its corollary, 
where it is assumed that
the sequence of base structures $\mbB$ satisfies both 
Assumption~\ref{properties of the base structures}
and~\ref{relative frequency of tau-closure types}.
In Theorem~\ref{main result about aggregations in general with extra assumption} we can relax
(compared to Theorem~\ref{main result about strongly unbounded aggregations})
the conditions on the formulas used by $PLA^*(\sigma)$-networks and to express queries in such a way that we only assume
that all aggregation functions that are used are continuous and that all conditioning subformulas that are used are
either bounded or ``positive'' $(\sigma, \lambda)$-closure types, for arbitrary $\lambda$ 
(where the notion of ``positive'' closure type is defined below).

We adopt the same assumptions as in 
Sections~\ref{Proving convergence in the inductive step}
and~\ref{Finding the balance in the inductive step}, 
so in particular, 
{\bf \em $\tau \subseteq \sigma$ are finite relational signatures,
$\mbB = (\mcB_n : n \in \mbbN^+)$ is a sequence of (finite) $\tau$-structures
which satisfies Assumption~\ref{properties of the base structures},
$\mbW_n$ is the set of all $\sigma$-structures which expand $\mcB_n$, $\mbbG$ is a $PLA^*(\sigma)$-network
based on $\tau$ and $\mbbP_n$ is the probability distribution on $\mbW_n$ which is induced by $\mbbG$.
Recall that if $R \in \sigma \setminus \tau$ then $\theta_R$ denotes the $PLA^*(\mr{par}(R) \cup \tau)$-formula 
which is associated to $R$ by $\mbbG$ as in Definition~\ref{definition of PLA-network},
where $\mr{par}(R)$ is the set of parents of $R$ in the DAG of $\mbbG$.}

\subsection{Results without further assumptions}\label{Results without further assumptions}

We now generalize the notion of complete $(\tau, \lambda)$-closure type $p(\bar{x}, \bar{y})$ to the notion
of (not necessarily complete) $\bar{y}$-positive $(\sigma, \lambda)$-closure type. 
The main results will consider aggregations conditioned on $\bar{y}$-positive closure types (for some sequence $\bar{y}$ 
of variables).

\begin{defin}\label{definition of y-positive} {\rm
Let $p(\bar{x}, \bar{y})$ be a $(\sigma, \lambda)$-closure type for some $\lambda \in \mbbN$.
We say that $p$ is {\bf \em $\bar{y}$-positive} (with respect to $(\mbB, \mbbG)$) if 
$p \uhrc \tau$ is cofinally satisfiable and
there is $\gamma \in \mbbN$ such that $(p, p \uhrc \tau, q)$ is positively balanced
(with respect to $(\mbB, \mbbG)$) if $q(\bar{x})$ is a complete $(\sigma, \gamma)$-closure type
such that $p \wedge q$ is cofinally satisfiable.
}\end{defin}

\noindent
Note that for every $\lambda \in \mbbN$, every complete $(\tau, \lambda)$-closure type is also
a $(\sigma, \lambda)$-closure type. The next lemma is a direct consequence of
Definition~\ref{definition of y-positive}.

\begin{lem}\label{tau-types are positive}
For all $\lambda \in \mbbN$, if $p(\bar{x}, \bar{y})$ is a complete $(\tau, \lambda)$-closure type
which is not $\bar{y}$-bounded (hence uniformly $\bar{y}$-unbounded),
then $p$ is $\bar{y}$-positive.
\end{lem}

\noindent
Consequently any statement that holds for all $\bar{y}$-positive $(\sigma, \lambda)$-closure types also holds 
for all $\bar{y}$-unbounded complete $(\tau, \lambda)$-closure types.

Recall the numbers $\kappa$ and $\kappa'$ from Assumption~\ref{induction hypothesis} which will be assumed to
be zero in the next couple of results.
We begin with a lemma essentially saying that we can always condition on a $\bar{y}$-positive 
uniformly $\bar{y}$-unbounded $(\sigma, 0)$-closure type.

\begin{lem}\label{(p-2, p-1, q) is balanced when kappa = kappa' = 0}
Suppose that $\kappa = \kappa' = 0$.
Let $p_1(\bar{x}, \bar{y})$ and $p_2(\bar{x}, \bar{y})$ be $(\sigma, 0)$-closure types and $q^*(\bar{x})$
a complete $(\sigma, 0)$-closure type.
If $p_2$ is uniformly $\bar{y}$-unbounded and $\bar{y}$-positive then $(p_1, p_2, q^*)$ is balanced with respect to $(\mbB,  \mbbG)$.
\end{lem}

\noindent
{\bf Proof.}
Let $p_1(\bar{x}, \bar{y})$, $p_2(\bar{x}, \bar{y})$ and $q^*(\bar{x})$ be as assumed.
Also assume that $\kappa = \kappa' = 0$.
Since $p_2$ is $\bar{y}$-positive there is $\gamma \in \mbbN$ such that if $q(\bar{x})$ is a 
complete $(\sigma, \gamma)$-closure type such that $p_2 \wedge q$ is cofinally satisfiable,
then $(p_2, p_2 \uhrc \tau, q)$ is $\beta_q$-balanced for some $\beta_q > 0$.
If $p_1 \wedge p_2 \wedge q^*$ is not cofinally satisfiable then 
(by Lemma~\ref{balance follows from inconsistency})
$(p_1, p_2, q^*)$ is 0-balanced.

Now suppose that $p_1  \wedge p_2 \wedge q^*$ is cofinally satisfiable.
Let $q(\bar{x})$ be a complete $(\sigma, \gamma)$-closure type such that $p_2 \wedge q$ is cofinally satisfiable
and $q \uhrc 0 = q^*$. Then $(p_2, p_2 \uhrc \tau, q)$ is $\beta_q$-balanced for some $\beta_q > 0$.
By Proposition~\ref{incomplete balanced triples, uniformly unbounded},
including its ``moreover'' part, $(p_2, p_2 \uhrc \tau, q^*)$ is $\beta$-balanced for some $\beta$.
Since $q \models q^*$ it follows from
Lemma~\ref{on stronger conditioning in the balance}
that $(p_2, p_2 \uhrc \tau, q)$ is $\beta$-balanced, so $\beta = \beta_q > 0$.
By Proposition~\ref{incomplete balanced triples, uniformly unbounded},
including its ``moreover'' part, $(p_1 \wedge p_2, p_2 \uhrc \tau, q^*)$ is $\alpha$-balanced for some $\alpha$.
It now follows from the definition of balanced triples that $(p_1, p_2, q^*)$ is $\alpha/\beta$-balanced.
\hfill $\square$

\bigskip

\noindent
Recall Definition~\ref{Definition of basic formula} of $(\sigma, \lambda)$-basic formula.

\begin{prop}\label{asymptotic elimination for strongly unbounded aggregations}
Suppose that $\kappa = \kappa' = 0$ where $\kappa$ and $\kappa'$ are as in 
Assumption~\ref{induction hypothesis}.
Let $\varphi(\bar{x}) \in PLA^*(\sigma)$ and suppose that if
\[
F\big(\varphi_1(\bar{y}, \bar{z}), \ldots, \varphi_m(\bar{y}, \bar{z}) : \bar{z} : 
\chi_1(\bar{y}, \bar{z}), \ldots, \chi_m(\bar{y}, \bar{z})\big)
\]
is a subformula of $\varphi(\bar{x})$ then,
for all $i = 1, \ldots, m$, $\chi_i(\bar{y}, \bar{z})$ is a strongly $\bar{z}$-unbounded and
$\bar{z}$-positive $(\sigma, 0)$-closure type
and $F$ is continuous.
Then $\varphi(\bar{x})$ is asymptotically equivalent to a $(\sigma, 0)$-basic formula.

If, in addition, conditions~(I) and~(II) of
Remark~\ref{remark about basic probability formulas}
hold, then we can replace `continuous' by (the weaker notion) `admissible' and the conclusion
that $\varphi(\bar{x})$ is asymptotically equivalent to a $(\sigma, 0)$-basic formula still follows.
\end{prop}

\noindent
{\bf Proof.}
We use induction on the complexity of formulas.
If $\varphi(\bar{x})$ is aggregation-free, then the conlusion follows from 
Lemma~\ref{connectives and basic formulas}~(ii).
Suppose that $\varphi(\bar{x})$ has the form
$\msfC(\psi_1(\bar{x}), \ldots, \psi_k(\bar{x}))$
where $\msfC : [0, 1]^k \to [0, 1]$ is a continuous
and each $\psi_i(\bar{x})$ is asymptotically equivalent
to a $(\sigma, 0)$-basic formula $\psi'_i(\bar{x})$.
Since $\msfC$ is continuous it follows that $\varphi(\bar{x})$ and
$\msfC(\psi'_1(\bar{x}), \ldots, \psi'_k(\bar{x}))$ are asymptotically equivalent.
Lemma~\ref{connectives and basic formulas}~(i) implies that
$\msfC(\psi'_1(\bar{x}), \ldots, \psi'_k(\bar{x}))$ is equivalent, and hence
asymptotically equivalent, to a $(\sigma, 0)$-basic formula $\varphi'(\bar{x})$.
By transitivity of asymptotic equivalence it follows that $\varphi$ and $\varphi'$ are asymptotically equivalent.

Now suppose that $\varphi(\bar{x})$ has the form
\[
F\big(\varphi_1(\bar{x}, \bar{y}), \ldots, \varphi_m(\bar{x}, \bar{y}) : \bar{y} : 
\chi_1(\bar{x}, \bar{y}), \ldots, \chi_m(\bar{x}, \bar{y})\big)
\]
where $F : \big([0, 1]^{<\omega}\big)^m \to [0, 1]$ is continuous, 
each $\varphi_i(\bar{x}, \bar{y})$ is asymptotically equivalent to a $(\sigma, 0)$-basic formula $\psi_i(\bar{x}, \bar{y})$,
and each $\chi_i(\bar{x}, \bar{y})$ is a strongly $\bar{y}$-unbounded and $\bar{y}$-positive $(\sigma, 0)$-closure type.

Let $L_0$ be the set of all complete $(\sigma, 0)$-closure types 
and let $L_1$ be the set of all (not necessarily complete) $(\sigma, 0)$-closure types.
For every $\varphi(\bar{x}, \bar{y}) \in L_0$ let $L_{\varphi(\bar{x}, \bar{y})}$ be the set of all
$(\sigma, 0)$-closure types in the variables $\bar{x}\bar{y}$ that are strongly $\bar{y}$-unbounded and $\bar{y}$-positive.
Due to Theorem~\ref{general asymptotic elimination}
it now suffices to show that 
Assumption~\ref{assumptions on the basic logic}
is satisfied.
Part~(1)
of Assumption~\ref{assumptions on the basic logic} 
follows from
Lemma~\ref{connectives and basic formulas}~(ii),
so we verify part~(2) of the same assumption.

Let $p_1(\bar{x}, \bar{y}), \ldots, p_k(\bar{x}, \bar{y}) \in L_0$ and 
$\chi_i(\bar{x}, \bar{y}) \in L_{p_i(\bar{x}, \bar{y})}$ for $i = 1, \ldots, k$.
Let $q_1(\bar{x}), \ldots, q_s(\bar{x})$ enumerate, up to equivalence, all complete $(\sigma, 0)$-closure types
in the variables $\bar{x}$.
By 
Lemma~\ref{(p-2, p-1, q) is balanced when kappa = kappa' = 0},
$(p_j, \chi_j, q_i)$ is balanced for all $i = 1, \ldots, s$ and $j = 1, \ldots, k$.
This means that there are $\alpha_{i, j} \in [0, 1]$, for $i = 1, \ldots, s$ and $j = 1, \ldots, k$ such that for every
$\varepsilon > 0$ there are $\mbY_n^\varepsilon \subseteq \mbW_n$ such that
$\lim_{n\to\infty} \mbbP_n(\mbY_n^\varepsilon) = 1$ and for all $i, j$, all $\mcA \in \mbY_n^\varepsilon$,
and all $\bar{a} \in (B_n)^{|\bar{x}|}$, if $\mcA \models q_i(\bar{a})$ then
\[
(\alpha_{i, j} - \varepsilon)|\chi_j(\bar{a}, \mcA)| \ \leq \ 
|p_j(\bar{a}, \mcA) \cap \chi_j(\bar{a}, \mcA)| \ \leq \ (\alpha_{i, j} + \varepsilon)|\chi_j(\bar{a}, \mcA)|.
\]
Thus condition~(d) of part~(2) of Assumption~\ref{assumptions on the basic logic} is satisfied.

Moreover, for all $n \in \mbbN^+$ and $\mcA \in \mbW_n$,
$\mcA \models \forall \bar{x} \bigvee_{i=1}^s q_i(\bar{x})$ and
$\mcA \models \forall \bar{x} \neg (q_i(\bar{x}) \wedge q_j(\bar{x}))$ if $i \neq j$,
so conditions~(a) and~(b) of part~(2) of Assumption~\ref{assumptions on the basic logic} are also satisfied.

Let $\chi'_1(\bar{x}), \ldots, \chi'_t(\bar{x})$ enumerate, up to equivalence, 
all complete $(\sigma, 0)$-closure types $\chi'(\bar{x})$
such that, for some $i \in \{1, \ldots, t\}$, $\chi' \uhrc \tau$ is not equivalent to $\chi_i \uhrc \bar{x}$.
Since each $\chi_i(\bar{x}, \bar{y})$ is strongly $\bar{y}$-unbounded it follows
that $\chi^*_i = \chi_i \uhrc \tau$ is strongly $\bar{y}$-unbounded and it follows from
Lemma~\ref{for strongly y-unbounded the x-part determines if y exists}
that for all sufficiently large $n$ and all $\mcA \in \mbW_n$
\[
\mcA \models \forall \bar{x} \Big(\Big(\bigvee_{i=1}^m \neg \exists \bar{y} \chi^*_i(\bar{x}, \bar{y}) \Big) \leftrightarrow
\Big(\bigvee_{i=1}^t \chi'_i(\bar{x})\Big)\Big).
\]
Since we also assume that each $\chi_i$ is $\bar{y}$-positive it follows that 
if $\mbZ_n$ is the set of $\mcA \in \mbW_n$ such that 
\[
\mcA \models \forall \bar{x} \Big(\Big(\bigvee_{i=1}^m \neg \exists \bar{y} \chi_i(\bar{x}, \bar{y}) \Big) \leftrightarrow
\Big(\bigvee_{i=1}^t \chi'_i(\bar{x})\Big)\Big),
\]
then $\lim_{n\to\infty} \mbbP_n(\mbZ_n) = 1$ (and note the occurence of $\chi_i$ in the last formula
instead of $\chi^*_i$ in the one before).
It follows that also condition~(c) of part~(2) of Assumption~\ref{assumptions on the basic logic} is satisfied.

Suppose that, in addition, conditions~(I) and~(II) of
Remark~\ref{remark about basic probability formulas} hold and recall that we already before assumed that $\kappa = \kappa' = 0$.
From that remark and
Remark~\ref{remark on eventually constant pairs in the general case}
it follows that if
$\lambda, \gamma \in \mbbN$,
$p(\bar{x})$ is a complete $(\sigma, \lambda)$-closure type, and 
$p_\tau(\bar{x})$ is a complete $(\tau, \lambda + \gamma)$-closure type,
then $(p, p_\tau)$ is eventually constant with respect to $(\mbB, \mbbG)$.
As mentioned in 
Remark~\ref{remark on eventual constancy and balance},
the additional condition in part~(ii) of 
Theorem~\ref{general asymptotic elimination}
is now satisfied, so $F$ can be asymptotically eliminated
(by the same argument as above)  if it is admissible.
\hfill $\square$

\begin{rem}\label{remark on kappa without further assumptions} {\rm
Suppose that $\kappa = \kappa' = 0$,
$\sigma^+$ is a nonempty finite relational signature, $\sigma \subset \sigma^+$, and that
$\mbbG^+$ is a $PLA^*(\sigma^+)$-network such that $\mr{mp}(\mbbG^+) = \mr{mp}(\mbbG) + 1$.
Furthermore, suppose that $\sigma = \{R \in \sigma^+ : \mr{mp}(R) < \mr{mp}(\mbbG^+)\}$
and that $\mbbG$ is the subnetwork of $\mbbG^+$ which is induced by $\sigma$.
Then, for each $R \in \sigma^+ \setminus \sigma$, the corresponding formula $\theta_R$ (of $\mbbG^+$)
belongs to $PLA^*(\sigma)$.
Suppose that, for every $R \in \sigma^+ \setminus \sigma$
and every subformula of $\theta_R$ of the form
\[
F\big(\varphi_1(\bar{y}, \bar{z}), \ldots, \varphi_m(\bar{y}, \bar{z}) : \bar{z} : 
\chi_1(\bar{y}, \bar{z}), \ldots, \chi_m(\bar{y}, \bar{z})\big)
\]
it holds that,
for all $i = 1, \ldots, m$, $\chi_i(\bar{y}, \bar{z})$ is a strongly $\bar{z}$-unbounded and
$\bar{z}$-positive $(\sigma, 0)$-closure type
and $F$ is continuous.

By Proposition~\ref{asymptotic elimination for strongly unbounded aggregations}, for every 
$R \in \sigma^+ \setminus \sigma$,
$\theta_R(\bar{x})$ is asymptotically equivalent to a $(\sigma, 0)$-basic formula.
Thus part~(1) of Assumption~\ref{induction hypothesis} holds if we replace 
$\sigma, \sigma', \kappa, \mbbG$, and $\mbbG'$ by
$\sigma^+, \sigma, 0, \mbbG^+$, and $\mbbG$, respectively.
This concludes the proof of the inductive step for part~(1)  of Assumption~\ref{induction hypothesis}
(provided that the assumptions of this remark are satisfied).
Note that in Remarks~\ref{the new kappa'} 
and~\ref{the new kappa' for balance},
$\kappa_1$ was defined as $\kappa + \kappa'$, so 
under the assumptions of this remark we get $\kappa_1 = 0$.
}\end{rem}

\begin{theor}\label{main result about strongly unbounded aggregations}
Suppose that Assumption~\ref{properties of the base structures} holds.
Suppose that for every $R \in \sigma \setminus \tau$
and every subformula of $\theta_R$ of the form
\begin{equation}\label{a subformula in the main results without additional assumptions}
F\big(\varphi_1(\bar{y}, \bar{z}), \ldots, \varphi_m(\bar{y}, \bar{z}) : \bar{z} : 
\chi_1(\bar{y}, \bar{z}), \ldots, \chi_m(\bar{y}, \bar{z})\big)
\end{equation}
it holds that,
for all $i = 1, \ldots, m$, $\chi_i(\bar{y}, \bar{z})$ is a strongly $\bar{z}$-unbounded and
$\bar{z}$-positive $(\sigma, 0)$-closure type
and $F$ is continuous.

Let $\varphi(\bar{x}) \in PLA^*(\sigma)$ and suppose that for every subformula of $\varphi(\bar{x})$ of the
form~(\ref{a subformula in the main results without additional assumptions}) the same conditions as stated above hold.
Then:\\
(i) $\varphi(\bar{x})$ is asymptotically equivalent to a $(\sigma, 0)$-basic formula.\\
(ii) For every complete $(\tau, 0)$-closure type $p(\bar{x})$ there are $k \in \mbbN^+$, 
$c_1, \ldots, c_k \in [0, 1]$, and $\beta_1, \ldots, \beta_k \in [0, 1]$ such that 
for every $\varepsilon > 0$ there is $n_0$ such that if $n \geq n_0$, 
and $\mcB_n \models p(\bar{a})$ then 
\begin{align*}
&\mbbP_n \Big( \big\{ \mcA \in \mbW_n : \mcA(\varphi(\bar{a})) \in 
\bigcup_{i = 1}^k [c_i - \varepsilon, c_i + \varepsilon] \big\}\Big)
\geq 1 - \varepsilon  \text{ and, for all $i = 1, \ldots, k$,} \\
&\mbbP_n\big(\{\mcA \in \mbW_n : \mcA(\varphi(\bar{a})) \in [c_i - \varepsilon, c_i + \varepsilon] \}\big) \in 
[\beta_i - \varepsilon, \beta_i + \varepsilon].
\end{align*}
\end{theor}

\noindent
{\bf Proof.}
(i) We use induction on $\mr{mp}(\mbbG)$ where the base case is when $\mr{mp}(\mbbG) = -1$, or equivalently,
when $\sigma = \tau$ and the DAG of $\mbbG$ is empty.
As noted in Remark~\ref{remark on induction hypothesis and base case},
Assumption~\ref{induction hypothesis}
follows  from 
Lemma~\ref{base case}
if $\sigma = \tau$ and if $\kappa = \kappa' = 0$.
From Remarks~\ref{the new kappa'}, \ref{the new kappa' for balance} and \ref{remark on kappa without further assumptions}
it follows that if Assumption~\ref{induction hypothesis} holds for every $PLA^*(\sigma)$-network with
mp-rank $\rho$ with the choice $\kappa = \kappa' = 0$, then Assumption~\ref{induction hypothesis} also holds for every $PLA^*(\sigma)$-network
with mp-rank $\rho + 1$ and with the choice $\kappa = \kappa' = 0$. 
Hence, Assumption~\ref{induction hypothesis} holds, with the choice $\kappa = \kappa' = 0$, 
for every finite relational signature $\sigma \supseteq \tau$
and every $PLA^*(\sigma)$-network subject to the conditions of the theorem.
Proposition~\ref{asymptotic elimination for strongly unbounded aggregations}
now implies that 
if $\varphi(\bar{x}) \in PLA^*(\sigma)$ is as assumed in the present theorem,
then $\varphi(\bar{x})$ is asymptotically equivalent to a $(\sigma, 0)$-basic formula.

(ii) Let $\varphi(\bar{x})$ be asymptotically equivalent to the $(\sigma, 0)$-basic formula
$\bigwedge_{i = 1}^m (\varphi_i(\bar{x}) \to c_i)$, so each $\varphi_i(\bar{x})$ is a 
complete $(\sigma, 0)$-closure type. 
Without loss of generality we can assume that $\varphi_i(\bar{x})$, $i = 1, \ldots, m$, enumerate
all, up to equivalence, complete $(\sigma, 0)$-closure types that are cofinally satisfiable and
that $\varphi_i \leftrightarrow \varphi_j$ is not cofinally satisfiable if $i \neq j$.
Let $p(\bar{x})$ be a complete $(\tau, 0)$-closure type.
By Proposition~\ref{convergence of conditional probabilities, part 3}, for all $i$,
$(\varphi_i, p)$ converges to some $\alpha_i$.
So for every $\varepsilon > 0$, 
if $n$ is large enough, $\bar{a} \in (B_n)^{|\bar{x}|}$ and $\mcB_n \models p(\bar{a})$, then 
\[
\mbbP_n\big(\mbE_n^{\varphi_i(\bar{a})} \ | \ \mbE_n^{p(\bar{a})}\big) \in [\alpha_i - \varepsilon, \alpha_i + \varepsilon].
\]
and (as $\varphi(\bar{x})$ and $\bigwedge_{i = 1}^m (\varphi_i(\bar{x}) \to c_i)$ are asymptotically equivalent)
\[
\mbbP_n\Big(\big\{\mcA \in \mbW_n : 
|\mcA(\varphi(\bar{a})) - \mcA(\bigwedge_{i = 1}^m (\varphi_i(\bar{a}) \to c_i))| \leq \varepsilon \big\}\Big) \geq 1 - \varepsilon.
\]
Note that if $\mcA \models \varphi_i(\bar{a})$ then $\mcA(\bigwedge_{i = 1}^m (\varphi_i(\bar{a}) \to c_i)) = c_i$.
Since we assume that $\varphi_i(\bar{x})$, $i = 1, \ldots, m$, enumerate all, up to equivalence,
complete $(\sigma, 0)$-closure types that are cofinally satisfiable it follows that for all sufficiently large $n$
\[
\mbbP_n \Big( \big\{ \mcA \in \mbW_n : \mcA(\varphi(\bar{a})) \in 
\bigcup_{i = 1}^k [c_i - \varepsilon, c_i + \varepsilon] \big\}\Big)
\geq 1 - \varepsilon.
\]
Let $c \in \{c_1, \ldots, c_m\}$ and for simplicity of notation suppose that, for some $1 \leq s \leq m$,
$c = c_i$ if $i \leq s$ and $c \neq c_i$ if $i > s$.
Let $\beta_c = \alpha_1 + \ldots + \alpha_s$.
It now follows that if $n$ is large enough and $\mcB_n \models p(\bar{a})$ then 
\[
\mbbP_n\big(\{\mcA \in \mbW_n : \mcA(\varphi(\bar{a})) \in [c - \varepsilon, c + \varepsilon]\}\big)
\in  [\beta_c - s\varepsilon, \beta_c + s\varepsilon].
\]
The claim now follows since $\varepsilon > 0$ can be chosen as small as we like.
\hfill $\square$

\begin{cor}\label{corollary to main result about strongly unbounded aggregations}
Suppose that Assumption~\ref{properties of the base structures} holds.
Suppose that for every $R \in \sigma \setminus \tau$, $\theta_R$ is an {\rm  aggregation-free} formula.
Also suppose that $\varphi(\bar{x}) \in PLA^*(\sigma)$ 
and if
\[
F\big(\varphi_1(\bar{y}, \bar{z}), \ldots, \varphi_m(\bar{y}, \bar{z}) : \bar{z} : 
\chi_1(\bar{y}, \bar{z}), \ldots, \chi_m(\bar{y}, \bar{z})\big)
\]
is a subformula of $\varphi(\bar{x})$ then,
for all $i = 1, \ldots, m$, $\chi_i(\bar{y}, \bar{z})$ is a strongly $\bar{z}$-unbounded and
$\bar{z}$-positive $(\sigma, 0)$-closure type
and $F$ is {\rm admissible}.
Then the conclusions~(i) and~(ii) of 
Theorem~\ref{main result about strongly unbounded aggregations}
hold.
\end{cor}

\noindent
{\bf Proof.}
Part~(ii) follows from part~(i) in the same way as 
part~(ii) of 
Theorem~\ref{main result about strongly unbounded aggregations}
follows from part~(i) of that theorem.

So we consider part~(i). 
Suppose that, for every $R \in \sigma \setminus \tau$, $\theta_R$ is aggregation-free.
By 
Lemma~\ref{connectives and basic formulas}
each $\theta_R$ is equivalent 
to a complete $(\mr{par}(R) \cup \tau, 0)$-basic formula.

Since each $\theta_R$ satisfies the condition formulated in 
Theorem~\ref{main result about strongly unbounded aggregations}
we can argue as in the proof of that theorem and it follows that
Assumption~\ref{induction hypothesis} holds with the choice $\kappa = \kappa' = 0$ for
every finite relational signature $\sigma \supseteq \tau$ and every $PLA^*(\sigma)$-network based on $\tau$
subject to the conditions of the present corollary.
Since each $\theta_R$ satisifies condition~(I) of
Remark~\ref{remark about basic probability formulas} we can, via induction and
Remark~\ref{remark on eventually constant pairs in the general case},
conclude that, for every complete $(\tau, 0)$-closure type $p_\tau(\bar{x})$ and complete $(\sigma, 0)$-closure type
$p(\bar{x})$, $(p, p_\tau)$ is eventually constant.
As explained in Remark~\ref{remark on eventual constancy and balance} it follows that
the additional condition in part~(ii) of 
Theorem~\ref{general asymptotic elimination}
is now satisfied, so it follows (by the same argument as in the proof of
Proposition~\ref{asymptotic elimination for strongly unbounded aggregations})
that all admissible aggregation functions of $\varphi(\bar{x})$ can be asymptotically eliminated
provided that $\varphi(\bar{x})$ is as assumed in the present corollary.
\hfill $\square$

\begin{exam}\label{example for corollary to theorem with no extra assumption}{\rm
Suppose that $\sigma = \tau \cup \{R\}$ where $R \notin \tau$ and $R$ has arity $k$. Also suppose
that $\theta_R(\bar{x})$, where $\bar{x} = (x_1, \ldots, x_k)$, has the form $\bigwedge_{i = 1}^s (q_i(\bar{x}) \to c_i)$ where 
$q_1(\bar{x}), \ldots, q_s(\bar{x})$ enumerate all, up to equivalence, $(\sigma, 0)$-closure types in the variables $\bar{x}$.
Also suppose that $\varphi(\bar{y})$ is a first-order formula such that every quantification in $\varphi(\bar{y})$
has the form
\begin{enumerate}
\item[\textbullet] $\exists u \big(p(u) \wedge $ ``there is no tuple $\bar{z}$ containg $u$ and some member of $\bar{w}$ such that $E(\bar{z})$ for some $E \in \tau$'' $\wedge \ \psi(u, \bar{w})\big)$, or the form
\item[\textbullet] $\forall u \big((p(u) \wedge $
``there is no tuple $\bar{z}$ containg $u$ and some member of $\bar{w}$ such that $E(\bar{z})\big)$ for some $E \in \tau$''
$) \to
\psi(u, \bar{w})\big)$
\end{enumerate}
where $p(u)$ is a complete $(\tau, 0)$-closure type.

Since the existential and universal quantifiers can be expressed by the aggregation functions maximum and minimum,
which are admissible, it follows from 
Corollary~\ref{corollary to main result about strongly unbounded aggregations} that $\varphi(\bar{y})$
is asymptotically equivalent to a $(\sigma, 0)$-basic formula $\varphi'(\bar{y})$. Since $\varphi(\bar{y})$ is 0/1-valued it 
follows that $\varphi'(\bar{y})$ is equivalent to a quantifier-free first-order formula.
}\end{exam}

\subsection{Results with an additional assumption}\label{Results with an additional assumption}

Now we add Assumption~\ref{relative frequency of tau-closure types} below,
about the sequence $\mbB = (\mcB_n : n \in \mbbN^+)$ of base structures, to our previous assumptions
(which are Assumption~\ref{properties of the base structures} and Assumption~\ref{induction hypothesis}).
This will allow us to prove results about asymptotic elimination of aggregation functions
and the distribution of values of formulas with {\em less restrictive} conditions
than in Section~\ref{Results without further assumptions}
on formulas used
by $PLA^*(\sigma)$-networks and on formulas used to express queries.
All our examples in Section~\ref{Examples of sequences of base structures}
satisfy Assumption~\ref{relative frequency of tau-closure types} as shown below in the most difficult case.

\begin{assump}\label{relative frequency of tau-closure types} {\rm
Suppose that if $\lambda, \mu \in \mbbN$, $\lambda \leq \mu$, 
$r(\bar{x})$ is a complete $(\tau, \mu)$-neigh\-bour\-hood type, $p(\bar{x}) = r \uhrc \lambda$,
$p(\bar{x})$ and $r(\bar{x})$ are strongly unbounded, and $\dim(p) = \dim(r) = 1$, then
\[
\lim_{n\to\infty} \frac{|r(\mcB_n)|}{|p(\mcB_n)|} \quad \text{ exists.}
\]
}\end{assump}

\noindent
Observe that if $p$ and $r$ are as in 
Assumption~\ref{relative frequency of tau-closure types}, 
then $|r(\mcB_n)| \leq |p(\mcB_n)|$ so the assumption stipulates
that the limit is a real number in the interval $[0, 1]$.

\begin{exam}\label{easy examples satisfy a stronger property}{\rm
It can be proved in a straightforward way that 
$\mbB = (\mcB_n : n \in \mbbN^+)$ defined in Examples~\ref{example with empty tau}--~\ref{example of grids}
satisfy Assumption~\ref{relative frequency of tau-closure types}.
}\end{exam}

\begin{exam}\label{Galton-Watson trees satisfy a stronger property} {\bf (Galton-Watson trees)} {\rm
Here we continue the reasoning begun in
Example~\ref{example of trees}.
So we let $\tau = \{E, \sqsubset\}$ where both relation symbols are binary.
We consider ordered rooted trees (often just called trees) as $\tau$-structures as explained in
Example~\ref{example of trees}
and we adopt other definitions from that example.

Fix some $\delta \in \mbbN^+$ and let $\Delta = \delta + 1$.
Also let $\mbB = (\mcB_n : n \in \mbbN^+)$ where each $\mcB_n$ is a $\delta$-bounded tree.

Let $\lambda \in \mbbN$, $\bar{x} = (x_1, \ldots, x_k)$, 
and let $p(\bar{x})$ be a complete $(\tau, \lambda)$-neighbourhood type such that for all $i, j = 1, \ldots, k$,
$x_i \sim_p x_j$. 
Also suppose that $p$ is cofinally satisfiable in $\mbB$.
If, for some $i \in \{1, \ldots, k\}$,
\[
p(\bar{x}) \models \text{ ``the distance from $x_i$ to the root'' is less than $\lambda$''}
\]
then (as $x_i \sim_p x_j$ for all $i$ and $j$) it follows that there is $m \in \mbbN$ such that
$|p(\mcB_n)| \leq m$ for all $n$ so in this case $p$ is $p$ is bounded with respect to $\mbB$.
Otherwise, 
\begin{align}\label{p implies that no x-i is close}
\text{for all } i \in \{1, \ldots, k\}, \
p(\bar{x}) \models \text{ ``the distance from $x_i$ to the root is $\geq \lambda$''}
\end{align}
and for the rest of this example we assume~(\ref{p implies that no x-i is close}).
Observe that if $p$ is strongly unbounded with respect to $\mbB$ then $p$ is cofinally satisfiable in $\mbB$
and~(\ref{p implies that no x-i is close}) holds.
Condition~(4) of Assumption~\ref{properties of the base structures} and
Assumption~\ref{relative frequency of tau-closure types}
are consequences of the following condition:
\begin{enumerate}
\item[($\dagger$)] If $\lambda \in \mbbN$, $\bar{x} = (x_1, \ldots, x_k)$,  
and $p(\bar{x})$ is a cofinally satisfiable complete $(\tau, \lambda)$-neighbourhood type such that for all $i, j = 1, \ldots, k$,
$x_i \sim_p x_j$, and~(\ref{p implies that no x-i is close}) holds, then there is $\alpha_p > 0$ such that
$\lim_{n\to\infty} |p(\mcB_n)|/n = \alpha_p$.
\end{enumerate}
Hence it suffices to prove that $\mbB$ can be chosen so that ($\dagger$) holds.

We first establish a connection between a cofinally satisfiable complete $(\tau, \lambda)$-neigh\-bour\-hood type
$p$ satisfying the assumptions of ($\dagger$) and a $\delta$-bounded tree with a specified subset of its vertices.
For a tree $\mcT$ a {\em subtree} of $\mcT$ is a substructure of $\mcT$ that is also a tree
(i.e. an ordered rooted tree, but the root of the subtree need not coincide with the root of the tree that it is embedded in).
Suppose that $\mcT_0$ and $\mcT$ are $\delta$-bounded trees and that $A \subseteq T_0$.
We say that $\mcT_0$ is an {\em $A$-full subtree} of $\mcT$ if $\mcT_0$ is a subtree of $\mcT$ and
for all $a \in A$ every child of $a$ in $\mcT$ is also a child of $a$ in $\mcT_0$.
An embedding $f$ (in the model theoretic sense) of $\mcT_0$ into $\mcT$ is called {\em $A$-full} if
$\mcT \uhrc f(T_0)$ is an $f(A)$-full subtree of $\mcT$.

Now let $p(\bar{x})$ be as in ($\dagger$) and suppose that $\mcT \models p(\bar{a})$ where $\mcT$ is a 
$\delta$-bounded tree and
$\bar{a} \in T^k$. Then $\mcT_0 = \mcT \uhrc N_\lambda^{\mcT}(\bar{a})$ is a subtree of $\mcT$.
But there may be subtrees $\mcT' \subseteq \mcT$ that are isomorphic to $\mcT_0$ but where there is no $\bar{a}' \in T'$
that satisfies $p(\bar{x})$. 
The reason is that it $p(\bar{x})$ may express that some vertices in the $\lambda$-neigh\-bour\-hood
are leaves; then a subtree $\mcT' \subseteq \mcT$ isomorphic to $\mcT_0$ has a 
tuple of vertices satisfying $p$ if and only if every leaf of $\mcT'$ is also a leaf of $\mcT$.
It is not hard to see that one can associate a $\delta$-bounded tree $\mcT_p$ to $p$ and 
a (possibly empty) set of leaves $A \subseteq T_p$ such that 
\begin{enumerate}
\item[($\ddagger$)] for every $\delta$-bounded tree $\mcT$ and $\bar{a} \in T^k$,
$\mcT \models p(\bar{a})$ if and only if there is an $A$-full embedding $f$
of $\mcT_p$ into $\mcT$ such that $f(T_p) = N_\lambda^{\mcT}(\bar{a})$.
\end{enumerate}

Suppose that $\mcT_0$ and $\mcT$ are $\delta$-bounded trees.
Then let $\msfN_{\mcT_0}(\mcT)$ be the number of embeddings of $\mcT_0$ into $\mcT$
(i.e. the number of subtrees of $\mcT$ that are isomorphic to $\mcT_0$).
If $A \subseteq T_0$ then let $\msfN_{\mcT_0}^A(\mcT)$ be the number of $A$-full embeddings of $\mcT_0$ into $\mcT$.
From~($\ddagger$) it follows that~($\dagger$) is a consequence of the following condition:
\begin{align}\label{number of full embeddings}
&\text{For every $\delta$-bounded tree $\mcT_0$ and $A \subseteq T_0$ there is $\alpha_0 > 0$ such that}\\
&\lim_{n\to\infty} \frac{\msfN_{\mcT_0}^A(\mcB_n)}{n} = \alpha_0. \nonumber
\end{align} 

Next we will show that for all $\delta$-bounded trees $\mcT_0$ and $\mcT$ and any $A \subseteq T_0$,
$\msfN_{\mcT_0}^A(\mcT)$ can be expressed by using only addition, subtraction and 
terms $\msfN_{\mcT_1}(\mcT), \ldots, \msfN_{\mcT_m}(\mcT)$ where $\mcT_1, \ldots, \mcT_m$ are certain
trees that are determined (up to isomorphism) by $\mcT_0$ and $A$.
Given $\delta$-bounded trees $\mcT$, $\mcT'$ and $A \subseteq \mcT$, we call $\mcT'$ an 
{\em (proper) $A$-extension of $\mcT$}
if $\mcT$ is a (proper) subtree of $\mcT'$ and
every vertex of $\mcT'$ is a vertex of $\mcT$ or a child (in $\mcT'$) to a vertex in $A$.

Let $\mcT$ and $\mcT_0$ be $\delta$-bounded trees and let $A \subseteq T_0$.
Let $\mcT_1, \ldots, \mcT_m$ enumerate all,
up to isomorphism (and without repeating isomorphic trees), $\delta$-bounded trees that are proper $A$-extensions of $\mcT_0$.
(There are only finitely many such $A$-extensions.)
If some $a \in A$ has fewer than $\delta$ children in $\mcT_0$ then 
 the sequence $\mcT_1, \ldots, \mcT_m$ is nonempty and
\begin{equation}\label{an expression for N-A}
\msfN_{\mcT_0}^A(\mcT) = N_{\mcT_0}(\mcT) - \sum_{i=1}^m \msfN_{\mcT_i}^A(\mcT).
\end{equation}
On the other hand,
\begin{align}\label{the base case for number of copies of T-0}
&\text{if every $a \in A$ has $\delta$ children in $\mcT_0$ then $\msfN_{\mcT_0}^A(\mcT) = \msfN_{\mcT_0}(\mcT)$}\\
&\text{and there is {\em no} $\delta$-bounded tree which is a proper $A$-extension of $\mcT_0$.} \nonumber
\end{align}
Also note that if $\mcT'$ is a $\delta$-bounded tree and a proper $A$-extension of $\mcT_0$ then
$\mcT'$ has {\em fewer} (counting up to isomorphism) proper $A$-extensions (that are $\delta$-bounded) than $\mcT_0$.
Therefore it follows from~(\ref{an expression for N-A}), (\ref{the base case for number of copies of T-0}) and induction 
on the number (up to isomorphism) of $\delta$-bounded trees that are $A$-extensions of $\mcT_0$ that
$\msfN_{\mcT_0}^A(\mcT)$ can be written as an expression involving only
$\msfN_{\mcT_1}(\mcT), \ldots, \msfN_{\mcT_m}(\mcT)$ and `$+$' and `$-$'.
It follows that~(\ref{number of full embeddings}) is a consequence of the following condition:
\begin{align}\label{number of embeddings}
&\text{For every $\delta$-bounded tree $\mcT_0$ there is $\alpha_0 > 0$ such that}\\
&\lim_{n\to\infty} \frac{\msfN_{\mcT_0}(\mcB_n)}{n} = \alpha_0. \nonumber
\end{align} 

Thus it remains to show that the sequence $\mbB = (\mcB_n : n \in \mbbN^+)$ of $\delta$-bounded trees can be
chosen so that~(\ref{number of embeddings}) holds.
Note that such $\mbB$ also satisfies condition~(4) of Assumption~\ref{properties of the base structures}. 
Let $\mbT_n$ be the set which contains exactly one representative from every
isomorphism class of ordered trees with exactly $n$ vertices.
Let $X$ be a random variable with values in $\mbbN$.
We can think of $X$ as being the number of children that a vertex in a tree has.
Choose reals $c_i \in [0, 1]$ for $i \in \mbbN$ such that $c_i > 0$ if $i \leq \delta$ and $c_i = 0$ otherwise.
Define the so-called {\em offspring distribution} by
$\mbbP(X = i) = c_i$ for all $i \in \mbbN$.
To every tree $\mcT$ we associate the {\em weight}
$\phi(\mcT) = \sum_{v \in T} c(v)$ where $c(v) = c_i$ if $v$ has exactly $i$ children,
and note that only $\delta$-bounded trees have positive weight.
Define 
$\gamma_n = \sum_{\mcT \in \mbT_n} \phi(\mcT)$ and
$\mbbP_n(\mcT) = \phi(\mcT)/\gamma_n$ for every $\mcT \in \mbT_n$.
Then $\mbbP_n$ is a probability distribution on $\mbT_n$.
Suppose that $\mbbE(X) = 1$,
the so called critical case for Galton-Watson trees.
(As remarked upon in \cite[Remark~3.1]{Jan16} this is only a minor restriction from a probability theoretic point of view.)
Now the main result in \cite{Jan21} implies the following:
{\em For every $\delta$-bounded tree $\mcT_0$ there is $\alpha_0 > 0$ such that for all $\varepsilon > 0$,
if $n$ is large enough then}
\[
\mbbP_n\Big(\Big\{\mcT \in \mbT_n : \alpha_0 - \varepsilon \leq 
\frac{\msfN_{\mcT_0}(\mcT)}{n} 
\leq \alpha_0 + \varepsilon \Big\}\Big)
\geq 1 - \varepsilon. 
\]
From this it follows that it is possible to choose $\mcB_n$, for all $n \in \mbbN^+$, 
so that~(\ref{number of embeddings}) holds. (In fact, with respect to the sequence of 
distributions $\mbbP_n$, most choices will lead to this).
}\end{exam}

\noindent
Before stating and proving the main results of this section we use
Assumption~\ref{relative frequency of tau-closure types} to prove a sequence of increasingly more general
results about balanced triples, culminating 
with Proposition~\ref{balance under stronger conditions for sigma-types}, 
which will be used in the proofs of the main results of this section.

\begin{lem}\label{result about balance for tau-closure types}
Let $\lambda \leq \mu \in \mbbN$, 
let $p(\bar{x}, \bar{y})$ be a complete $(\tau, \lambda)$-closure type,
let $r(\bar{x}, \bar{y})$ be a complete $(\tau, \mu)$-closure type which is consistent with $p$
 and suppose that
$p(\bar{x}, \bar{y})$ is cofinally satisfiable.
Then there is $\gamma \in \mbbN$ such that if $q(\bar{x})$ is a complete $(\tau, \gamma)$-closure type,
then $(r, p, q)$ is balanced.
\end{lem}

\noindent
{\bf Proof.}
Let $\lambda \leq \mu \in \mbbN$ and let $p(\bar{x}, \bar{y})$ and  $r(\bar{x}, \bar{y})$ be as assumed in the lemma.
Then we can as well assume that $p = r \uhrc \lambda$.
We use induction on $\dim_{\bar{y}}(p)$.
If $\dim_{\bar{y}}(p) = 0$, that is, if $p$ is $\bar{y}$-bounded then the conclusion of the lemma follows
from Lemma~\ref{balanced bounded triples}.

Now suppose that $\dim_{\bar{y}}(p) \geq k+1$ where $k \in \mbbN$.
By Lemma~\ref{getting a strongly unbounded type of dimension 1},
we may assume that $\bar{y} = \bar{u}\bar{v}$,
$p$ is strongly $\bar{v}$-unbounded, $\dim_{\bar{v}}(p) = 1$.
By Definition~\ref{definition of dimension, part 2} of dimension
and Lemma~\ref{another characterization of strong unboundedness},
if $p_1(\bar{x}, \bar{u}) = p \uhrc \bar{x}\bar{u}$ then $\dim_{\bar{u}}(p_1) = k$.
Let $r_1(\bar{x}, \bar{u}) = r \uhrc \bar{x}\bar{u}$.
By the induction hypothesis there is $\gamma_1 \in \mbbN$ such that
if $q(\bar{x})$ is a complete $(\tau, \gamma_1)$-closure type,
then $(r_1, p_1, q)$ is $\alpha$-balanced for some $\alpha \in [0, 1]$.

So let $q(\bar{x})$ be a complete $(\tau, \gamma_1)$-closure type and suppose that
$(r_1, p_1, q)$ is $\alpha$-balanced.
If $(r, p, r_1)$ is $\beta$-balanced then it follows that $(r, p, q)$ is $\alpha\beta$-balanced,
so it suffices to show that $(r, p, r_1)$ is balanced (although in the first case below we show directly
that $(r, p, q)$ is balanced).

First suppose that $r$ is $\bar{v}$-bounded, so there is $m \in \mbbN$ such that
for all $n$, $\bar{a} \in (B_n)^{|\bar{x}|}$, and $\bar{c} \in (B_n)^{|\bar{u}|}$,
$r(\bar{a}, \bar{c}, \mcB_n) \leq m$.
Since $p$ is strongly $\bar{v}$-unbounded (hence uniformly $\bar{v}$-unbounded)
it follows 
(from Definition~\ref{definition of unbounded formula})
that there is $f : \mbbN \to \mbbN$ such that $\lim_{n\to\infty}f(n) = \infty$ and for all $n$,
$\bar{a} \in (B_n)^{|\bar{x}|}$, and $\bar{c} \in (B_n)^{|\bar{u}|}$,
if $|p(\bar{a}, \bar{c}, \mcB_n)| \neq \es$ then
$|p(\bar{a}, \bar{c}, \mcB_n)| \geq f(n)$.
Since $r$ implies $p$ it follows that 
for every $\varepsilon > 0$ we have
\[
|r(\bar{a}, \bar{c}, \mcB_n)| \leq \varepsilon |p(\bar{a}, \bar{c}, \mcB_n)| \quad \text{  if $n$ is large enough}.
\]
Hence $(r, p, q)$ is $0$-balanced.

Now suppose that $r$ is not $\bar{v}$-bounded (and we will show that $(r, p, r_1)$ is balanced).
By Lemmas~\ref{not bounded implies uniformly unbounded}
and~\ref{dimension and unboundedness},
$r$ is then uniformly $\bar{v}$-unbounded and $\dim_{\bar{v}}(r) \geq 1$.
Since  $r$  implies
$p$ it follows (from the definitions of $\sim_p$ and $\bar{v}$-dimension)
that $\dim_{\bar{v}}(r) \leq \dim_{\bar{v}}(p) = 1$,
so $\dim_{\bar{v}}(r) = 1$.

Next we show that $r$ is strongly $\bar{v}$-unbounded.
Since $p$ is strongly $\bar{v}$-unbounded and $\dim_{\bar{v}}(p) = 1$ 
it follows 
(from Definition~\ref{definition of dimension, part 2} of dimension
and Lemma~\ref{another characterization of strong unboundedness}) 
that all variables in $\bar{v}$ belong to the
same $\sim_p$-class. Since $\mu \geq \lambda$ and $p = r \uhrc \lambda$ it follows that
all variables in $\bar{v}$ belong to the same $\sim_r$-class.
Suppose for a contradiction that $r$ is not strongly $\bar{v}$-unbounded.
Then (by definition of strong unboundedness) there is a subsequence $\bar{v}'$ of $\bar{v}$ such that
$r \uhrc \bar{v}'$ is not uniformly $\bar{v}'$-unbounded.
By Lemma~\ref{not bounded implies uniformly unbounded},
$r$ is $\bar{v}'$-bounded. 
As all variables of $\bar{v}$ belong to the same $\sim_r$-class it follows that $r$ is $\bar{v}$-bounded, contradicting
our assumption.

Let $p_2(\bar{v}) = p \uhrc \bar{v}$ and $r_2(\bar{v}) = r \uhrc \bar{v}$.
Since $p$ and $r$ are strongly $\bar{v}$-unbounded 
it follows from 
Lemma~\ref{elementary properties of closure types}~(ii)
that $p_2(\bar{v})$ and $r_2(\bar{v})$ are strongly unbounded.
Since $p_2 = p \uhrc \bar{v}$ and all variables in $\bar{v}$ are in the same $\sim_p$-class it follows 
that all variables in $\bar{v}$ are in the same $\sim_{p_2}$-class.
By similar reasoning all variables in $\bar{v}$ are in the same $\sim_{r_2}$-class.

As $p$ and $r$ are strongly $\bar{v}$-unbounded it follows from
Lemma~\ref{strong unboundedness implies large distance}
that 
\begin{align*}
&p(\bar{x}, \bar{u}, \bar{v}) \models \dist(\bar{x}\bar{u}, \bar{v}) > 2\lambda \ \wedge \
\forall z \big(\text{``$z$ is $\lambda$-rare''} \rightarrow \dist(z, \bar{v}) > 2\lambda \big)
\text{ and } \\
&r(\bar{x}, \bar{u}, \bar{v}) \models \dist(\bar{x}\bar{u}, \bar{v}) > 2\mu  \ \wedge \
\forall z \big(\text{``$z$ is $\mu$-rare''} \rightarrow \dist(z, \bar{v}) > 2\mu \big).
\end{align*}
By Lemma~\ref{splitting a closure type}
there are a complete $(\tau, \lambda)$-{\em neigh\-bour\-hood type} $p'_2(\bar{v})$ and a
complete $(\tau, \mu)$-{\em neigh\-bour\-hood type} $r'_2(\bar{v})$ such that $p_2 \models p'_2$,
$r_2 \models r'_2$,
\begin{align}\label{characterization of satisfaction of p}
&\text{$\mcB_n \models p(\bar{a}, \bar{c}, \bar{d})$ if and only if
$\mcB_n \models p_1(\bar{a}, \bar{c}) \wedge p'_2(\bar{d})$ and
$\dist(\bar{e}\bar{a}\bar{c}, \bar{d}) > 2\lambda$,} \\
&\text{where $\bar{e}$ enumerates all $\lambda$-rare elements, and} \nonumber
\end{align}
\begin{align}\label{characterization of satisfaction of r}
&\text{$\mcB_n \models r(\bar{a}, \bar{c}, \bar{d})$ if and only if
$\mcB_n \models r_1(\bar{a}, \bar{c}) \wedge r'_2(\bar{d})$ and
$\dist(\bar{e}\bar{a}\bar{c}, \bar{d}) > 2\mu$,} \\
&\text{where $\bar{e}$ enumerates all $\mu$-rare elements.} \nonumber
\end{align}

\medskip

\noindent
{\bf Claim:} $p'_2(\bar{v})$ and $r'_2(\bar{v})$ are strongly unbounded and both have dimension 1.

\medskip
\noindent
{\bf Proof of the claim:}
We prove the claim only for $p'_2$ since the proof for $r'_2$ is identical, besides replacing $p_2$ and $p'_2$
by $r_2$ ad $r'_2$, respectively.
Suppose for a contradiction that $p'_2(\bar{v})$ is not strongly unbounded.
Then there is a subsequence $\bar{v}'$ of $\bar{v}$ such that $p'_2 \uhrc \bar{v}'$ is not uniformly unbounded,
and hence $p'_2 \uhrc \bar{v}'$ is bounded. 
From $p_2(\bar{v}) \models p'_2(\bar{v})$ we get $p_2 \uhrc \bar{v}' \models p'_2 \uhrc \bar{v}'$
and it follows that $p_2 \uhrc \bar{v}'$ is bounded which contradicts that $p_2$ is strongly unbounded
(which we have already proved).

By Lemma~\ref{connection between closure types and neighbourhood types}
there is a complete 
$(\tau, \lambda)$-neighbourhood type $p^+_2(\bar{z}, \bar{v})$ such that, if
$\varphi_\lambda(w)$ expresses that ``$w$ is $\lambda$-rare'', then
\begin{align}\label{p-2 equivalent to exists z such that p-+-2}
&p_2(\bar{v}) \text{ is equivalent to } 
\exists \bar{z} \big( \forall w\big(\varphi_\lambda(w) \to w \in \rng(\bar{z})\big) 
\wedge p^+_2(\bar{z}, \bar{v})\big), \text{ and }\\
&p^+_2 \uhrc \bar{z} \text{ is bounded.}  \nonumber
\end{align}
Since $p_2(\bar{v}) \models p'_2(\bar{v})$ it follows from~(\ref{p-2 equivalent to exists z such that p-+-2})
that $p'_2$ and $p^+_2 \uhrc \bar{v}$ are equivalent.
Hence we can assume that $p'_2(\bar{v}) = p^+_2 \uhrc \bar{v}$.
Suppose, towards a contradiction, that $z_i \sim_{p^+_2} v_j$ for some $z_i \in \rng(\bar{z})$ and $v_j \in \rng(\bar{v})$.
Then $p^+_2 \uhrc \bar{z}v_j$ is $v_j$-bounded and as $p^+_2 \uhrc \bar{z}$ is bounded
(by~(\ref{p-2 equivalent to exists z such that p-+-2}))
it follows that $p_2 \uhrc v_j$ is bounded which contradicts that $p_2$ is strongly unbounded (proved above).
Hence $z_i \not\sim_{p^+_2} v_j$ for all $z_i \in \rng(\bar{z})$ and $v_j \in \rng(\bar{v})$.

Above we proved that all variables in $\bar{v}$ belong to the same $\sim_{p_2}$-class.
From~(\ref{p-2 equivalent to exists z such that p-+-2})
it now follows that all variables in $\bar{v}$ belong to the same $\sim_{p'_2}$-class, so $\dim(p'_2) \leq 1$.
Since $p'_2$ is strongly unbounded it follows from the definition of dimension that 
$\dim(p'_2) \geq 1$, so $\dim(p'_2) = 1$.
\hfill $\square$

\medskip

\noindent
By the choices of $p'_2$ and $r'_2$, $p'_2$ must be equivalent to $r'_2 \uhrc \lambda$.
Hence we can assume that $p'_2 = r'_2 \uhrc \lambda$.
The claim and Assumption~\ref{relative frequency of tau-closure types} implies that there is $\beta \in [0, 1]$
such that 
\begin{equation}\label{r-1 divided by p-1}
\lim_{n\to\infty} \frac{|r'_2(\mcB_n)|}{|p'_2(\mcB_n)|} = \beta.
\end{equation}

Recall that by Assumption~\ref{properties of the base structures},
there is a fixed number $\Delta$ such that all $\mcB_n$ have degree at most $\Delta$.
Let $\bar{e}$ enumerate all $\mu$-rare (and hence also $\lambda$-rare) elements of $\mcB_n$.
Suppose that $\mcB_n \models r_1(\bar{a}, \bar{c})$.
Since (by the claim) $\dim(p'_2) = \dim(r'_2) = 1$
it follows that there is a constant $K \in \mbbN$ (independent of $n$) such that 
at most $K$ different members of $p'_2(\mcB_n) \cup r'_2(\mcB_n)$ ($= p'_2(\mcB_n)$) are within distance $2\mu$ 
from any member of $\bar{e}\bar{a}\bar{c}$.
Since (by the claim) $p'_2$ and $r'_2$ are strongly unbounded (complete neighbourhood types) with dimension 1 it follows from 
Assumption~\ref{properties of the base structures} that
$\lim_{n\to\infty}|p'_2(\mcB_n)| = \lim_{n\to\infty}|r'_2(\mcB_n)| = \infty$.
From this together with~(\ref{characterization of satisfaction of p}), (\ref{characterization of satisfaction of r})
and~(\ref{r-1 divided by p-1}) it follows that for every $\varepsilon > 0$
\[
(\beta - \varepsilon)|p(\bar{a}, \bar{c}, \mcB_n)| \leq |r(\bar{a}, \bar{c}, \mcB_n)|
\leq (\beta + \varepsilon)|p(\bar{a}, \bar{c}, \mcB_n)| \ \ \text{ if $n$ is large enough.}
\]
It follows that $(r, p, r_1)$ is $\beta$-balanced. 
This completes the proof.
\hfill $\square$

\begin{lem}\label{balance under stronger assumptions} 
Let $\lambda, \mu \in \mbbN$, 
let $p(\bar{x}, \bar{y})$ be a complete $(\sigma, \lambda)$-closure type
and let $p_\tau(\bar{x}, \bar{y})$ be a uniformly $\bar{y}$-unbounded complete $(\tau, \mu)$-closure type.
There is $\gamma \in \mbbN$ such that if $q(\bar{x})$ is a complete $(\sigma, \gamma)$-closure type,
then $(p, p_\tau, q)$ is balanced.
\end{lem}

\noindent
{\bf Proof.}
Let $p(\bar{x}, \bar{y})$ and $p_\tau(\bar{x}, \bar{y})$ be as assumed in the lemma.
We assume that $\kappa$ and $\kappa'$ are as assumed in 
Assumption~\ref{induction hypothesis}.
If $\mu \geq \lambda + \max(\lambda, \kappa + \kappa')$ then the conclusion of the lemma follows from
Lemma~\ref{balanced triples, uniformly unbounded}.

So suppose that $\mu < \lambda + \max(\lambda, \kappa + \kappa')$.
Let $p_{\tau, i}(\bar{x}, \bar{y})$, $i = 1, \ldots, s$, enumerate all, up to equivalence,
complete $(\tau, \lambda + \max(\lambda, \kappa + \kappa'))$-closure types that imply $p_\tau$
(or in other words such that $p_\tau = p_{\tau, i} \uhrc \mu$).
Note that since $p_\tau$ is assumed to be uniformly $\bar{y}$-unbounded 
it follows that $p_\tau$ is cofinally satisfiable.
By 
Lemma~\ref{result about balance for tau-closure types}
and~\ref{on stronger conditioning in the balance},
there is $\gamma \in \mbbN$ such that if $q(\bar{x})$ is a complete $(\tau, \gamma)$-closure type,
then $(p_{\tau, i}, p_\tau, q)$ is balanced for all $i = 1, \ldots s$.
So there are $\alpha_1, \ldots, \alpha_s \in [0, 1]$ such that
$(p_{\tau, i}, p_\tau, q)$ is $\alpha_i$-balanced for all $i = 1, \ldots, s$.
By Lemma~\ref{on stronger conditioning in the balance} and
Lemma~\ref{balanced triples, uniformly unbounded},
we may assume that $\gamma$ is large enough such that there are 
$\beta_1, \ldots, \beta_s$ such that if $q(\bar{x})$ is a complete 
$(\tau, \gamma)$-closure type, then $(p, p_{\tau, i}, q)$ is $\beta_i$-balanced.
It follows that if $\alpha = \sum_{i=1}^s \alpha_i\beta_i$,
then $(p, p_\tau, q)$ is $\alpha$-balanced.
\hfill $\square$

\begin{lem}\label{general statement about balance under stronger assumptions}
Let $\lambda, \mu \in \mbbN$, 
let $p(\bar{x}, \bar{y})$ be a complete $(\sigma, \lambda)$-closure type
and let $p_\tau(\bar{x}, \bar{y})$ be a  complete $(\tau, \mu)$-closure type.
There is $\gamma \in \mbbN$ such that if $q(\bar{x})$ is a complete $(\sigma, \gamma)$-closure type,
then $(p, p_\tau, q)$ is balanced.
\end{lem}

\noindent
{\bf Proof.}
If $p_\tau$ is $\bar{y}$-bounded then the result follows from
Lemma~\ref{balanced bounded triples}.
If $p_\tau$ is not $\bar{y}$-bounded, then 
(by Lemma~\ref{not bounded implies uniformly unbounded})
it is uniformly $\bar{y}$-unbounded and the result now follows from 
Lemma~\ref{balance under stronger assumptions}.
\hfill $\square$

\begin{prop}\label{balance under stronger conditions for noncomplete types}
Let $\lambda, \mu \in \mbbN$, 
let $p(\bar{x}, \bar{y})$ be a (not necessarily complete) $(\sigma, \lambda)$-closure type
and let $p_\tau(\bar{x}, \bar{y})$ be a  complete $(\tau, \mu)$-closure type.
Then there is $\gamma \in \mbbN$ such that if $q(\bar{x})$ is a complete $(\sigma, \gamma)$-closure type,
then $(p, p_\tau, q)$ is balanced.
\end{prop}

\noindent
{\bf Proof.}
We can argue essentially as in the proof of
Proposition~\ref{incomplete balanced triples, uniformly unbounded}
and use 
Lemma~\ref{general statement about balance under stronger assumptions}
when, in the proof of Proposition~\ref{incomplete balanced triples, uniformly unbounded},
Lemma~\ref{balanced triples, uniformly unbounded}
was used.
The details are left for the reader.
\hfill $\square$

\begin{prop}\label{balance under stronger conditions for sigma-types}
Suppose that $\lambda, \mu \in \mbbN$,
$p_1(\bar{x}, \bar{y})$ is a $(\sigma, \lambda)$-closure type, and
$p_2(\bar{x}, \bar{y})$ is a $(\sigma, \mu)$-closure type that is $\bar{y}$-positive.
Then there is $\gamma \in \mbbN$ such that if $q(\bar{x})$ is a complete $(\sigma, \gamma)$-closure type,
then $(p_1, p_2, q)$ is balanced.
\end{prop}

\noindent
{\bf Proof.}
Let $p_1(\bar{x}, \bar{y})$ and $p_2(\bar{x}, \bar{y})$ be as assumed in the lemma.
Since $p_2$ is $\bar{y}$-positive there is $\gamma_1$ such that if $q(\bar{x})$ is a complete 
$(\sigma, \gamma_1)$-closure type and $p_2 \uhrc \tau \wedge q$ is cofinally satisfiable, then 
$(p_2, p_2 \uhrc \tau, q)$ is $\beta$-balanced for some $\beta > 0$.
By Proposition~\ref{balance under stronger conditions for noncomplete types},
there is $\gamma_2$ such that if $q(\bar{x})$ is a complete $(\sigma, \gamma_2)$-closure type,
then $(p_1 \wedge p_2, p_2 \uhrc \tau, q)$ is $\alpha$-balanced for some $\alpha$.
Let $\gamma = \max\{\gamma_1, \gamma_2\}$.
It follows from Lemma~\ref{on stronger conditioning in the balance}
that if $q(\bar{x})$ is a complete $(\sigma, \gamma)$-closure type such that $p_2 \uhrc \tau \wedge q$ is cofinally satisfiable,
then $(p_2, p_2 \uhrc \tau, q)$ is $\beta$-balanced where $\beta > 0$
and $(p_1 \wedge p_2, p_2 \uhrc \tau, q)$ is $\alpha$-balanced;
hence $(p_1, p_2, q)$ is $\alpha/\beta$-balanced.
If $p_2 \uhrc \tau \wedge q$ is not cofinally satisfiable then 
(by Lemma~\ref{balance follows from inconsistency})
$(p_1, p_2, q)$ is 0-balanced.
\hfill $\square$

\bigskip

\noindent
Now we are ready to show (again) that aggregation functions can be asymptotically eliminated under certain assymptions.
Proposition~\ref{asymptotic elimination under stronger conditions} below has stronger asumptions on the base sequence
of structures than 
Proposition~\ref{asymptotic elimination for strongly unbounded aggregations},
but Proposition~\ref{asymptotic elimination under stronger conditions}
applies to more $PLA^*(\sigma)$-networks and more formulas (for expressing queries) than
Proposition~\ref{asymptotic elimination for strongly unbounded aggregations}.

\begin{prop}\label{asymptotic elimination under stronger conditions}
Suppose that 
Assumptions~\ref{properties of the base structures}, \ref{induction hypothesis},
and~\ref{relative frequency of tau-closure types}
are satisfied.
Let $\varphi(\bar{x}) \in PLA^*(\sigma)$ and suppose that if
\[
F\big(\varphi_1(\bar{y}, \bar{z}), \ldots, \varphi_m(\bar{y}, \bar{z}) : \bar{z} : 
\chi_1(\bar{y}, \bar{z}), \ldots, \chi_m(\bar{y}, \bar{z})\big)
\]
is a subformula of $\varphi(\bar{x})$ (where $F$ is an aggregation function) then either
\begin{enumerate}
\item for all $i = 1, \ldots, m$, $\varphi_i(\bar{y}, \bar{z})$ is finite valued and there is $\xi_i \in \mbbN$
such that for all $n$, $\mcA \in \mbW_n$, $\bar{a} \in (B_n)^{|\bar{y}|}$, and $\bar{b} \in (B_n)^{|\bar{z}|}$,
if $\mcA \models \chi_i(\bar{a}, \bar{b})$ then $\rng(\bar{b}) \subseteq C_{\xi_i}^{\mcB_n}(\bar{a})$, or

\item for all $i = 1, \ldots, m$, $\chi_i(\bar{y}, \bar{z})$ is a $\bar{z}$-positive $(\sigma, \lambda_i)$-closure type
for some $\lambda_i \in \mbbN$ and
$F$ is continuous.
\end{enumerate}
Then there is $\xi \in \mbbN$ such that $\varphi(\bar{x})$ is asymptotically equivalent to a
$(\sigma, \xi)$-basic formula.

If, in addition, conditions~(I) and~(II) of
Remark~\ref{remark about basic probability formulas}
hold, then we can weaken condition~(2) above by replacing `continuous' by `admissible' and the conclusion
that $\varphi(\bar{x})$ is asymptotically equivalent to a $(\sigma, \xi)$-basic formula still follows.
\end{prop}

\noindent
{\bf Proof.}
We use induction on the complexity of formulas.
If $\varphi(\bar{x})$ is aggregation-free, then the conlusion follows from 
Lemma~\ref{connectives and basic formulas}~(ii).
Suppose that $\varphi(\bar{x})$ has the form
$\msfC(\psi_1(\bar{x}), \ldots, \psi_k(\bar{x}))$
where $\msfC : [0, 1]^k \to [0, 1]$ is a continuous
and each $\psi_i(\bar{x})$ is asymptotically equivalent (with respect to $(\mbB, \mbbG)$)
to a $(\sigma, \lambda)$-basic formula $\psi'_i(\bar{x})$.
Since $\msfC$ is continuous it follows that $\varphi(\bar{x})$ and
$\msfC(\psi'_1(\bar{x}), \ldots, \psi'_k(\bar{x}))$ are asymptotically equivalent.
By Lemma~\ref{connectives and basic formulas}~(i),
$\msfC(\psi'_1(\bar{x}), \ldots, \psi'_k(\bar{x}))$ is equivalent, hence asymptotically equivalent,
to a $(\sigma, \lambda)$-basic formula $\varphi'(\bar{x})$.
By transitivity of asymptotic equivalence it follows that $\varphi$ and $\varphi'$ are asymptotically equivalent.

Now suppose that $\varphi(\bar{x})$ has the form
\begin{equation}\label{an aggregation function to be eliminated}
F\big(\varphi_1(\bar{x}, \bar{y}), \ldots, \varphi_m(\bar{x}, \bar{y}) : \bar{y} : 
\chi_1(\bar{x}, \bar{y}), \ldots, \chi_m(\bar{x}, \bar{y})\big).
\end{equation}
By the induction hypothesis, 
for each $i = 1, \ldots, m$, there is $\lambda_i \in \mbbN$ and a
$(\sigma, \lambda_i)$-basic formula $\psi_i(\bar{x}, \bar{y})$ that is
asymptotically equivalent to $\varphi_i(\bar{x}, \bar{y})$.

First suppose that~(1) holds, that is, for all $i$, $\varphi_i(\bar{x}, \bar{y})$ is finite valued and there is $\xi_i \in \mbbN$
such that for all $n$, $\mcA \in \mbW_n$, $\bar{a} \in (B_n)^{|\bar{x}|}$, and $\bar{b} \in (B_n)^{|\bar{y}|}$,
if $\mcA \models \chi_i(\bar{a}, \bar{b})$ then $\rng(\bar{b}) \subseteq C_{\xi_i}^{\mcB_n}(\bar{a})$.
Note that since $\psi_i(\bar{x}, \bar{y})$ is a $(\sigma, \lambda_i)$-basic formula it is finitely valued.
Since $\psi_i$ and $\varphi_i$ are asymptotically equivalent and both are finite valued it follows that there are
$\mbY_{n, i} \subseteq \mbW_n$ for all $n$ such that $\lim_{n\to\infty}\mbbP_n(\mbY_{n, i}) = 1$ and
for all $\mcA \in \mbY_{n, i}$, $\bar{a} \in (B_n)^{|\bar{x}|}$ and $\bar{b} \in (B_n)^{|\bar{y}|}$,
$\mcA(\varphi_i(\bar{a}, \bar{b})) = \mcA(\psi_i(\bar{a}, \bar{b}))$.

Let $\varphi'(\bar{x})$ denote the formula
\begin{equation}\label{another aggregation to be eliminated}
F\big(\psi_1(\bar{x}, \bar{y}), \ldots, \psi_m(\bar{x}, \bar{y}) : \bar{y} : 
\chi_1(\bar{x}, \bar{y}), \ldots, \chi_m(\bar{x}, \bar{y})\big).
\end{equation}
Let $\xi = \max(\xi_1 + \lambda_1, \ldots, \xi_m + \gamma_m)$ and
let $q_1(\bar{x}), \ldots, q_s(\bar{x})$ enumerate all, up to equivalence, 
complete $(\sigma, \xi)$-closure types in the variables $\bar{x}$.
Then, for all $i = 1, \ldots, s$, all finite $\sigma$-structures $\mcA$, $\mcA'$ 
and all $\bar{a} \in A^{|\bar{x}|}$ and $\bar{a}' \in (A')^{|\bar{x}|}$,
if $\mcA \models q_i(\bar{a})$ and $\mcA' \models q_i(\bar{a}')$, then $\mcA(\varphi'(\bar{a})) = \mcA'(\varphi'(\bar{a}'))$.
In other words, the value $\mcA(\varphi'(\bar{a}))$ depends only on which $q_i$ $\bar{a}$ satisfies.
For all $i = 1, \ldots, m$ let $c_i = \mcA(\varphi(\bar{a}))$ where $\mcA$ and $\bar{a}$ are such that 
$\mcA \models q_i(\bar{a})$.
Then $\varphi'(\bar{x})$ is equivalent to the $(\sigma, \xi)$-basic formula $\bigwedge_{i=1}^s (q_i(\bar{x}) \to c_i)$.
From the choice of $\mbY_{n, i}$ it follows that for all $n$, $\mcA \in \bigcap_{i=1}^m \mbY_{n, i}$,
and $\bar{a} \in (B_n)^{|\bar{x}|}$, $\mcA(\varphi(\bar{a})) = \mcA(\varphi'(\bar{a}))$.
Since $\lim_{n\to\infty}\mbbP_n\big(\bigcap_{i=1}^m \mbY_{n, i}\big) = 1$ it follows that 
$\varphi(\bar{x})$ is asymptotically equivalent to $\varphi'(\bar{x})$ and hence also to 
$\bigwedge_{i=1}^s (q_i(\bar{x}) \to c_i)$.

Now suppose that~(2) holds, that is,
for all $i = 1, \ldots, m$, $\chi_i(\bar{x}, \bar{y})$ is a $\bar{y}$-positive $(\sigma, \lambda_i)$-closure type
for some $\lambda_i \in \mbbN$ and
$F$ is continuous.
Let $L_0$ be the set of all $\varphi \in PLA^*(\sigma)$ such that, for some $\lambda \in \mbbN$, 
$\varphi$ is a {\em complete} $(\sigma, \lambda)$-closure type in some sequence of variables.
Let $L_1$ be the set of all $\varphi \in PLA^*(\sigma)$ such that  for some $\lambda \in \mbbN$, 
$\varphi$ is a $(\sigma, \lambda)$-closure type in some sequence of variables.
For every $\varphi(\bar{x}, \bar{y}) \in L_0$ let $L_{\varphi(\bar{x}, \bar{y})}$ be the set
of all $\chi(\bar{x}, \bar{y}) \in L_1$ such that, for some $\gamma \in \mbbN$, $ \chi$ is a $\bar{y}$-positive
$(\sigma, \gamma)$-closure type.

Due to Theorem~\ref{general asymptotic elimination}
it now suffices to show that 
Assumption~\ref{assumptions on the basic logic}
is satisfied.
Part~(1)
of Assumption~\ref{assumptions on the basic logic} 
follows from
Lemma~\ref{connectives and basic formulas}~(ii),
so we verify part~(2) of the same assumption.
So, for $i = 1, \ldots, k$, 
let $p_i(\bar{x}, \bar{y})$ be a complete $(\sigma, \mu_i)$-closure type and
let $\chi_i(\bar{x}, \bar{y})$ be a $\bar{y}$-positive $(\sigma, \xi_i)$-closure type.
From
Lemma~\ref{on stronger conditioning in the balance}
and Proposition~\ref{balance under stronger conditions for sigma-types}
it follows that there is a $\mu\in \mbbN$ such that if $q(\bar{x})$ is a complete $(\sigma, \mu)$-closure type,
then, for each $j$, $(p_j, \chi_j, q)$ is balanced.

Let $q_1(\bar{x}), \ldots, q_s(\bar{x})$ enumerate, up to equivalence, all complete $(\sigma, \mu)$-closure types.
Then $(p_j, \chi_j, q_i)$ is balanced for every choice of $i$ and $j$. 
It follows that there are
$\alpha_{i, j} \in [0, 1]$ such that for every $\varepsilon > 0$ and $n$ there is $\mbY_n^\varepsilon \subseteq \mbW_n$
such that $\lim_{n\to\infty}\mbbP_n(\mbY_n^\varepsilon) = 1$ and 
for all $i, j$, all $\mcA \in \mbY_n^\varepsilon$,
and all $\bar{a} \in (B_n)^{|\bar{x}|}$, if $\mcA \models q_i(\bar{a})$ then
\[
(\alpha_{i, j} - \varepsilon)|\chi_j(\bar{a}, \mcA)| \ \leq \ 
|p_j(\bar{a}, \mcA) \cap \chi_j(\bar{a}, \mcA)| \ \leq \ (\alpha_{i, j} + \varepsilon)|\chi_j(\bar{a}, \mcA)|.
\]
Thus condition~(d) of part~(2) of Assumption~\ref{assumptions on the basic logic} is satisfied.
It is also clear from the choice of the $q_i$ that for all $n$ and $\mcA \in \mbW_n$,
\[
\mcA \models \forall \bar{x} \bigvee_{i=1}^s q_i(\bar{x}) \ \text{ and } \
\mcA \models \forall \bar{x} \neg (q_i(\bar{x}) \wedge q_j(\bar{x})) \ \text{ if $i \neq j$}.
\]
So conditions~(a) and~(b) of part~(2) of Assumption~\ref{assumptions on the basic logic} are satisfied.
Let $\hat{\chi}_i(\bar{x}, \bar{y}) = \chi_i \uhrc \tau$ for $i = 1, \ldots, k$.
It follows from 
Lemma~\ref{an x-type that determines if y exists}
that there are $\xi \in \mbbN$ and complete $(\tau, \xi)$-closure types $\chi^*_1(\bar{x}), \ldots, \chi^*_{t_0}(\bar{x})$
such that for all sufficiently large $n$, 
\[
\mcB_n \models \Big(\bigvee_{i=1}^m \neg \exists \bar{y} \hat{\chi}_i(\bar{x}, \bar{y}) \Big) \leftrightarrow
\Big(\bigvee_{i=1}^{t_0} \chi^*_i(\bar{x})\Big).
\]
By considering all, up to equivalence, complete $(\sigma, \xi)$-types in the variables $\bar{x}$ which are consistent with 
some $\chi^*_i(\bar{x})$ we find a sequence of complete $(\sigma, \xi)$-types $\chi'_1(\bar{x}), \ldots, \chi'_t(\bar{x})$
such that for all sufficiently large $n$ and all $\mcA \in \mbW_n$,
\[
\mcA \models \Big(\bigvee_{i=1}^m \neg \exists \bar{y} \hat{\chi}_i(\bar{x}, \bar{y}) \Big) \leftrightarrow
\Big(\bigvee_{i=1}^t \chi'_i(\bar{x})\Big).
\]
Since we assume that each $\chi_i$ is $\bar{y}$-positive it follows that 
if $\mbZ_n$ is the set of $\mcA \in \mbW_n$ such that
\[
\mcA \models \Big(\bigvee_{i=1}^m \neg \exists \bar{y} \chi_i(\bar{x}, \bar{y}) \Big) \leftrightarrow
\Big(\bigvee_{i=1}^t \chi'_i(\bar{x})\Big),
\]
then $\lim_{n\to\infty} \mbbP_n(\mbZ_n) = 1$.
Hence condition (c) of part~(2) of Assumption~\ref{assumptions on the basic logic} is satisfied and we have
verified that part~(2) of Assumption~\ref{assumptions on the basic logic} holds.

The final statement, regarding when conditions~(I) and~(II) of
Remark~\ref{remark about basic probability formulas} hold,
follows by the same kind of argument as in (the end of) the proof of 
Proposition~\ref{asymptotic elimination for strongly unbounded aggregations},
so we leave out the details.
\hfill $\square$

\begin{rem}\label{remark on kappa with an additional assumption} {\rm
Suppose that
$\sigma^+$ is a nonempty finite relational signature, $\sigma \subset \sigma^+$, and that
$\mbbG^+$ is a $PLA^*(\sigma^+)$-network such that $\mr{mp}(\mbbG^+) = \mr{mp}(\mbbG) + 1$.
Furthermore, suppose that $\sigma = \{R \in \sigma^+ : \mr{mp}(R) < \mr{mp}(\mbbG^+)\}$
and that $\mbbG$ is the subnetwork of $\mbbG^+$ which is induced by $\sigma$.
Then, for each $R \in \sigma^+ \setminus \sigma$, the corresponding formula $\theta_R$ (of $\mbbG^+$)
belongs to $PLA^*(\sigma)$.
Suppose that for every $R \in \sigma \setminus \tau$ and every subformula of $\theta_R$ of the form
\begin{equation}
F\big(\varphi_1(\bar{y}, \bar{z}), \ldots, \varphi_m(\bar{y}, \bar{z}) : \bar{z} : 
\chi_1(\bar{y}, \bar{z}), \ldots, \chi_m(\bar{y}, \bar{z})\big)
\end{equation}
(where $F$ is an aggregation function) it holds either that
\begin{enumerate}
\item for all $i = 1, \ldots, m$, $\varphi_i(\bar{y}, \bar{z})$ is finite valued and there is $\xi_i \in \mbbN$
such that for all $n$, $\mcA \in \mbW_n$, $\bar{a} \in (B_n)^{|\bar{y}|}$, and $\bar{b} \in (B_n)^{|\bar{z}|}$,
if $\mcA \models \chi_i(\bar{a}, \bar{b})$ then $\rng(\bar{b}) \subseteq C_{\xi_i}^{\mcB_n}(\bar{a})$, or that

\item for all $i = 1, \ldots, m$, $\chi_i(\bar{y}, \bar{z})$ is a $\bar{z}$-positive $(\sigma, \lambda_i)$-closure type
for some $\lambda_i \in \mbbN$ and
$F$ is continuous.
\end{enumerate}
By Proposition~\ref{asymptotic elimination under stronger conditions}, for
$R \in \sigma^+ \setminus \sigma$,
$\theta_R(\bar{x})$ is asymptotically equivalent to a $(\sigma, \xi_R)$-basic formula for some $\xi_R \in \mbbN$.
Let $\kappa_0 = \max\{\xi_R : R \in \sigma^+ \setminus \sigma\}$.
Then part~(1) of Assumption~\ref{induction hypothesis} holds if we replace 
$\sigma, \sigma', \kappa, \mbbG$, and $\mbbG'$ by
$\sigma^+, \sigma, \kappa_0, \mbbG^+$, and $\mbbG$, respectively.
This concludes the proof of the inductive step for part~(1)  of Assumption~\ref{induction hypothesis}
(provided that the assumption above on each $\theta_R$ holds).
}\end{rem}

\begin{theor}\label{main result about aggregations in general with extra assumption}
Suppose that Assumptions~\ref{properties of the base structures} and~\ref{relative frequency of tau-closure types}
are satisfied.
Suppose that for every $R \in \sigma \setminus \tau$ and every subformula of $\theta_R$ of the form
\begin{equation}\label{an aggregation to be eliminated in theorem with extra assumption}
F\big(\varphi_1(\bar{y}, \bar{z}), \ldots, \varphi_m(\bar{y}, \bar{z}) : \bar{z} : 
\chi_1(\bar{y}, \bar{z}), \ldots, \chi_m(\bar{y}, \bar{z})\big)
\end{equation} 
it holds either that
\begin{enumerate}
\item for all $i = 1, \ldots, m$, $\varphi_i(\bar{y}, \bar{z})$ is finite valued and there is $\xi_i \in \mbbN$
such that for all $n$, $\mcA \in \mbW_n$, $\bar{a} \in (B_n)^{|\bar{y}|}$, and $\bar{b} \in (B_n)^{|\bar{z}|}$,
if $\mcA \models \chi_i(\bar{a}, \bar{b})$ then $\rng(\bar{b}) \subseteq C_{\xi_i}^{\mcB_n}(\bar{a})$, or that

\item for all $i = 1, \ldots, m$, $\chi_i(\bar{y}, \bar{z})$ is a $\bar{z}$-positive $(\sigma, \lambda_i)$-closure type
for some $\lambda_i \in \mbbN$ and
$F$ is continuous.
\end{enumerate}
Let $\varphi(\bar{x}) \in PLA^*(\sigma)$ and suppose that for every subformula of $\varphi(\bar{x})$ of the 
form~(\ref{an aggregation to be eliminated in theorem with extra assumption})
condition~(1) or condition~(2) holds.
Then:\\
(i) $\varphi(\bar{x})$ is asymptotically equivalent to a $(\sigma, \xi)$-basic formula for some $\xi \in \mbbN$.\\
(ii) For the same $\xi \in \mbbN$ as in part~(i) the following holds:
For every complete $(\tau, \xi)$-closure type $p(\bar{x})$ there are $k \in \mbbN^+$, 
$c_1, \ldots, c_k \in [0, 1]$, and $\beta_1, \ldots, \beta_k \in [0, 1]$, depending only
on $p$ and $\varphi$, such that 
for every $\varepsilon > 0$ there is $n_0$ such that if $n \geq n_0$, 
and $\mcB_n \models p(\bar{a})$ then 
\begin{align*}
&\mbbP_n \Big( \{ \mcA \in \mbW_n : \mcA(\varphi(\bar{a})) \in \bigcup_{i = 1}^k [c_i - \varepsilon, c_i + \varepsilon] \Big)
\geq 1 - \varepsilon \text{ and, for all $i = 1, \ldots, k$,} \\
&\mbbP_n\big(\{\mcA \in \mbW_n : \mcA(\varphi(\bar{a})) \in [c_i - \varepsilon, c_i + \varepsilon] \}\big) \in 
[\beta_i - \varepsilon, \beta_i + \varepsilon].
\end{align*}
\end{theor}

\noindent
{\bf Proof.}
(i) We just modify the proof of Theorem~\ref{main result about strongly unbounded aggregations}.
If $\sigma = \tau$ then the claims of 
Assumption~\ref{induction hypothesis} (with the choice $\kappa = \kappa' = 0$)
follows  from 
Lemma~\ref{base case},
as pointed out in 
Remark~\ref{remark on induction hypothesis and base case}.
From Remarks~\ref{the new kappa'}, \ref{the new kappa' for balance} and \ref{remark on kappa with an additional assumption}
it follows that if Assumption~\ref{induction hypothesis} holds for every $PLA^*(\sigma)$-network with
mp-rank $\rho$ for some choice of $\kappa$ and $\kappa'$, then Assumption~\ref{induction hypothesis} also holds for every $PLA^*(\sigma)$-network
with mp-rank $\rho + 1$ for some (possibly other) choice of $\kappa$ and $\kappa'$. 
Thus, Assumption~\ref{induction hypothesis} holds for every finite relational signature $\sigma$
and every $PLA^*(\sigma)$-network subject to conditions~(1) and~(2) of the theorem.
Proposition~\ref{asymptotic elimination under stronger conditions}
now implies that 
if $\varphi(\bar{x}) \in PLA^*(\sigma)$ and for every subformula of $\varphi$ of the 
form~(\ref{an aggregation to be eliminated in theorem with extra assumption}) either condition~(1) or~(2) holds,
then $\varphi(\bar{x})$ is asymptotically equivalent to a $(\sigma, \xi)$-basic formula for some $\xi \in \mbbN$.

Part~(ii) is proved just like part~(ii) of 
Theorem~\ref{main result about strongly unbounded aggregations}
except that we consider complete $(\sigma, \xi)$-closure types here. The details are left for the reader.
\hfill $\square$

\begin{cor}\label{corollary to main result about aggregations in general with extra assumption}
Suppose that Assumptions~\ref{properties of the base structures} and~\ref{relative frequency of tau-closure types}
are satisfied.
Suppose that for every $R \in \sigma \setminus \tau$, there is $\lambda_R \in \mbbN$
such that $\theta_R$ is a $(\mr{par}(R) \cup \tau, \lambda_R)$-closure type.
Let $\varphi(\bar{x}) \in PLA^*(\sigma)$ be such that if
\[
F\big(\varphi_1(\bar{y}, \bar{z}), \ldots, \varphi_m(\bar{y}, \bar{z}) : \bar{z} : 
\chi_1(\bar{y}, \bar{z}), \ldots, \chi_m(\bar{y}, \bar{z})\big)
\]
is a subformula of $\varphi(\bar{x})$ then either
\begin{enumerate}
\item for all $i = 1, \ldots, m$, $\varphi_i(\bar{y}, \bar{z})$ is finite valued and there is $\xi_i \in \mbbN$
such that for all $n$, $\mcA \in \mbW_n$, $\bar{a} \in (B_n)^{|\bar{y}|}$, and $\bar{b} \in (B_n)^{|\bar{z}|}$,
if $\mcA \models \chi_i(\bar{a}, \bar{b})$ then $\rng(\bar{b}) \subseteq C_{\xi_i}^{\mcB_n}(\bar{a})$, or

\item for all $i = 1, \ldots, m$, $\chi_i(\bar{y}, \bar{z})$ is a $\bar{z}$-positive $(\sigma, \lambda_i)$-closure type
for some $\lambda_i \in \mbbN$ and
$F$ is {\rm admissible}.
\end{enumerate}
Then the conclusions~(i) and~(ii) of 
Theorem~\ref{main result about aggregations in general with extra assumption}
hold.
\end{cor}

\noindent
{\bf Proof.}
The argument follows the pattern of the proof of 
Corollary~\ref{corollary to main result about strongly unbounded aggregations}
but uses Theorem~\ref{main result about aggregations in general with extra assumption}
and Proposition~\ref{asymptotic elimination under stronger conditions}
and that under the present assumptions $(p, p_\tau)$ is eventually constant whenever 
$p(\bar{x})$ is a complete $(\sigma, \lambda)$-closure type and $p_\tau$ a complete $(\tau, \gamma)$-closure type
for some $\lambda, \gamma \in \mbbN$.
The details are left for the reader.
\hfill $\square$

\begin{exam}\label{example of results with additional assumption}{\rm
Suppose that $\sigma = \tau \cup \{E\}$ where $E \notin \tau$ is a binary relation symbol.
Let $\mbbG$ be a $PLA^*(\sigma)$-network based on $\tau$ such that, for some $\lambda \in \mbbN$,
$\theta_E(x_1, x_2)$ has the form $\bigwedge_{i=1}^s (\chi_i(x_1, x_2) \to c_i)$ where $\chi_i(x_1, x_2)$,
$i = 1, \ldots, s$, enumerates all, up to equivalence, $(\tau, \lambda)$-closure types 
in the variables $\bar{x} = (x_1, x_2)$ and $c_i \in [0, 1]$.
Then $\mbbG$ ``expresses'' that if $\chi_i(a_1, a_2)$ holds (i.e. if the $\tau$-structure of the
$\lambda$-closure of $(a_1, a_2)$ is a described by $\chi_i(x_1, x_2)$), then the probability that
$R(a_1, a_2)$ holds is $c_i$.
Without loss of generality assume that, for some $t \leq s$,  $\chi_1(x_1, x_2), \ldots, \chi_t(x_1, x_2)$
enumerates all $\chi_i(x_1, x_2)$ such that $c_i > 0$.
Then all aggregations of $\theta_E$ satisfy condition~(1) of 
Theorem~\ref{main result about aggregations in general with extra assumption}
(if $E$ takes the role of $R$ in that theorem).

Recall Example~\ref{example of page rank}
where it was shown that, for every $k \in \mbbN$, the $k^{th}$ stage of the PageRank
can be expressed by a $PLA^*(\sigma)$-formula $\varphi_k(x)$.
However, we cannot apply any of our results directly to show that the $\varphi_k(x)$ from 
Example~\ref{example of page rank} is asymptotically equivalent to a $(\sigma, \lambda)$-basic formula
for some $\lambda \in \mbbN$.
This is because the aggregations used in the formula $\varphi_k(x)$ from
Example~\ref{example of page rank}
do not have the properties required by conditions~(1) and (2) of
Theorem~\ref{main result about aggregations in general with extra assumption}
or Corollary~\ref{corollary to main result about aggregations in general with extra assumption}.

However, we can construct a variant of the of the formula $\varphi_k(x)$ from 
Example~\ref{example of page rank}
which captures the idea of the PageRank and which is asymptotically equivalent to a $(\sigma, \gamma)$-basic
formula for some $\gamma \in \mbbN$.
For every $k \in \mbbN^+$ let $\mr{Av}_k : [0, 1]^k \to [0, 1]$ be the continuous connective defined by
letting $\mr{Av}_k(r_1, \ldots, r_k)$ be the average of $r_1, \ldots, r_k$.
For each $i = 1, \ldots, s$ let $\varphi_0^i(x)$ be the formula
$\mr{length}^{-1}(x=x : y : \chi_i(x, y))$.
Then, for every finite $\sigma$-structure $\mcA$ and $a \in A$, $\mcA(\varphi_0^i(a)) = (m_i)^{-1}$
where $m_i$ is the number of $b \in A$ such that $\mcA \models \chi_i(a, b)$, if such $b$ exists, and otherwise
$\mcA(\varphi_0^i(a)) = 0$.
Let $\varphi_0(x)$ 
be the formula
$\mr{Av}_t(\varphi_0^1(x), \ldots, \varphi_0^t(x))$.
For each $i = 1, \ldots, s$, let $\psi_i(y)$ be the formula
$\mr{length}^{-1}(y=y : z : \chi_i(y, z) \wedge E(y, z))$ and let 
$\psi(y)$ be the formula
$\mr{Av}_t(\psi_1(y), \ldots, \psi_t(y))$.

Now let $\varphi_{k+1}^i(x)$ be the formula
$\mr{tsum}\big(x=x \wedge (\varphi_k(y) \cdot \psi(y)) : y : \chi_i(y, x) \wedge E(y, x)\big)$
and let
$\varphi_{k+1}(x)$ be the formula
$\mr{Av}_t\big(\varphi_{k+1}^1(x), \ldots, \varphi_{k+1}^t(x)\big)$
(where the aggregation function tsum was defined in Example~\ref{example of page rank}.)

Since the aggregation functions $\mr{length}^{-1}$ and tsum are continuous
and, for $i = 1, \ldots, t$, $\chi_i(x_1, x_2) \wedge R(x_1, x_2)$ is $x_2$-positive
it follows 
that, for every $k \in \mbbN$, the formula $\varphi_k(x)$ satisfies the conditions of the formula called $\varphi(\bar{x})$
in Theorem~\ref{main result about aggregations in general with extra assumption}.
Consequently 
(by Theorem~\ref{main result about aggregations in general with extra assumption}), 
for each $k$, $\varphi_k(x)$ is asymptotically equivalent to a $(\sigma, \lambda_k)$-basic formula for some $\lambda_k  \in  \mbbN$.
Moreover, the distribution of the truth values of $\varphi_k(x)$ tends to a limit (as $n\to\infty$)
in the sense of Theorem~\ref{main result about aggregations in general with extra assumption}.

Above we used the average connective $\mr{Av}_t : [0, 1]^t \to [0, 1]$ when defining $\varphi_k(x)$.
But if we have information about the ratio $|\chi_i(\bar{a}, \mcB_n)| / |\chi_j(\bar{a}, \mcB_n)|$ 
as $n$ tends to infinity, for all (or some) $i \neq j$,
then it may be more suitable to use a {\em weighted} average to define an estimate of the
$k^{th}$ approximation of the PageRank.
}\end{exam}

\noindent
{\bf Acknowledgement.}
The author was partially supported by the Swedish Research Council, grant 2023-05238\_VR.

\end{document}